\documentclass[12pt,oneside]{book}
\usepackage{geometry}                		
\usepackage{graphicx}									
\usepackage{amssymb}\vspace{0.25cm}
\usepackage{mathtools}
\usepackage{setspace}
\usepackage{tikz-cd}
\usepackage[mathscr]{euscript}
\newcommand{\abs}[1]{\left\vert#1\right\vert}
\usepackage[normalem]{ulem}
\usepackage{manfnt}
\usepackage{pgfplots}
\usepackage{pifont}
\usepackage[safe]{tipa}
\usepackage{marvosym}
\usepackage{scalerel,stackengine}
\usepackage{titlesec}
\usepackage{makeidx}
\usepackage{centernot}
\usepackage{xspace}
\usepackage[noauto]{chappg}
\usepackage[OT2,T1]{fontenc}
\usepackage{hyperref}
\usepackage{amsmath}
\usepackage[amsmath,thmmarks]{ntheorem}
\usepackage{tabularx}
\usepackage{yfonts}
\usepackage{epstopdf}
\usepackage{pdfpages}
\newcolumntype{Y}{>{\raggedright\arraybackslash}X}

\newcommand\N{%
\mathbb{N}
}

\newcommand\R{%
\mathbb{R}
}

\newcommand\XX{%
\mathbb{X}
}

\newcommand\sD{%
\mathcal{D}
}

\newcommand\sF{%
\mathcal{F}
}

\newcommand\sM{%
\mathcal{M}
}

\newcommand\sP{%
\mathcal{P}
}

\newcommand\sV{%
\mathcal{V}
}

\newcommand\ap{%
\text{ap}
}

\newcommand\REG{%
\text{REG}
}

\newcommand\loc{%
\text{loc}
}
\newcommand\locx{%
\text{$\ell$oc}
}

\newcommand\td{%
\text{d}
}

\newcommand\tc{%
\text{c}
}

\newcommand\iI{
\text{I}\hspace{0.05cm}
}

\newcommand\Lm{%
\text{L}
}
\newcommand\Lp{
\text{L}
}

\newcommand\subsetx{%
\ \subset \ 
}

\newcommand\capx{%
\hsx \cap \hsx 
}

\newcommand\CantorSet{%
\text{P}
}
\newcommand\LebesgueFunction{%
\text{L}
}

\newcommand\tD{%
\text{D}
}

\newcommand\tS{%
\text{S}
}

\newcommand\tT{%
\text{T}
}

\def\thereforex{\boldsymbol{\text{ }
\leavevmode
\lower0.4ex\hbox{\textbullet}
\kern-.9em\raise1.1ex\hbox{\textbullet}
\kern-0.9em\lower0.4ex\hbox{\textbullet}
\hspace{0.1cm}\thinspace\text{ }}}
\def\thereforez{\boldsymbol{\text{ }
\leavevmode
\lower0.4ex\hbox{$\circ$}
\kern-.9em\raise1.1ex\hbox{$\circ$}
\kern-0.9em\lower0.4ex\hbox{$\circ$}
\hspace{0.1cm}\thinspace\text{ }}}
\newcommand\ra{%
\rightarrow
}
\newcommand\lra{%
\longrightarrow
}

\newcommand\ds{%
\displaystyle
}

\newcommand\sptx{%
\text{spt}\hspace{0.05cm}
}

\newcommand\diam{%
\text{diam}\hspace{0.05cm}
}
\newcommand\osc{%
\text{osc}\hspace{0.05cm}
}

\newcommand\sgn{%
\text{sgn}\hspace{0.05cm}
}

\newcommand\un[1]{%
\underline{#1}\xspace
}

\newcommand\restr[2]{%
{#1}|{#2}
}

\newcommand\hsx{%
\hspace{0.05cm}
}
\newcommand\hsy{%
\hspace{0.03cm}
}

\newcommand{\norm}[1]{\left\lVert #1 \right\rVert}

\newcommand{\normx}[1]{\big|\hspace{-.05cm}\big| #1 \big|\hspace{-.05cm}\big|}

\newcommand\hps{%
\hspace{0.7cm}
}

\newcommand\BV{%
 \text{BV}
}
\newcommand\BVC{%
\text{BVC}
}
\newcommand\BVL{%
\text{BVL}
}
\newcommand\NBV{%
\text{NBV}
}
\newcommand\AC{%
\text{AC}
}
\newcommand\AT{%
\text{AT}
}

\newcommand\ov[1]{%
\overline{#1}
}

\newcommand\chisubt{\chi\raisebox{-.1cm}{$\scaleto{_t}{5.5pt}$}}
\newcommand\chisubtsubone{\chi\raisebox{-.1cm}{$\scaleto{_{t_{1}}}{6.5pt}$}}
\newcommand\chisubtsubtwo{\chi\raisebox{-.1cm}{$\scaleto{_{t_{2}}}{6.5pt}$}}

\newcommand\chisubn{\chi\raisebox{-.1cm}{$\scaleto{_n}{4.5pt}$}}
\newcommand\chisubni{\chi\raisebox{-.1cm}{$\scaleto{_{n \hsy i}}{6.5pt}$}}
\newcommand\chisubnplusone{\chi\raisebox{-.1cm}{$\scaleto{_{n + 1}}{6.5pt}$}}

\newcommand\chisubI{\chi\raisebox{-.1cm}{$\scaleto{_I}{4.5pt}$}}

\definecolor{ultramarine}{RGB}{0, 32, 96}

\definecolor{darkcerulean}{rgb}{0.3, 0.27, 0.49}

\definecolor{forestgreen}{rgb}{0.0, 0.27, 0.13}
\definecolor{forestgreenweb}{rgb}{0.13, 0.55, 0.13}
\definecolor{deepjunglegreen}{rgb}{0.0, 0.29, 0.29}

\definecolor{midnightblue}{rgb}{0.1, 0.1, 0.44}
\definecolor{midnightgreen}{rgb}{0.0, 0.29, 0.33}
\definecolor{myrtle}{rgb}{0.13, 0.26, 0.12}
\definecolor{darkviolet}{rgb}{0.58, 0.0, 0.83}
\definecolor{darkgreen}{rgb}{0.0, 0.2, 0.13}
\definecolor{officegreen}{rgb}{0.0, 0.5, 0.0}

\definecolor{harvardcrimson}{rgb}{0.79, 0.0, 0.09}
\definecolor{hollywoodcerise}{rgb}{0.96, 0.0, 0.63}
\definecolor{debianred}{rgb}{0.84, 0.04, 0.33}
\definecolor{darkturquoise}{rgb}{0.0, 0.81, 0.82}

\definecolor{darktangernine}{rgb}{1.0, 0.66, 0.07}
\definecolor{aureolin}{rgb}{0.99, 0.93, 0.0}
\definecolor{canaryyellow}{rgb}{1.0, 0.94, 0.0}
\definecolor{amber}{rgb}{1.0, 0.75, 0.0}

\definecolor{urobilin}{rgb}{0.88, 0.68, 0.13}
\definecolor{uscgold}{rgb}{1.0, 0.8, 0.0}

\usepackage{scalerel,stackengine}
\stackMath
\newcommand\reallywidehat[1]{%
\savestack{\tmpbox}{\stretchto{%
  \scaleto{%
    \scalerel*[\widthof{\ensuremath{#1}}]{\kern-.6pt\bigwedge\kern-.6pt}%
    {\rule[-\textheight/2]{1ex}{\textheight}}
  }{\textheight}%
}{0.5ex}}%
\stackon[1pt]{#1}{\tmpbox}%
}
\makeatletter
\DeclareRobustCommand\widecheck[1]{{\mathpalette\@widecheck{#1}}}
\def\@widecheck#1#2{%
    \setbox\z@\hbox{\m@th$#1#2$}%
    \setbox\tw@\hbox{\m@th$#1%
       \widehat{%
          \vrule\@width\z@\@height\ht\z@
          \vrule\@height\z@\@width\wd\z@}$}%
    \dp\tw@-\ht\z@
    \@tempdima\ht\z@ \advance\@tempdima2\ht\tw@ \divide\@tempdima\thr@@
    \setbox\tw@\hbox{%
       \raise\@tempdima\hbox{\scalebox{1}[-1]{\lower\@tempdima\box
\tw@}}}%
    {\ooalign{\box\tw@ \cr \box\z@}}}
\makeatother

\newtheoremstyle{xx}
  {4pt}
  {0pt}
  {\upshape}
  {\bfseries}
  {}
  { }
  {}
  
\makeatletter 
 \newtheoremstyle{myu}%
  {\upshape\item[ \indent\indent\bf\underline{\theorem@headerfont ##2:}]}%
\makeatother
\makeatletter 
 \newtheoremstyle{myn}%
  {\item[\hskip\labelsep \ \bf ##1 \theorem@headerfont ##2.]}%
\makeatother
\theoremstyle{myn}
\newtheorem{theoremn}{Theorem} 
\theoremstyle{myu}
{\upshape}
\newtheorem{x}[theoremn]{}

 \newtheoremstyle{mr}%
  {\upshape\item[ \indent{\theorem@headerfont ##2. \hspace{.2cm}}]}%
\makeatother
\theoremstyle{mr}
{\upshape}
\newtheorem{rf}[theoremn]{}

\title{\textbf{Analysis 101:\\
Functions of a Single Variable}}
\author{Garth Warner\\
Department of Mathematics\\
University of Washington}
\date{}	

\titleformat{\chapter}[display]
{\normalfont\filcenter\huge\bfseries}{}{0pt}{\large}

\titleformat{\chapter}[display]
{\normalfont\filcenter\huge\bfseries}{}{0pt}{\large}

\usepackage[OT2,OT1]{fontenc} 
\newcommand\cyr
{
\renewcommand\rmdefault{wncyr} 
\renewcommand\sfdefault{wncyss} 
\renewcommand\encodingdefault{OT2} 
\normalfont
\selectfont
}

\DeclareTextFontCommand{\textcyr}{\cyr}

\makeindex 
\begin{document}

\maketitle                              

\titlespacing*{\chapter}{0pt}{-50pt}{40pt}
\setlength{\parskip}{0.1em}
\pagenumbering{bychapter}
\setcounter{chapter}{0}
\pagenumbering{bychapter}

\begingroup
\fontsize{11pt}{11pt}\selectfont

\[
\textbf{ABSTRACT}
\]
\\

These notes are a chapter in Real Analysis.  While primarily standard, the reader will find a discussion of certain topics that are ordinarily not covered in the usual accounts.
For example, the notion of bounded variation in the sense of Cesari is introduced.
This is Volume 1 of 4, to be followed by Curves and Length, Functions of Several Variables, and Surfaces and Area.
\\[2cm]

\[
\textbf{ACKNOWLEDGEMENT}
\]

Many thanks to David Clark for his rendering the original transcript into AMS-LaTeX.  
Both of us also thank Judith Clare for her meticulous proofreading.
\newpage

\[
\textbf{FUNCTIONS OF A SINGLE VARIABLE}
\]
\\


\hspace{2.1cm} \ \S0. \quad RADON MEASURES%
\\[-.26cm]

\hspace{2.1cm} \ \S1. \quad VARIATION OF A FUNCTION%
\\[-.26cm]

\hspace{2.1cm} \ \S2. \quad LIMIT AND OSCILLATION%
\\[-.26cm]

\hspace{2.1cm} \ \S3. \quad FACTS AND EXAMPLES %
\\[-.26cm]

\hspace{2.1cm} \ \S4. \quad PROPERTIES %
\\[-.26cm]

\hspace{2.1cm} \ \S5. \quad REGULATED FUNCTIONS %
\\[-.26cm]

\hspace{2.1cm} \ \S6. \quad POSITIVE AND NEGATIVE %
\\[-.26cm]

\hspace{2.1cm} \ \S7. \quad CONTINUITY %
\\[-.26cm]

\hspace{2.1cm} \ \S8. \quad ABSOLUTE CONTINUITY I%
\\[-.26cm]

\hspace{2.1cm} \ \S9. \quad DINI DERIVATIVES %
\\[-.26cm]

\hspace{2.1cm} \S10. \quad DIFFERENTIATION %
\\[-.26cm]

\hspace{2.1cm} \S11. \quad ESTIMATE OF THE IMAGE 
\\[-.26cm]

\hspace{2.1cm} \S12. \quad ABSOLUTE CONTINUITY II
\\[-.26cm]

\hspace{2.1cm} \S13. \quad MULTIPLICITIES 
\\[-.26cm]

\hspace{2.1cm} \S14. \quad LOWER SEMICONTINUITY 
\\[-.26cm]

\hspace{2.1cm} \S15. \quad FUNCTIONAL ANALYSIS 
\\[-.26cm]

\hspace{2.1cm} \S16. \quad DUALITY 
\\[-.26cm]

\hspace{2.1cm} \S17. \quad INTEGRAL MEANS 
\\[-.26cm]

\hspace{2.1cm} \S18. \quad ESSENTIAL VARIATION 
\\[-.26cm]

\hspace{2.1cm} \S19. \quad BVC 
\\[-.26cm]

\hspace{2.1cm} \S20. \quad ABSOLUTE CONTINUITY III
\\

\hspace{3.35cm}   REFERENCES
\\

\[
\]



\endgroup 

\chapter{
$\boldsymbol{\S}$\textbf{0}.\quad  RADON MEASURES}
\setlength\parindent{2em}
\setcounter{theoremn}{0}
\renewcommand{\thepage}{\S0-\arabic{page}}

\indent
\indent 
Let $X$ be a locally compact Hausdorff space.
\\

\begin{x}{\small\bf NOTATION} \ 
$C(X)$ is the set of real valued continuous functions on $X$ and
$BC(X)$ is the set of bounded real valued continuous functions on $X$.
\\[-.55cm]
\end{x}

\begin{x}{\small\bf DEFINITION} \ 
Given $f \in C(X)$, its \un{support}, denoted \sptx $(f)$, is the smallest closed subset of $X$ outside of which 
$f$ vanishes, 
i.e., the closure of $\{x: f(x) \neq 0\}$, and $f$ is said to be 
\un{compactly supported} provided \sptx $(f)$ is compact.
\\[-.55cm]
\end{x}

\begin{x}{\small\bf NOTATION} \ 
$C_c(X)$ is the subset of $C(X)$ whose elements are compactly supported.
\\[-.55cm]
\end{x}

\begin{x}{\small\bf DEFINITION} \ 
A function $f \in C(X)$ is said to \un{vanish at infinity} if $\forall \varepsilon > 0$, the set
\[
\{x: \abs{f(x)} 
\ \geq \ 
\varepsilon\}
\]
is compact.
\\[-.55cm]
\end{x}

\begin{x}{\small\bf NOTATION} \ 
$C_0(X)$ is the subset of $C(X)$ whose elements vanish at infinity.
\\[-.55cm]
\end{x}

\begin{x}{\small\bf \un{N.B.}} \ 
$C_c(X) \subset C_0(X) \subset C_c(X)$.
\\[-.55cm]
\end{x}

\begin{x}{\small\bf LEMMA} \ 
$C_0(X)$ is the closure of $C_c(X)$ in the uniform metric:
\[
\td (f, g) 
\ = \ 
\norm{f - g}_\infty.
\]
\end{x}

\begin{x}{\small\bf DEFINITION} \ 
A linear functional $\iI : C_c(X) \ra \R$ is 
\un{positive} if 
\[
f \geq 0 
\implies 
\iI \hsx (f) \geq 0.
\]
\end{x}


\begin{x}{\small\bf LEMMA} \ 
If $\iI$ is a positive linear functional on $C_c(X)$, then for each compact set $K \subset X$ 
there is a constant $C_K \geq 0$ such that
\[
\abs{\iI (f)}
\ \leq \ 
C_K \hsx \norm{f}_\infty
\]
for all $f \in C_c(X)$ such that $\sptx(f) \subset K$.
\\[-.55cm]
\end{x}

\begin{x}{\small\bf DEFINITION} \ 
A \un{Radon measure} on $X$ is a Borel measure $\mu$ that is 
finite on compact sets, 
outer regular on Borel sets, 
and inner regular on open sets.
\\[-.55cm]
\end{x}

\begin{x}{\small\bf EXAMPLE} \ 
Take $X = \R^n$ $-$then the restriction of Lebesgue measure $\lambda$ to the Borel sets in $X$ is a Radon measure. 
\\[-.55cm]
\end{x}

Every Radon measure $\mu$ on $X$ gives rise to a positive linear functional on $C_c(X)$, viz. the assignment 
\[
f \ra \int\limits_X \ f \ \td \mu.
\]
And all such arise in this fashion. 
\\

\begin{x}{\small\bf RIESZ REPRESENTATION THEOREM} \ 
If $\iI$ is a positive linear functional on $C_c(X)$, then there exists a unique Radon measure $\mu$ on $X$ such that
\[
\iI(f) 
\ = \ 
\int\limits_\R \ f \ \td \mu
\]
for all $f \in C_c(X)$.
\\[-.55cm]
\end{x}

\begin{x}{\small\bf EXAMPLE} \ 
Take $X = \R$ and define $\iI$ by the rule
\[
\iI(f) 
\ = \ 
\int\limits_X \ f \ \td x
\hspace{0.7cm} \text{(Riemann integral).}
\]
Then the Radon measure in this setup per the RRT is the restriction of Lebesgue measure $\lambda$ 
on the line to the Borel sets.
\\[-.55cm]
\end{x}


\begin{x}{\small\bf RAPPEL} \ 
$C_c(X)$ is a complete topological vector space when equipped with the 
\un{inductive topology}, i.e., the topology of uniform convergence on compact sets. 
\\[-.55cm]
\end{x}

\begin{x}{\small\bf DEFINITION} \ 
A \un{distribution of order 0} is a continuous linear functional 
\[
\tT : C_c(X) \ra \R.
\]
\end{x}

\begin{x}{\small\bf LEMMA} \ 
A linear functional $\tT : C_c(X) \ra \R$ is a distribution of order 0 iff for each compact set $K \subset X$ 
there is a constant $C_K > 0$ such that 
\[
\abs{\tT(f)} 
\ \leq \ 
C_K \hsx \norm{f}_\infty
\]
for all $f \in C_c(X)$ such that $\sptx (f) \subset K$.
\\[-.55cm]
\end{x}

Therefore a positive linear function $\iI : C_c(X) \ra \R$ is a distribution of order 0, 
hence is continuous in the inductive topology.
\\[-.25cm]

Denote the set of distributions of order 0 by the symbol $\sD^0$.
\\[-.25cm]

\begin{x}{\small\bf LEMMA} \ 
$\sD^0$ is a vector lattice.
\\[-.55cm]
\end{x}

If $\tT \in \sD^0$, then its Jordan decomposition is given by
\[
\tT
\ = \ 
\tT^+ - \tT^-,
\]
where
\[
\begin{cases}
\ \tT^+ \hsx (f) \ = \hspace{0.5cm} \sup\limits_{0 \leq g \leq f} \ \tT(g)\\
\ \tT^- \hsx (f) \ = \ -\inf\limits_{0 \leq g \leq f} \ \tT(g)
\end{cases}
.
\]
Here $\tT^+$, $\tT^- \in \sD^0$ are positive linear functionals and 
\[
\tT
\ = \ 
\tT^+ - \tT^-.
\]
Therefore
\[
\begin{cases}
\ \tT^+ \longleftrightarrow \mu^+\\
\ \tT^- \longleftrightarrow \mu^-
\end{cases}
\hspace{0.7cm} \text{(Radon),}
\]
so $\forall \ f \in C_c(X)$, 
\[
\tT \hsx (f)
\ = \ 
\int\limits_X \ f \ \td \mu^+ - \int\limits_X \ f \ \td \mu^-
\]
and
\[
\abs{\tT} \hsx (f)
\ = \ 
\int\limits_X \ f \ \td \hsy (\mu^+ + \mu^-).
\]

\begin{x}{\small\bf \un{N.B.}} \ 
Both $\mu^+$ and $\mu^-$ might have infinite measure, thus in general their difference is not defined.
\\[-.55cm]
\end{x}

\begin{x}{\small\bf REMARK} \ 
As we have seen, the positive linear functionals on $C_c(X)$ can be identified with the Radon measures.  
Bearing in mind that $C_0(X)$ is the uniform closure of $C_c(X)$, the positive linear functionals on $C_0(X)$ 
can be identified with the finite Radon measures. 
\end{x}

\[
* \ * \ * \ * \ * \ * \ * \ * \ * \ * \ 
\]

Let $X$ be a compact Hausdorff space.
\\[-.25cm]

\begin{x}{\small\bf \un{N.B.}} \ 
It is clear that in this situation $C_c(X) = C(X)$.
\\[-.55cm]
\end{x}

Equip $C(X)$ with the uniform norm:
\[
\norm{f}_\infty 
\ = \ 
\sup\limits_X \ \abs{f}.
\]
Then the pair $(C(X), \norm{\hsx \cdot \hsx}_\infty)$ is a Banach space.  
Let $C(X)^*$ be the dual space of $C(X)$, 
i.e., the linear space of all continuous linear functionals $\Lambda$ on $C(X)$ $-$then the prescription
\[
\norm{\Lambda}^* 
\ = \ 
\inf \ \{M \geq 0 : \abs{\Lambda \hsy (f)} \leq M \hsx \norm{f}_\infty 
\hspace{0.25cm} (f \in C(X))\}
\]
is a norm on $C(X)^*$ under which the pair $(C(X)^*, \norm{\hsx \cdot \hsx})$ is a Banach space.

\begin{x}{\small\bf \un{N.B.}} \ 
$\forall \ f \in C(X)$, $\forall \ \Lambda \in C(X)^*$, 
\[
\abs{\Lambda \hsy (f)} 
\ \leq \ 
\norm{\Lambda}^*  \hsx \norm{f}_\infty.
\]
\end{x}

\begin{x}{\small\bf RAPPEL} \ 
A \un{signed Radon measure} is a signed Borel measure $\mu$ whose positive variation $\mu^+$ is Radon 
and whose negative variation $\mu^-$ is Radon.
\\[-.5cm]

[Note: \ 
As usual, $\mu = \mu^+ - \mu^-$ is the Jordan decomposition of $\mu$ and its total variation, 
denoted $\abs{\mu}$, is by definition $\abs{\mu} \mu^+ + \mu^-$.  
In addition, $\mu$ is finite if $\abs{\mu}$ is finite, i.e., if $\abs{\mu} \hsy (X) < +\infty$.]
\\[-.55cm]
\end{x}

\begin{x}{\small\bf RIESZ REPRESENTATION THEOREM} \ 
Given a $\Lambda \in C(X)^*$, there exists a unique finite signed Radon measure $\mu$ such that 
$\forall \ f \in C(X)$, 
\[
\Lambda \hsy (f) 
\ = \ 
\int\limits_X \ f \ \td \mu.
\]
And
\[
\norm{\Lambda}^* 
\ = \ 
\abs{\mu} \hsy (X).
\]
\end{x}

\begin{x}{\small\bf NOTATION} \ 
$\sM(X)$ is the set of finite signed Radon measures on $X$.
\\[-.55cm]
\end{x}

\begin{x}{\small\bf LEMMA} \ 
$\sM(X)$ is a vector space of $\R$.
\\[-.55cm]
\end{x}

\begin{x}{\small\bf NOTATION} \ 
Given $\mu \in \sM(X)$, put
\[
\norm{\mu}_{\sM(X)}
\ = \ 
\abs{\mu}  (X).
\]
\end{x}


\begin{x}{\small\bf LEMMA} \ 
$\norm{\hsx \cdot \hsx}_{\sM(X)}$ is a norm on $\sM(X)$ under which the pair 
$(\sM(X), \norm{\hsx \cdot \hsx}_{\sM(X)})$ is a Banach space. 
\\[-.55cm]
\end{x}

\begin{x}{\small\bf THEOREM} \ 
Define an arrow
\[
\Lambda : \sM(X) \ra C(X)^*
\]
by the rule
\[
\Lambda\hsy (\mu) \hsy (f)
\ = \ 
\int\limits_X \ f \td \mu.
\]
Then $\Lambda$ is an isometric isomorphism.
\\[-.15cm]

[E.g.: 
\begin{align*}
\abs{\Lambda\hsy (\mu) \hsy (f)}\ 
&= \
\bigg|\int\limits_X \ f \td \mu \hsy\bigg|\ 
\\[15pt]
&\leq\ 
\int\limits_X \ \abs{f} \td \abs{\mu} \ 
\\[15pt]
&\leq\ 
\norm{f}_\infty \norm{\mu} (X) \ 
\\[15pt]
&=\ 
\norm{f}_\infty \abs{\mu}_{\sM(X)}.
\end{align*}
Therefore
\[
\Lambda\hsy (\mu) \in C(X)^*.]
\]
\end{x}

\[
* \ * \ * \ * \ * \ * \ * \ * \ * \ * \ 
\]

If $X$ is not compact, then the story for $C_0(X)$ is the same as that for $C(X)$ when $X$ is compact.  
Without stopping to spell it all out, once again the bounded linear functionals are in a one-to-one 
correspondence with the finite signed Radon measures and 
\[
\norm{\Lambda}^*
\ = \ 
\abs{\mu} (X).
\]

\chapter{
$\boldsymbol{\S}$\textbf{1}.\quad VARIATION OF A FUNCTION}
\setlength\parindent{2em}
\setcounter{theoremn}{0}
\renewcommand{\thepage}{\S1-\arabic{page}}

Let $[a,b] \subset \R$ be a closed interval $(a < b, \ -\infty < a < b < +\infty)$.
\\

\begin{x}{\small\bf DEFINITION} \ 
A \un{partition} of $[a,b]$ is a finite set $P = \{x_0, \ldots, x_n\} \subset [a,b]$, where
\[
a 
\ = \ 
x_0 
\ < \ 
x_1 
\ < \ 
\cdots 
\ < \ 
x_n 
\ = \ 
b.
\]
\end{x}

\begin{x}{\small\bf NOTATION} \ 
The set of all partitions of $[a,b]$ is denoted by $\sP [a,b]$.
\end{x}

\begin{x}{\small\bf EXAMPLE} \ 
\[
\{a, b\} \in \sP [a,b].
\]
\end{x}

Let $(\XX, \td)$ be a metric space and let $f:[a,b] \ra \XX$ be a function.
\\

\begin{x}{\small\bf DEFINITION} \ 
Given a partition $P \in \sP [a,b]$, put
\[
\bigvee\limits_a^b \ (f; P) 
\ = \ 
\sum\limits_{i = 1}^n \ \td(f(x_i), \hsy f(x_{i - 1})), 
\]
the \un{variation} of $f$ in $P$.
\end{x}

\begin{x}{\small\bf NOTATION} \ 
Put
\[
\tT_f [a,b] 
\ = \ 
\sup\limits_{P \in \sP [a,b]} \ \bigvee\limits_a^b \ (f; P),
\]
the \un{total variation} of $f$ in $[a,b]$.
\end{x}

\begin{x}{\small\bf \un{N.B.}} \ 
Here, $(\XX, \td)$ is implicit $\ldots$ \hsx .
\end{x}

One can then develop the basics at this level of generality but we shall
instead specialize immediately and take 
\[
\XX 
\ = \ 
\R, 
\quad 
\td (x, y) 
\ = \ 
\abs{x - y},
\]
thus now $f:[a,b] \ra \R$ .  
Later on, we shall deal with the situation when the domain $[a,b]$ is replaced by the open interval $]a,b[$ 
(or in principle, by any nonempty open set $\Omega \subset \R$ 
(recall that such an $\Omega$ can be written as an at most countable union of pairwise disjoint open intervals), 
e.g., $\Omega = \R$).  
As for the range, we shall stick with $\R$ for the time being but will eventually consider matters when $\R$ is replaced by 
$\R^M$ $(M = 1, 2, \ldots, )$ (curve theory).

\chapter{
$\boldsymbol{\S}$\textbf{2}.\quad  LIMITS AND OSCILLATION}
\setlength\parindent{2em}
\setcounter{theoremn}{0}
\renewcommand{\thepage}{\S2-\arabic{page}}

Let $f:[a,b] \ra \R$.
\\

\begin{x}{\small\bf DEFINITION} \ 
Given a closed subinterval $I = [x,y] \subset [a,b]$, put
\[
v(f; I) 
\ = \ 
\abs{f(y) - f(x)}, 
\]
the \un{variation} of $f$ in $I$.
\end{x}

\begin{x}{\small\bf DEFINITION} \ 
Given a partition $P \in \ \sP[a,b]$, put
\begin{align*}
\bigvee\limits_a^b \ (f; P) \ 
&=\ 
\sum\limits_{i = 1}^n \ \abs{f(x_i) - f(x_{i - 1})}
\\[15pt]
&=\ 
\sum\limits_{i = 1}^n \ v(f; I_i) 
\hspace{0.7cm} (I_i = [x_{i - 1}, x_i]),
\end{align*}
the \un{variation} of $f$ in $P$. 
\end{x}

\begin{x}{\small\bf NOTATION} \ 
Put
\[
\tT_f[a,b] 
\ = \ 
\sup\limits_{P \in \ \sP[a,b]} \ \bigvee\limits_a^b \ (f;P),
\]
the \un{total variation} of $f$ in $[a,b]$. 
\end{x}

\begin{x}{\small\bf DEFINITION} \ 
A function $f:[a,b] \ra \R$ is of \un{bounded variation} in $[a,b]$ provided
\[
\tT_f[a,b] 
\ < \ 
+\infty.
\]
\end{x}

\begin{x}{\small\bf NOTATION} \ 
$\BV[a,b]$ is the set of functions of bounded variation in $[a,b]$.
\end{x}


\begin{x}{\small\bf EXAMPLE} \ 
Take $[a,b] = [0,1]$ and define $f:[0,1] \ra \R$ by the rule
\[
f(x) \quad = \quad 
\begin{cases}
\ \text{0 \quad if $x$ is irrational}\\
\ \text{1 \quad if $x$ is rational}
\end{cases}
.
\]
Then $f \notin \ \BV[0,1]$.
\end{x}

\begin{x}{\small\bf NOTATION} \ 
Given $P \in \ \sP[a,b]$, put
\[
\norm{P} 
\ = \ 
\max (x_i - x_{i - 1}) 
\hspace{0.7cm} (i = 1, \ldots, n).
\]
\end{x}

\begin{x}{\small\bf THEOREM} \ 
Let $f \in \ \BV[a,b]$.  
Assume: \ $f$ is continuous $-$then
\[
\tT_f [a,b] 
\ = \ 
\lim\limits_{\norm{P} \ra 0} \ \bigvee\limits_a^b \ (f; P).
\]

[Note: \ 
The continuity assumption is essential.  
E.g., take $[a,b] = [-1, +1]$ and consider 
$f(0) = 1$, $f(x) = 0$ $(x \neq 0)$.]
\end{x}

Let $f:[a,b] \ra \R$.
\\

\begin{x}{\small\bf DEFINITION} \ 
Given a closed subinterval $I = [x,y] \subset [a,b]$, denote by $M$ and $m$ the supremum and infimum of $f$ in $I$ 
and put
\[
\osc(f; I) 
\ = \ 
M - m,
\]
the \un{oscillation} of $f$ in $I$.

[Note: \ 
Since the diameter of $f(I)$ is the supremum of the distances between pairs of points of $f(I)$, it follows that
\[
M - m 
\ = \ 
\diam f(I)
\]
or still, 
\[
\osc (f; I) 
\ = \ 
\diam f(I).
\]
And, of course,
\[
v(f; I) 
\ \leq \ 
\diam f(I).]
\]
\end{x}

Let
\[
v(f; [a,b]) 
\ = \ 
\sup\limits_{P \in \sP[a,b]} \ \sum\limits_{i = 1}^n \ \osc(f; I_i).
\]
\\[-.5cm]

\begin{x}{\small\bf THEOREM} \ 
\[
\tT_f [a,b] 
\ = \ 
v(f; [a,b]).
\]

PROOF \ 
It is obvious that 
\[
\tT_f [a,b] 
\ \leq \ 
v(f; [a,b]).
\]
To go the other way, fix $\varepsilon > 0$.  
Choose a partition $P$ of $[a,b]$ such that if 
$\Delta_i = \osc(f;I_i)$, then
\[
\sigma
\ = \ 
\sum\limits_{i = 1}^n \ \Delta_i
\]
is greater than $v(f; [a,b]) - \varepsilon$ or $\varepsilon^{-1}$ according to whether 
$v(f; [a,b]) < +\infty$ or $v(f; [a,b]) = +\infty$.  
To deal with the first possibility, note that in each interval $I_i = [x_{i-1}, x_i]$ there are two points 
$\xi_i^\prime$, $\xi_i^{\prime\prime}$ with 
\[
\abs{f(\xi_i^{\prime\prime}) - f(\xi_i^\prime)}
\ > \ 
\Delta_i - \frac{\varepsilon}{n}.
\]
The points $\xi_i^\prime$, $\xi_i^{\prime\prime}$ divide $I_i$ into one or two or three subintervals.  
Call
\[
Q 
\ = \ 
\{y_0, \ldots, y_m\} 
\hspace{0.7cm} (n \leq m \leq 3n)
\]
the partition of $[a,b]$ thereby determined $-$then the sum
\[
(i) \hsx \sum\ \abs{f(y_j) - f(y_{j-1})}
\hspace{0.7cm} \text{$([y_{j-1}, y_j]$ contained in $[x_{i-1}, x_i]$)}
\]
is $\Delta_i - \ds\frac{\varepsilon}{n}$.
Therefore
\begin{align*}
\sum\limits_{j = 1}^m \ \abs{f(y_j) - f(y_{j-1})} \ 
&=\ 
\sum\limits_{i = 1}^n \ (i) \sum \abs{f(y_j) - f(y_{j-1})}
\\[15pt]
&> \ 
\sum\limits_{i = 1}^n \ \bigg(\Delta_i -\frac{\varepsilon}{n}\bigg)
\\[15pt]
&=\ 
\sum\limits_{i = 1}^n \ \Delta_i \ - \ 
\frac{\varepsilon}{n} \ \sum\limits_{i = 1}^n \ 1
\\[15pt]
&=\ 
\sigma - \varepsilon
\\[15pt]
&> \ 
v(f; [a,b]) - \varepsilon - \varepsilon,
\end{align*}
from which
\[
\tT_f [a,b] 
\ \geq \ 
v(f; [a,b]).
\]
\end{x}

\chapter{
$\boldsymbol{\S}$\textbf{3}.\quad  FACTS AND EXAMPLES}
\setlength\parindent{2em}
\setcounter{theoremn}{0}
\renewcommand{\thepage}{\S3-\arabic{page}}

\begin{x}{\small\bf FACT} \ 
Suppose that $f \in \BV[a,b]$ $-$then $f$ is bounded on $[a,b]$.
\end{x}

[Given $a \leq x \leq b$, write
\begin{align*}
\abs{f(x)} \ 
&=\ 
\abs{f(x) - f(a) + f(a)}
\\[11pt]
&\leq \ 
\abs{f(x) - f(a)} + \abs{f(a)}
\\[11pt]
&\leq \ 
\abs{f(x) - f(a)} + \abs{f(b) - f(x)} + \abs{f(a)}
\\[11pt]
&\leq \ 
\tT_f\hsx [a,b] + \abs{f(a)} 
\\[11pt]
&< \ 
+\infty.]
\end{align*}
\\[-1cm]

\begin{x}{\small\bf FACT} \ 
A function $f:[a,b] \ra \R$ is constant iff $\tT_f\hsx [a,b] = 0$.
\end{x}

[A constant function certainly has the stated property.  
Conversely, if $f$ is not constant on $[a,b]$, then the claim is that $\tT_f\hsx [a,b] \neq 0$.  
Thus choose $x_1 \neq x_2 \in \ [a,b]$ such that $f(x_1) \neq f(x_2)$, say $x_1 < x_2$ $-$then
\[
\tT_f\hsx [a,b] 
\ \geq  \ 
\abs{f(x_1) - f(a)} + \abs{f(x_2) - f(x_1)} + \abs{f(b) - f(x_2)}
\]

$\implies$
\[
\tT_f\hsx [a,b] 
\ \geq  \ 
\abs{f(x_2) - f(x_1)} 
\ > \ 
0.]
\]
\\[-1cm]

\begin{x}{\small\bf FACT} \ 
If $f:[a,b] \ra \R$  is increasing, then $f \in \BV[a,b]$ and 
\[
\tT_f\hsx [a,b] 
\ = \ 
f(b) - f(a).
\]
\end{x}

[If $P = \{x_0, \ldots, x_n\}$ is a partition of $[a,b]$, then
\begin{align*}
\bigvee\limits_a^b \ (f;P) \ 
&=\ 
\sum\limits_{i = 1}^n \ \abs{f(x_i) - f(x_{i - 1})}
\\[11pt]
&=\  
\sum\limits_{i = 1}^n \ (f(x_i) - f(x_{i - 1}))
\\[11pt]
&=\  
f(b) - f(a).]
\end{align*}

\begin{x}{\small\bf FACT} \ 
If $f:[a,b] \ra \R$ satisfies a Lipschitz condition, then $f \in \BV[a,b]$.
\end{x}

[To say that $f$ satisfies a Lipschitz condition means that there exists a constant $K > 0$ such that for all $x$, $y \in [a,b]$, 
\[
\abs{f(x)  - f(y)} 
\ \leq \ 
K \hsx \abs{x - y}.]
\]
\\[-.75cm]

\begin{x}{\small\bf FACT} \ 
If $f:[a,b] \ra \R$ is differentiable on $[a,b]$ and if its derivative 
$f^\prime:[a,b] \ra \R$ is bounded on $[a,b]$, then $f \in \BV[a,b]$.
\end{x}

[The mean value theorem implies that $f$ satisfies a Lipschitz condition on $[a,b]$.]
\\[-.25cm]

[Note: \ 
Therefore polynomials on $[a,b]$ are in $\BV[a,b]$.]
\\[-.25cm]

\begin{x}{\small\bf FACT} \ 
If $f:[a,b] \ra \R$ has finitely many relative maxima and minima, say at the points
\[
a < \xi_1 < \cdots < \xi_n < b,
\]
then
\begin{align*}
\tT_f\hsx [a,b]  \ 
&=\ 
\abs{f(a) - f(\xi_1)} + \cdots + \abs{f(\xi_n) - f(b)} 
\\[11pt]
&< \ 
+\infty,
\end{align*}
so $f \in \BV[a,b]$.
\end{x}

\begin{x}{\small\bf EXAMPLE} \ 
Take $f(x) = \sin x \ (0 \leq x \leq 2 \hsy \pi)$ $-$then $\tT_f\hsx [0, 2 \hsy \pi] = 4$.
\end{x}

Neither continuity and/or boundedness on $[a,b]$ suffices to force bounded variation.
\\[-.25cm]

\begin{x}{\small\bf EXAMPLE} \ 
Take $[a,b] = [0,1]$ and let
\[
f(x) = \quad 
\begin{cases}
\ x \hsx \sin (1/x) \hspace{0.7cm} (0 < x \leq 1)\\
\hspace{1cm} 0 \hspace{1.5cm} (x = 0)
\end{cases}
.
\]
Then $f(x)$ is continuous and bounded but $f \notin \BV[0,1]$.
\\[-.25cm]

[Note: \ 
On the other hand,
\[
f(x) = \quad 
\begin{cases}
\ x^2 \hsx \sin (1/x) \hspace{0.7cm} (0 < x \leq 1)\\
\hspace{1cm} 0 \hspace{1.63cm} (x = 0)
\end{cases}
\]
is continuous and of bounded variation in $[0,1]$.]
\end{x}

The composition of two functions of bounded variation need not be of bounded variation.
\\

\begin{x}{\small\bf EXAMPLE} \ 
Work on $[0,1]$ and take $f(x) = \sqrt{x}$, 
\[
g(x) = \quad 
\begin{cases}
\ x^2 \hsx \sin^2 (1/x) \hspace{0.7cm} (0 < x \leq 1)\\
\hspace{1cm} 0 \hspace{1.8cm} (x = 0)
\end{cases}
.
\]
Then $f:[0,1] \ra \R$, $g:[0,1] \ra [0,1]$ are of bounded variation but $f \circ g : [0,1] \ra \R$ is not of bounded variation.
\end{x}

\begin{x}{\small\bf FACT} \ 
Suppose that If $f:[a,b] \ra [a,b]$  $-$then the composition $f \circ g \in \BV[a,b]$ for all $g:[a,b] \ra [a,b]$ of bounded variation iff 
$f$ satisfies a Lipschitz condition.
\end{x}

[In one direction, suppose that
\[
\abs{f(x) - f(y)} 
\ \leq \ 
K \abs{x - y} 
\hspace{0.7cm} (x, \ y \in [a,b]).
\]

Let $P \in \sP[a,b]$:
\begin{align*}
\bigvee\limits_a^b \ (f \circ g;P) \ 
&=\ 
\sum\limits_{i = 1}^n \ \abs{(f \circ g) (x_i) - (f \circ g)(x_{i - 1})}
\\[11pt]
&\leq\  
\sum\limits_{i = 1}^n \ K \hsx  \abs{g(x_i) - g(x_{i - 1})}
\\[11pt]
&\leq\   
K \ \bigvee\limits_a^b \ (g;P) 
\\[11pt]
&\leq\  
K \hsy \tT_g [a,b]
\\[11pt]
&< \   
+\infty.]
\end{align*}

\chapter{
$\boldsymbol{\S}$\textbf{4}.\quad  PROPERTIES}
\setlength\parindent{2em}
\setcounter{theoremn}{0}
\renewcommand{\thepage}{\S4-\arabic{page}}

\begin{x}{\small\bf THEOREM} \ 
If $f$, $g \in \BV [a,b]$, then $f + g \in \BV [a,b]$ and 
\[
\tT_{f + g} [a,b] 
\ \leq \ 
\tT_f [a,b] + \tT_g [a,b].
\]
\end{x}

\begin{x}{\small\bf THEOREM} \ 
If $f \in \BV [a,b]$ and $c \in \R$, then $c \hsy f \in \BV [a,b]$ and 
\[
\tT_{c \hsy f} [a,b] 
\ = \ 
\abs{c} \ \tT_f [a,b].
\]
\end{x}

\begin{x}{\small\bf SCHOLIUM} \ 
$\BV [a,b]$ is a linear space.
\\
\end{x}

\begin{x}{\small\bf THEOREM} \ 
If $f$, $g \in \BV [a,b]$, then $f \hsy g \in \BV [a,b]$ and 
\[
\tT_{f \hsy g} [a,b] 
\ \leq \ 
\big(
\sup\limits_{[a,b]} \ \abs{g} 
\big)
\tT_f [a,b] 
 + 
\big(
\sup\limits_{[a,b]} \ \abs{f} 
\big)
\tT_g [a,b].
\]
\\[-1cm]
\end{x}

\begin{x}{\small\bf SCHOLIUM} \ 
$\BV [a,b]$ is an algebra.
\\
\end{x}

\begin{x}{\small\bf THEOREM} \ 
Let $f \in \BV [a,b]$ and let $a < c < b$ $-$then
\[
\begin{cases}
&f \in \BV [a,c]\\
&f \in \BV [c,b]
\end{cases}
\]
and 
\[
\tT_f [a,b] 
\ = \ 
\tT_f [a,c] + \tT_f [c,b].
\]
\\[-1cm]
\end{x}

\begin{x}{\small\bf CRITERION} \ 
Suppose given a function $f : [a,b] \ra \R$ with the property that $[a,b]$ can be divided into a finite number of 
subintervals on each of which $f$ is monotonic $-$then $f \in \BV [a,b]$.
\\
\end{x}

\begin{x}{\small\bf EXAMPLE} \ 
A function of bounded variation need not be monotonic in any subinterval of its domain.
\\[-.5cm]

[Take $[a,b] = [0,1]$ and let $r_1, r_2, \ldots$ be an ordering of the rational numbers in $]0, 1[$.  
Fix $0 < c < 1$ and define $f[0,1] \ra \R$ by 
\[
f(x) \ = \ 
\begin{cases}
& c^k \hspace{0.4cm} (x = r_k)\\
& \ 0 \hspace{0.5cm} \text{otherwise}
\end{cases}
.
\]
Then $f$ is nowhere monotonic but it is of bounded variation in $[0,1]$:
\[
\tT_f [0,1] 
\ = \ 
\frac{2 \hsy c}{1 - c}.]
\]
\end{x}

\begin{x}{\small\bf THEOREM} \ 
\[f \in \BV [a,b]
\implies 
\abs{f} \in \BV [a,b].
\]

Therefore $\BV [a,b]$ is closed under the formation of the conbinations
\[
\begin{cases}
& \ds\frac{1}{2} (f + \abs{f})\\[15pt]
& \ds\frac{1}{2} (f - \abs{f})
\end{cases}
.
\]
\end{x}

\chapter{
$\boldsymbol{\S}$\textbf{5}.\quad  REGULATED FUNCTIONS}
\setlength\parindent{2em}
\setcounter{theoremn}{0}
\renewcommand{\thepage}{\S5-\arabic{page}}

\qquad Given a function $f:[a,b] \ra \R$ and a point $c \in ]a,b[\hsy,$
\[
\begin{cases}
\ f(c+) \ = \  \text{limit from the right $ \ = \ \lim\limits_{x \downarrow c} \ f(x)$}\\[8pt]
\ f(c-)  \ = \ \text{limit from the left $ \ = \ \lim\limits_{x \uparrow c} \ f(x)$}
\end{cases}
.
\]

[Note: \ 
Define $f(a+)$ and $f(b-)$ in the obvious way.]
\\[-.25cm]

\begin{x}{\small\bf DEFINITION} \ 
$f$ is said to be \un{regulated} if 
\[
\begin{cases}
\textbullet \text{$f(c+)$ exists for all $a \leq c < b$.}\\[8pt]
\textbullet \text{$f(c-)$ exists for all $a < c \leq b$}.
\end{cases}
\]
\end{x}

\begin{x}{\small\bf NOTATION} \ 
$\REG[a,b]$ is the set of regulated functions in $[a,b]$.
\end{x}

\begin{x}{\small\bf THEOREM} \ 
$\REG[a,b]$ is a linear space.
\\[-.5cm]

[Sums and scalar multiples of regulated functions are regulated.]
\end{x}

\begin{x}{\small\bf \un{N.B.}} \ 
Continuous functions $f:[a,b] \ra \R$ are regulated, i.e., 
\[
C [a,b] 
\subset 
\REG [a,b].
\]
\end{x}

\begin{x}{\small\bf THEOREM} \ 
Let $f \in \REG [a,b]$ $-$then the discontinuity set of $f$ is at most countable.
\end{x}

\begin{x}{\small\bf DEFINITION} \ 
A function $f:[a,b] \ra \R$ is  \un{right continuous} if for all $a \leq c < b$, 
\[
f(c)  
\ = \ 
f(c+).
\]
\end{x}

\begin{x}{\small\bf DEFINITION} \ 
Let $f \in \REG [a,b]$ $-$then the \un{right continuous modification}
$f_r$ of $f$ is defined by
\[
f_r(x) 
\ = \ 
f(x+) 
\qquad (a \leq x < b).
\]
\end{x}

\begin{x}{\small\bf LEMMA} \ 
Up to an at most countable set, $f_r = f$.
\\[-.5cm]

[The set of points at which $f$ is not right continuous is a subset of the set of points at which $f$ is not continous.]
\end{x}

\begin{x}{\small\bf LEMMA} \ 
$f_r$ is right continuous.
\\[-.25cm]

[For
\[
f_r(c+) 
\ = \ 
\lim\limits_{x \downarrow c} \ f_r(x) 
\ = \ 
\lim\limits_{x \downarrow c} \ f(x) 
\ = \ 
f(c+)
\ = \ 
f_r(c).]
\]
\end{x}

\begin{x}{\small\bf DEFINITION} \ 
Let $f:[a,b] \ra \R$.
\\[-.25cm]

\qquad \textbullet \ If $f(x) = \chisubI (x)$, where $I = [a,b]$ or $]a,b[$, or $[a,b[$, or $]a,b]$, then
$f$ is said to be a \un{single step function}.
\\[-.5cm]

\qquad \textbullet \ 
If $f$ is a finite linear combination of single step functions, then $f$ is said to be a \un{step function}.
\end{x}

\begin{x}{\small\bf LEMMA} \ 
A function $f:[a,b] \ra \R$ is a step function iff there are points
\[
a 
\ = \ 
x_0 
\ < \ 
x_1 
\ < \ 
\cdots 
\ < \ 
x_n 
\ = \ 
b
\]
such that $f$ is constant on each open interval $]x_{i - 1}, x_i[$ \ $(i = 1, \ldots, n)$.
\end{x}

\begin{x}{\small\bf THEOREM} \ 
Let $f:[a,b] \ra \R$ $-$then $f$ is regulated iff $f$ is a uniform limit of a sequence of step functions.
\end{x}

\begin{x}{\small\bf \un{N.B.}} \ 
Regulated functions are bounded.
\\[-.25cm]

[Take an $f \in \REG [a,b]$ and choose a step function $g$ such that $\norm{f - g}_\infty \leq 1$, 
hence $\forall \ x \in [a,b]$, 
\[
\abs{f(x)} 
\ \leq \ 
\norm{f - g}_\infty + \norm{g}_\infty 
\ \leq \ 
1 + \norm{g}_\infty.]
\]
\end{x}

\begin{x}{\small\bf THEOREM} \ 
Let $f \in \BV [a,b]$ $-$then $f$ is regulated.
\\[-.5cm]

PROOF \ 
Suppose that $a < c \leq b$ and $f(c-)$ does not exist $-$then there is a positive number 
$\varepsilon$ and a sequence of real numbers $c_k$ increasing to $c$ such that for all $k$, 
\[
f(c_k) - f(c_{k + 1}) 
\ < \ 
- \varepsilon
\ < \ 
\varepsilon
\ < \ 
f(c_{k + 2}) - f(c_{k +1}). 
\]
It therefore follows that for all $n$, 
\[
+\infty 
\ > \ 
\tT_f [a,b] 
\ \geq \ 
\sum\limits_{k = 1}^n \ \abs{f(c_k) - f(c_{k + 1}) } 
\ > \ 
n \hsy \varepsilon,
\]
an impossibility.  
In the same vein, $f(c+)$ must exist for all $a \leq c < b$.
\end{x}

\begin{x}{\small\bf SCHOLIUM} \ 
\[
\BV [a,b]
\subset 
\REG [a,b].
\]
\end{x}

In particular: 
The discontinuity set of an $f \in \BV[a,b]$ is at most countable.
\\[-.25cm]

\begin{x}{\small\bf THEOREM} \ 
$\REG [a,b]$ is a Banach space in the uniform norm and $\BV [a,b]$ is a dense linear subspace of $\REG [a,b]$, thus
\[
\ov{\BV [a,b]} 
\ = \ 
\REG [a,b] 
\]
per $\norm{\hsx \cdot \hsx}_\infty$.
\end{x}

\chapter{
$\boldsymbol{\S}$\textbf{6}.\quad  POSITIVE AND NEGATIVE}
\setlength\parindent{2em}
\setcounter{theoremn}{0}
\renewcommand{\thepage}{\S6-\arabic{page}}

\begin{x}{\small\bf NOTATION} \ 
Given a real number $x$, put
\[
\begin{cases}
\ x^+ = \max (x, 0) = \ds\frac{1}{2} \ (\abs{x} + x)\\[8pt]
\ x^- = \max (-x, 0) = \ds\frac{1}{2} \ (\abs{x} - x)
\end{cases}
.
\]

Given a function $f:[a,b] \ra \R$, let
\[
\begin{cases}
\ \tT_f^+ [a,b]  \ = \ \sup\limits_{P \in \sP [a,b]} \ \sum\limits_{i = 1}^n \ \left(f(x_i) - f(x_{i - 1}) \right)^+\\[15pt]
\ \tT_f^- [a,b]  \ = \ \sup\limits_{P \in \sP [a,b]} \ \sum\limits_{i = 1}^n \ \left(f(x_i) - f(x_{i - 1}) \right)^-
\end{cases}
,
\]
the 
\[
\begin{cases}
\ \text{\un{positive}}\\
\ \text{\un{negative}}
\end{cases}
\ \text{\un{total variation}} 
\hspace{2cm}
\]
of $f$ in $[a,b]$.
\\[-.25cm]

Obviously, 
\[
\begin{cases}
\ 0 \leq \tT_f^+ [a,b] \leq \tT_f [a,b] \leq +\infty\\[8pt]
\ 0 \leq \tT_f^- [a,b] \leq \tT_f [a,b] \leq +\infty
\end{cases}
,
\]
so $\tT_f^+ [a,b]$, $\tT_f^- [a,b]$, $\tT_f [a,b]$ are all finite if $f \in \BV [a,b]$.
\\[-.5cm]
\end{x}

\begin{x}{\small\bf \un{N.B.}} \ 
Abbreviate
\[
\sum\limits_{i = 1}^n \ \left(f(x_i) - f(x_{i - 1})\right) 
\ \text{to} \ \Sigma, 
\hspace{0.7cm}
\]
\[
\sum\limits_{i = 1}^n \ \left(f(x_i) - f(x_{i - 1})\right)^+ 
\ \text{to} \ \Sigma^+,
\]
\[
\sum\limits_{i = 1}^n \ \left(f(x_i) - f(x_{i - 1})\right)^- 
\ \text{to} \ \Sigma^-.
\]
Then
\[
\Sigma^+ + \Sigma^- 
\ = \ 
\Sigma, 
\quad 
\Sigma^+ - \Sigma^- 
\ = \ 
f(b) - f(a)
\]

$\implies$
\[
2 \hsy \Sigma^+ 
\ = \ 
\Sigma + f(b) - f(a), 
\quad
2 \hsy \Sigma^- 
\ = \ 
\Sigma - f(b) + f(a).
\]
\\[-1cm]
\end{x}

\begin{x}{\small\bf THEOREM} \ 
If $f \in \BV [a,b]$, then
\[
\begin{cases}
\ \tT_f^+ [a,b] + \tT_f^- [a,b]  = \tT_f [a,b] \\[8pt]
\ \tT_f^+ [a,b] - \tT_f^- [a,b]  = f(b) - f(a)
\end{cases}
.
\]
\end{x}

Replace ``$b$'' by ``$x$'' and assume that $f \in \BV [a,b]$.
\\[-.25cm]

\qquad \textbullet \quad 
$\tT_f^+ [ax] \ = \ 2^{-1} \hsx \left(\tT_f [a,x] + f(x) - f(a) \right) $
\\[-.5cm]

\hspace{2cm} $\implies$
\[
\frac{1}{2} \ \left(\tT_f [a,x] + f(x) \right)
\ = \ 
\tT_f^+ [a,x] + 2^{-1} f(a).
\]

\qquad \textbullet \quad 
$\tT_f^- [a,x] \ = \ 2^{-1} \hsx \left(\tT_f [a,x] - f(x) + f(a) \right) $
\\[-.5cm]

\hspace{2cm} $\implies$
\[
\frac{1}{2} \ \left(\tT_f [a,x] - f(x) \right)
\ = \ 
\tT_f^+ [a,x] - 2^{-1} f(a).
\]
\\[-.75cm]

\begin{x}{\small\bf LEMMA} \ 
The functions
\[
\begin{cases}
\ x \ra \ds\frac{1}{2} \ \left(\tT_f [a,x] + f(x) \right) \\[8pt]
\ x \ra \ds\frac{1}{2} \ \left(\tT_f [a,x] - f(x) \right) 
\end{cases}
, \tT_f [a,a] = 0
\]
are increasing.  
\\[-.5cm]

PROOF \ 
Let $a \leq x < y \leq b$.
\begin{align*}
\text{\textbullet} \quad \frac{1}{2} \ \left(\tT_f [a,y] + f(y) \right) - \frac{1}{2} \ \left(\tT_f [a,x] + f(x) \right) \ 
&= \ 
\frac{1}{2} \ \left(\tT_f [a,y] - \tT_f [a,x] + f(y) - f(x) \right)
\\[11pt]
&\geq \ 
\frac{1}{2} \ \left(\tT_f [x,y] - \abs{f(y) - f(x)}\right)
\\[11pt]
&\geq \ 
0.
\end{align*}
\begin{align*}
\text{\textbullet} \quad \frac{1}{2} \ \left(\tT_f [a,y] - f(y) \right) - \frac{1}{2} \ \left(\tT_f [a,x] - f(x) \right) \ 
&= \ 
\frac{1}{2} \ \left(\tT_f [a,y] - \tT_f [a,x] - f(y) +f(x) \right)
\\[11pt]
&\geq \ 
\frac{1}{2} \ \left(\tT_f [x,y] - \abs{f(y) - f(x)}\right)
\\[11pt]
&\geq \ 
0.
\end{align*}
\end{x}

\begin{x}{\small\bf DEFINITION} \ 
The representation
\[
f(x) 
\ = \ 
\frac{1}{2} \ \left(\tT_f [a,x] + f(x) \right) 
- 
\frac{1}{2} \ \left(\tT_f [a,x] - f(x) \right) 
\]
is the \un{Jordan decomposition} of $f$.
\end{x}

\begin{x}{\small\bf REMARK} \ 
To arrive at a representation of $f$ as the difference of two strictly increasing functions, write
\[
f(x) 
\ = \ 
\left(
\frac{1}{2} \ \left(\tT_f [a,x] + f(x) \right) + x \right) 
\hsx - \hsx 
\left(
\frac{1}{2} \ 
\left(\tT_f [a,x] - f(x) \right) 
+ x
\right).
\]
\\[-1.25cm]
\end{x}

\begin{x}{\small\bf THEOREM} \ 
Suppose that $f \in \BV [a,b]$ $-$then $f$ is Borel measurable.
\\[-.5cm]

[For this is the case of an increasing function.]
\end{x}

\chapter{
$\boldsymbol{\S}$\textbf{7}.\quad  CONTINUITY}
\setlength\parindent{2em}
\setcounter{theoremn}{0}
\renewcommand{\thepage}{\S7-\arabic{page}}

\begin{x}{\small\bf THEOREM} \ 
Let $f \in \BV [a,b]$.  
Suppose that $f$ is continuous at $c \in [a,b]$ $-$then $\tT_f[a,-]$ is continuous at $c \in [a,b]$.
\\[-.25cm]

PROOF \ 
The function $x \ra \tT_f [a,x]$ is increasing, hence both one sided limits exist at all points $c \in [a,c]$, the claim being that
\[
\lim\limits_{x \ra c} \ \tT_f [a,x]
\ = \ 
\tT_f [a,c].
\]
To this end, it will be shown that the right hand limit of $\tT_f [a,x]$ as $x \ra c$ is equal to $\tT_f [a,c]$, 
where $a \leq c < b$, the discussion for the left hand limit being analogous.  
So let $\varepsilon > 0$ and chooose $\delta > 0$ such that 
\[
0 \ < \ x - c \ < \delta 
\implies 
\abs{f(x) - f(c)} \ < \ \frac{\varepsilon}{2}.
\]
Partition $[c,b]$ by the scheme
\[
\tT_f [c,b] 
\ < \ 
 \sum\limits_{i = 1}^n \ \abs{f(x_i) - f(x_{i - 1})}
\ + \  \frac{\varepsilon}{2} 
 \qquad \left(x_0 = c, \ x_n = b\right).
\]
If $x_1 - c < \delta$, then

\begin{align*}
\tT_f [c, b] - \frac{\varepsilon}{2} 
&<\ 
\abs{f(x_1) - f(c)} \ + \  \sum\limits_{i = 2}^n \ \abs{f(x_i) - f(x_{i - 1})}
\\[15pt]
&< \ 
\frac{\varepsilon}{2} + \tT_f [x_1, b]
\end{align*}

$\implies$
\[
\tT_f [c, b] - \tT_f [x_1, b]
\ < \ 
\varepsilon.
\]
On the other hand, if $x_1 - c \hsx \geq \hsx  \delta$, add a point $x$ to the partition subject to
$x - c \hsx < \hsx  \delta$, thus
\begin{align*}
\tT_f [c, b] - \frac{\varepsilon}{2} 
&<\ 
\abs{f(x_1) - f(x_0)} \ + \   \sum\limits_{i = 2}^n \ \abs{f(x_i) - f(x_{i - 1})}
\\[15pt]
&\leq \ 
\abs{f(x_1) - f(x)} \ + \   \abs{f(x) - f(x_0)} \ + \   \sum\limits_{i = 2}^n \ \abs{f(x_i) - f(x_{i - 1})}
\\[15pt]
&<\ 
\abs{f(x_1) - f(x)} \ + \   \frac{\varepsilon}{2} \ + \   \sum\limits_{i = 2}^n \ \abs{f(x_i) - f(x_{i - 1})}.
\end{align*}
Since
\[
\{x, \hsx  x_1, \hsx \ldots, \hsx x_n \}
\]
is a partition of $[x,b]$, it follows that
\[
\tT_f [c,b] -  \frac{\varepsilon}{2}
\ < \ 
 \frac{\varepsilon}{2} + \tT_f [x, b]
\]

$\implies$
\[
\tT_f [c,b]  - \tT_f [x, b] 
\ < \ 
\varepsilon.
\]
Finally
\begin{align*}
\tT_f [a,b] - \tT_f [x, b] \ 
&=\ 
\tT_f [c,x] 
\\[4pt]
&=\ 
\tT_f [a,x] - \tT_f [a,c]
\\[4pt]
&<\ 
\varepsilon
\end{align*}
if $x - c \hsx < \hsx  \delta$.  
Therefore
\[
\tT_f [a,c+] 
\ = \ 
\tT_f [a,c],
\]
so $\tT_f [a,x]$ is right continuous at $c$.
\\[-.25cm]
\end{x}

\begin{x}{\small\bf SCHOLIUM} \ 
If $f \in \BV [a,b] \capx C [a,b]$, then 
\[
 \tT_f [a,-] \in C [a,b].
\]
\end{x}


\begin{x}{\small\bf REMARK} \ 
It is also true that
\[
\begin{cases}
\ \tT_f^+ [a,-] \in C [a,b]\\[4pt]
\ \tT_f^- [a,-] \in C [a,b]
\end{cases}
.
\]
\\[-1cm]

Proof: \ 

\[
\begin{cases}
\ \tT_f^+ [a,x] \ = \ 2^{-1} \left(\tT_f [a,x] + f(x) - f(a)\right) \\[4pt]
\ \tT_f^- [a,x]  \ = \ 2^{-1} \left(\tT_f [a,x] - f(x) + f(a)\right)
\end{cases}
.
\]
\\[-.75cm]
\end{x}

\begin{x}{\small\bf THEOREM} \ 
If $f \in \BV [a,b]$ is continuous, then $f$ can be written as the difference of two increasing continuous functions.
\\[-.5cm]

[In view of what has been said above, this is obvious.]
\\[-.25cm]
\end{x}

\begin{x}{\small\bf LEMMA} \ 
Let $f \in \BV [a,b]$.  
Assume: $\tT_f [a,-]$ is continuous at $c \in [a,b]$ $-$then $f$ is continuous at $c \in [a,b]$.
\\[-.25cm]

PROOF \ For
\[
\begin{cases}
\ c < x \implies \abs{f(x) - f(c)} \ \leq \  \tT_f [c, x] \ = \  \tT_f [a,x] - \tT_f [a,c]\\[4pt]
\ x < c \implies \abs{f(c) - f(x)} \ \leq \ \tT_f [x, c] \ = \  \tT_f [a,c] - \tT_f [a,x]
\end{cases}
.
\]
\\[-.75cm]
\end{x}

\begin{x}{\small\bf RAPPEL} \ 
Let $f : [a,b] \ra \R$ be increasing and let $x_1, \hsx x_2, \ldots$ be an enumeration of the interior points of discontinuity of $f$ 
$-$then the \un{saltus function} $s_f : [a,b] \ra \R$ attached to $f$ is defined by 
\[
s_f(a) \ = \ 0
\]
and if $a < x \leq b$, by
\[
s_f(x) 
\ = \ 
(f(a+) - f(a)) 
+ 
\sum\limits_{x_k < x} \ 
(f(x_k+) - f(x_k-) )
\hsx + \hsx (f(x) - f(x-)).
\]
\\[-1cm]
\end{x}

\begin{x}{\small\bf FACT} \ 
The difference \hsx $f - s_f$ \hsx is an increasing continuous function.  
\\[.25cm]
Assume again that $f \in \BV [a,b]$ and put
\[
V(x) 
\ = \ 
\tT_f [a,x], 
\quad 
F(x) \ = \ V(x) - f(x) 
\qquad (a \leq x \leq b).
\]

\end{x}

\begin{x}{\small\bf \un{N.B.}} \ 
$V(x)$ and $F(x)$ are increasing functions of $x$.
\end{x}

Let
\[
\{x_1, \hsx x_2, \hsx \ldots\} \quad (a < x_k < b)
\]
be the set comprised of the discontinuity points of $V$.
\\[-.25cm]

\begin{x}{\small\bf REMARK} \ 
The discontinuity set of $V$ coincides with the discontinuity set of $f$ and the discontinuity set of $F$ is contained 
in the discontinuity set of $f$.
\end{x}

Introduce
\[
s_V(x) 
\ = \ 
(V(a+) - V(a)) 
+ 
\sum\limits_{x_k < x} \ (V(x_k+) - V(x_k-))
\hsx + \hsx (V(x) - V(x-)),
\]
and 
\[
s_F(x) 
\ = \ 
(F(a+) - F(a)) 
+
\sum\limits_{x_k < x} \ (F(x_k+) - F(x_k-))
\hsx + \hsx (F(x) - F(x-)),
\]
where $a < x \leq b$ and take
\[
s_V(a) = 0, 
\quad
s_F(a) = 0.
\]
\\[-1cm]

\begin{x}{\small\bf LEMMA} \ 
$s_V$ is the saltus function of $V$ and $s_F$ is the saltus function of $F$.  
\\[-.5cm]

[Per $V$, this is true by its very construction.  
As for $F$, if $x_k$ is not a discontinuity point, then
\[
F(x_k+) - F(x_k-) 
\ = \ 0,
\]
thus such a term does not participate.]
\\[-.25cm]
\end{x}

\begin{x}{\small\bf DEFINITION} \ 
The \un{saltus function} $s_f : [a,b] \ra \R$ attached to $f$ is the difference
\[
s_f 
\ = \ 
s_V - s_F.
\]

Spelled out, 
\[
s_f(a) 
\ = \ 0
\]
and
\[
s_f(x) 
\ = \ 
(f(a+) - f(a)) 
\hsx + \hsx 
\sum\limits_{x_k < x} \ (f(x_k+) - f(x_k-))
\hsx + \hsx (f(x) - f(x-))
\]
subject to $a < x \leq b$.
\\[-.25cm]
\end{x}

\begin{x}{\small\bf SCHOLIUM} \ 
The functions
\[
x \ \ra \ 
\begin{cases}
\ V(x) - s_V(x)\\
\ F(x) - s_F(x)
\end{cases}
\]
are increasing and continuous.  
Therefore
\begin{align*}
f(x) - s_f(x) \ 
&=\ 
V(x) - F(x) - (s_V(x) - s_F(x))
\\[4pt]
&=\ 
(V(x) - s_V(x)) - (F(x) - s_F(x))
\end{align*}
is a continuous function of bounded variation.
\end{x}


\chapter{
$\boldsymbol{\S}$\textbf{8}.\quad  ABSOLUTE CONTINUITY I}
\setlength\parindent{2em}
\setcounter{theoremn}{0}
\renewcommand{\thepage}{\S8-\arabic{page}}

\begin{x}{\small\bf DEFINITION} \ 
A function $f:[a,b] \ra \R$ is \un{absolutely continuous} if $\forall$ $\varepsilon > 0$, 
$\exists$ $\delta > 0$ such that
\[
\sum\limits_{k = 1}^n \ \abs{f(b_k) - f(a_k)} 
\ < \ 
\varepsilon
\]
whenever
\[
a 
\ \leq \
 a_1 
\ < \ 
b_1 
\ \leq \
a_2 
\ < \ 
b_2
\ \leq \
\cdots
\ \leq \
a_n 
\ < \ 
b_n
\ \leq \
b
\]
for which 
\[
\sum\limits_{k = 1}^n \ (b_k - a_k) 
\ < \ 
\delta.
\]
\\[-1cm]
\end{x}

\begin{x}{\small\bf NOTATION} \ 
$\AC[a,b]$ is the set of absolutely continuous functions in $[a,b]$.
\\[-.25cm]
\end{x}

\begin{x}{\small\bf THEOREM} \ 
An absolutely continuous function is uniformly continuous.
\\[-.25cm]
\end{x}

\begin{x}{\small\bf THEOREM} \ 
\[
f \in \AC[a,b] \implies \abs{f} \in \AC[a,b].
\]
\\[-1.5cm]
\end{x}

\begin{x}{\small\bf THEOREM} \ 
If $f$, $g \in \AC[a,b]$, then so do their sum, difference, and product.
\\[-.25cm]
\end{x}

\begin{x}{\small\bf THEOREM} \ 
\[
\AC[a,b] \subset \BV[a,b].
\]
\\[-1.5cm]
\end{x}

\begin{x}{\small\bf SCHOLIUM} \ 
If $f \in C [a,b]$ but $f \notin \BV[a,b]$, then $f \notin \AC[a,b]$.
\\[-.25cm]
\end{x}

\begin{x}{\small\bf CRITERION} \ 
If $f$ is continuous in $[a,b]$ and if $f^\prime$ exists and is bounded in $]a,b[$, then $f$ is absolutely continuous in $[a,b]$.
\\[-.25cm]

[Define $M > 0$ by $\abs{f^\prime(x)} < M$ for all $x$ in $]a,b[$.  
Take $\varepsilon > 0$ and consider
\[
\sum\limits_{k = 1}^n \  \abs{f(b_k) - f(a_k)},
\]
where 
\[
\sum\limits_{k = 1}^n \  (b_k - a_k) 
\ < \ 
\frac{\varepsilon}{M}.
\]
Owing to the Mean Value Theorem, $\exists$ $x_k \in ]a_k, b_k[$ such that 
\[
\frac{{f(b_k) - f(a_k)}}{{b_k - a_k}} 
\ = \ 
f^\prime(x_k).
\]
Therefore 
\begin{align*}
\sum\limits_{k = 1}^n \  \abs{f(b_k) - f(a_k)}  \ 
&=\ 
\sum\limits_{k = 1}^n \  \abs{\frac{f(b_k) - f(a_k)}{b_k - a_k}} \hsx \abs{b_k - a_k}
\\[15pt]
&=\ 
\sum\limits_{k = 1}^n \  \abs{f^\prime(x_k)} \hsx \abs{b_k - a_k}
\\[15pt]
&< \ 
\sum\limits_{k = 1}^n \  M \hsx  \abs{b_k - a_k}
\\[15pt]
&=\ 
M \hsx \sum\limits_{k = 1}^n \  \abs{b_k - a_k}
\\[15pt]
&< \ 
M \hsx \frac{\varepsilon}{M}
\\[15pt]
&=\ 
\varepsilon.]
\end{align*}
\\[-1.5cm]
\end{x}

\begin{x}{\small\bf EXAMPLE} \ 
It can happen that a continuous function with an unbounded derivative is absolutely continuous.
\\[-.5cm]

[Consider $f(x) = \sqrt{x}$ \  $(0 \leq x \leq 1)$ $-$then $f \in \AC[0,1]$ but
\[
f^\prime(x) 
\ = \ 
\frac{1}{2\sqrt{x}}
\hps (0 < x < 1).]
\]
\\[-1.5cm]
\end{x}

\begin{x}{\small\bf EXAMPLE} \ 
Consider
\[
f(x) \ = \quad
\begin{cases}
\ x^2 \sin \bigg(\ds\frac{1}{x}\bigg) \hps (0 < x \leq 1) \\[8pt]
\ \quad 0 \hspace{2.0cm} (x = 0)
\end{cases}
.
\]
Then $f \in \BV[0,1]$.  
But more is true, viz.  $f \in \AC[0,1]$.  
In fact, in $]0,1[$, 
\[
f^\prime(x) 
\ = \ 
2 \hsy x \hsx  \sin \bigg(\ds\frac{1}{x}\bigg) -  \cos \bigg(\ds\frac{1}{x}\bigg)
\]
$\implies$
\begin{align*}
\abs{f^\prime(x)} \ 
&\leq \ 
2 \abs{x} \big|\sin \bigg(\ds\frac{1}{x}\bigg)\big|\ + \big|\cos \bigg(\ds\frac{1}{x}\bigg)\big|
\\[8pt]
&\leq 
3.
\end{align*}
\\[-1.25cm]
\end{x}

\begin{x}{\small\bf THEOREM}  \hsx 
\text{Let $f \in \BV[a,b]$ $-$then $f \in \AC[a,b]$ iff $T_f[a,-] \in \AC[a,b]$.}
\\[-.25cm]

PROOF \ 
Suppose first that $f$ is absolutely continuous.  
Given $\varepsilon > 0$, introduce the pairs
\[
\{(a_1, b_1), \ (a_2, b_2), \ldots, (a_n, b_n)\}
\]
subject to
\[
\sum\limits_{k = 1}^n \ (b_k - a_k) 
\ < \ 
\delta,
\]
thus
\[
\sum\limits_{k = 1}^n \ \abs{f(b_k) - f(a_k)}
\ < \ 
\varepsilon.
\]
For each $k$, let 
\[
P_k : \ a_k = x_{k_0} < x_{k_1} \hsx < \hsx \cdots \hsx < \hsx x_{k_{n_k}} = b_k
\]
be a partition of $[a_k, b_k]$ $-$then
\begin{align*}
\sum\limits_{k = 1}^n \
\sum\limits_{i = 1}^{n_k} \ 
\big(x_{k_i} - x_{k_{i - 1}}\big) \ 
&= \ 
\sum\limits_{k = 1}^n \ (b_k - a_k) 
\\[11pt]
&< \ 
\delta
\end{align*}
$\implies$
\[
\sum\limits_{k = 1}^n \
\sum\limits_{i = 1}^{n_k} \ 
\abs{f(x_{k_i}) - f(x_{k_{i - 1}})}
\ < \ 
\varepsilon.
\]
Vary now the $P_k$ through $\sP([a_k,b_k])$ and take the supremum, hence 
\[
\sum\limits_{k = 1}^n \ T_f[a_k, b_k]
\ < \ 
\varepsilon
\]
or still, 
\[
\sum\limits_{k = 1}^n \
T_f[a, b_k] - T_f[a, a_k]
\ < \ 
\varepsilon.
\]
So $T_f[a,-] \in \AC[a,b]$.  
In the other direction, simply note that
\[
\abs{f(b_k) - f(a_k)} 
\ \leq \ 
T_f[a, b_k] - T_f[a, a_k].
\]
\\[-1.5cm]
\end{x}

Recall that the Jordan decomposition of $f$ is the representation 
\[
f(x) 
\ = \ 
\frac{1}{2} (T_f[a, x] + f(x))
-
\frac{1}{2} (T_f[a, x] - f(x)).
\]
\\[-1.5cm]

\begin{x}{\small\bf SCHOLIUM} \ 
If $f \in \AC[a,b]$, then $f$ can be represented as the difference of two increasing absolutely continuous functions. 
\\[-.25cm]
\end{x}

Here is a useful technicality.
\\

\begin{x}{\small\bf LEMMA} \ 
Suppose that $f:[a,b] \ra \R$ is absolutely continuous $-$then 
$\forall$ $\varepsilon > 0$, $\exists$ $\delta > 0$ such that for any arbitrary finite or countable system of pairwise
disjoint open intervals $\{(a_k, b_k)\}$ with 
\[
\sum\limits_k \ (b_k - a_k) 
\ < \ 
\delta,
\]
the inequality
\[
\sum\limits_k \ \osc(f; [a_k, b_k])
\ < \ 
\varepsilon
\]
obtains.
\\[-.25cm]
\end{x}

\begin{x}{\small\bf DEFINITION} \ 
A function $f:[a,b] \ra \R$ is said to have \un{property (N)} 
if $f$ sends sets of Lebesgue measure 0 to sets of Lebesgue measure 0: 
\[
E \subset [a,b] \quad \& \quad \lambda(E) = 0 
\ \implies \ 
\lambda (f(E)) = 0.
\]
\\[-1.5cm]
\end{x}

\begin{x}{\small\bf THEOREM} \ 
If $f:[a,b] \ra \R$  is absolutely continuous, then $f$ has property (N).  
\\[-.25cm]

PROOF \ 
Suppose that $\lambda(E) = 0$ and assume that $a \notin E$, $b \notin E$ 
(this omission has no bearing on the final outcome).  
Notationally $\varepsilon$, $\delta$, and $\{[a_k,b_k]\}$ are per \#13, thus
\[
\sum\limits_k \  (b_k - a_k) < \delta
\implies
\sum\limits_k \ \osc(f; [a_k, b_k])
< 
\varepsilon.
\]
To fix the data and thereby pin matters down, start by putting
$
\begin{cases}
\ m_k = \min\limits_{[a_k, b_k]} f\\
\ M_k = \max\limits_{[a_k, b_k]} f
\end{cases}
, 
$
hence
\[
\osc(f; [a_k, b_k])
\ = \ 
M_k - m_k.
\]
Since $\lambda(E) = 0$, there exists an open set $S \subset [a,b]$ such that 
\[
E \subset S, \quad \lambda(S) < \delta.
\]
Decompose $S$ into its connected components $]a_k, b_k[$, so 
\[
\sum\limits_k \ (b_k - a_k) 
\ < \ 
\delta.
\]
Next
\begin{align*}
f(E) \subset f(S) \ 
&= \ 
\sum\limits_k \ f(]a_k, b_k[)
\\[15pt]
&\subset \
\sum\limits_k \ f([a_k, b_k])
\end{align*}
or still
\[
\lambda^*(f(E)) 
\ \leq \ 
\sum\limits_k \ \lambda^*(f([a_k, b_k])). 
\]
But
\[
f([a_k, b_k])
\ = \ 
[m_k, M_k].
\]
Therefore
\begin{align*}
\lambda^*(f(E)) \ 
&\leq \ 
\sum\limits_k \ (M_k - m_k) 
\\[15pt]
&< \ 
\varepsilon.
\end{align*}
Since $\varepsilon$ is arbitrary, it follows that
\[
\lambda (f(E)) 
\ = \ 
0.
\]
\\[-1.5cm]
\end{x}

\begin{x}{\small\bf THEOREM} \ 
If $f:[a,b] \ra \R$ is continuous, then $f$ has property (N) iff for every Lebesgue measurable set $E \subset [a,b]$, 
$f(E)$ is Lebesgue measurable. 
\\[-.25cm]

PROOF \ 
Assuming that $f$ has property (N), take an $E$ and write
\[
E 
\ = \ 
\bigg( \bigcup\limits_{j = 1}^\infty \ K_j \bigg) \hsx \cup \hsx S 
\hps (K_1 \subset K_2 \subset \cdots ),
\]
where each $K_j$ is compact and $S$ has Lebesgue measure 0.  
Since $f$ is continuous, $f(K_j)$ is compact, hence
\[
\bigcup\limits_{j = 1}^\infty \ f(K_j ) 
\]
is Lebesgue measurable.  
But  $f$  has property (N), hence  $f(S)$  has Lebesgue $\text{measure 0.}$ 
Therefore
\[
f(E) 
\ = \ 
\bigg( \bigcup\limits_{j = 1}^\infty \ f(K_j) \bigg) \hsx \cup \hsx f(S) 
\]
is Lebesgue measurable.  
In the other direction, suppose that $f$ does not possess property (N), thus that there exists a set $E \subset [a,b]$ of Lebesgue measure 0 
such that $f(E)$ is not a set of Lebesgue measure 0.
\\[-.25cm]

\qquad \textbullet \quad If $f(E)$ is Lebesgue measurable, then it contains a nonmeasurable subset.
\\[-.25cm]

\qquad \textbullet \quad If $f(E)$ is not Lebesgue measurable, then it contains (is \ldots) a nonmeasurable set.
\\[-.5cm]

So there exists a nonmeasurable set $A \subset f(E)$.  
Put $S = f^{-1}(A) \cap E$: $S$ is Lebesgue measurable (being a subset of $E$, a set of Lebesgue measure 0), 
yet $f(S) = A$ is not Lebesgue measurable.
\\[-.25cm]
\end{x}

\begin{x}{\small\bf SCHOLIUM} \ 
An absolutely continuous function sends Lebesgue measurable  sets to Lebesgue measurable sets.
\\[-.25cm]
\end{x}

\begin{x}{\small\bf REMARK} \ 
Let $E \subset [a,b]$ be Lebesgue measurable  $-$then its image $f(E)$ under a continuous function $f[a,b] \ra \R$ need not be Lebesgue measurable.
\\[-.25cm]
\end{x}

\begin{x}{\small\bf RAPPEL} \ 
If $E \subset \R$ is a set of Lebesgue measure 0, then its complement $E^\tc$ is a dense subset of $\R$.
\\[-.25cm]

[In fact, $E^\tc \cap I \neq \emptyset$ for every open interval $I$.]
\\[-.25cm]
\end{x}

\begin{x}{\small\bf LEMMA} \ 
Suppose that $f, \hsx g:[a,b] \ra \R$ are continuous and $f = g$ almost everywhere $-$then $f = g$.
\\[-.25cm]

[The set
\[
E 
\ = \ 
\{x \in [a,b] : f(x) \neq g(x) \}
\]
is a set of Lebesgue measure 0.]
\\[-.25cm]
\end{x}

\begin{x}{\small\bf APPLICATION} \ 
Two absolutely continuous functions which are equal almost everywhere are equal.
\end{x}

\chapter{
$\boldsymbol{\S}$\textbf{9}.\quad  DINI DERIVATIVES}
\setlength\parindent{2em}
\setcounter{theoremn}{0}
\renewcommand{\thepage}{\S9-\arabic{page}}

\begin{x}{\small\bf DEFINITION} \ 
Let $f:[a,b] \ra \R$.
\\[-.25cm]

\qquad \textbullet \ Given $x \in [a,b[$, 
\[
(\tD^+ \hsy f) (x) 
\ = \ 
\limsup\limits_{h \downarrow 0} \ \frac{f(x + h) - f(x)}{h}
\]
is the \un{upper right derivative} of $f$ at $x$ and
\[
(\tD_+ \hsy f) (x) 
\ = \ 
\liminf\limits_{h \downarrow 0} \ \frac{f(x + h) - f(x)}{h}
\]
is the \un{lower right derivative} of $f$ at $x$.
\\[-.25cm]

\qquad \textbullet \ Given $x \in ]a,b]$, 
\[
(\tD^- \hsy f) (x) 
\ = \ 
\limsup\limits_{h \uparrow 0} \ \frac{f(x + h) - f(x)}{h}
\]
is the \un{upper left derivative} of $f$ at $x$ and
\[
(\tD_- \hsy f) (x) 
\ = \ 
\liminf\limits_{h \uparrow 0} \ \frac{f(x + h) - f(x)}{h}
\]
is the \un{lower left derivative} of $f$ at $x$.
\\[-.25cm]
\end{x}

\begin{x}{\small\bf \un{N.B.}} \ 
Collectively, these are the \un{Dini derivatives}.
\\[-.25cm]
\end{x}

\begin{x}{\small\bf EXAMPLE} \ 
Suppose that $a < b$ and $c < d$.  
Let
\[
f(x) \ = \quad
\begin{cases}
\ a \hsy x \hsx \left(\sin \frac{1}{x}\right)^2 
+ \hsx b \hsy x \hsx \left(\cos \frac{1}{x}\right)^2
\hspace{0.7cm} (x > 0)\\[11pt]
\ \hps 0 
\hspace{4.45cm} (x = 0) \\[11pt]
\ c \hsy x \hsx \left(\sin \frac{1}{x}\right)^2 
+ \hsx d \hsy x \hsx \left(\cos \frac{1}{x}\right)^2
\hspace{0.75cm} (x < 0)
\end{cases}
.
\]
Tben
\[
\begin{cases}
\ (D^+ f) (0) \ = \  b  \ > \  a \ = \ (D_+ f)(0)\\[4pt]
\ (D^- f) (0) \ = \  d  \ > \  c \ = \ (D_- f)(0)
\end{cases}
.
\]
\\[-.75cm]
\end{x}

If $(D^+ f)(x) = (D_+ f)(x)$, then the common value is called the \un{right derivative} of $f$ at $x$, denoted $(D_r f)(x)$, 
and $f$ is said to be \un{right differentiable} at $x$ if this common value is finite. 
\\[-.5cm]

If $(D^- f)(x) = (D_- f)(x)$, then the common value is called the \un{left derivative} of $f$ at $x$, denoted $(D_\ell f)(x)$, 
and $f$ is said to be \un{left differentiable} at $x$ if this common value is finite.
\\

\begin{x}{\small\bf EXAMPLE} \ 
Take $f(x) = \abs{x}$ $-$then
\[
\begin{cases}
\ (D^+ f) (0) \ = \ 1\\[8pt]
\ (D_+ f) (0) \ = \ 1
\end{cases}
\implies (D_r f)(0) \ = \ 1
\]
\\[-.75cm]

\noindent
and

\[
\begin{cases}
\ (D^- f) (0) \ = \ -1\\[8pt]
\ (D_- f) (0) \ = \ -1
\end{cases}
\implies (D_\ell f)(0) \ = \ -1.
\]
\\[-.25cm]

If $(D_r f) (x)$ and $(D_\ell f)(x)$ exist and are equal, then their common value is denoted by 
$f^\prime (x)$ and is called the \un{derivative} of $f$ at $x$, $f$ being \un{differentiable} at $x$ 
if $f^\prime (x)$ is finite.
\\[-.25cm]

[So the relations 
\[
\pm \infty 
\ \neq \ 
(D^+ f) (x) 
\ = \ 
(D_+ f) (x) 
\ = \ 
(D^- f) (x) 
\ = \ 
(D_- f) (x) 
\ \neq \ 
\pm \infty 
\]
are tantamount to the differentiability of $f$ at $x$.]
\\[-.25cm]
\end{x}


\begin{x}{\small\bf EXAMPLE} \ 
Take $f(x) = \ds\frac{1}{x}$ $(x \neq 0)$, $f(0) = 0$ $-$then 
\[
\begin{cases}
\  (D_r f)(0) = +\infty \\[4pt]
\ (D_\ell f)(0) =  +\infty
\end{cases}
.
\]
Therefore $f^\prime(0) = +\infty$ but $f$ is not differentiable at 0.
\\[-.25cm]
\end{x}

There is much that can be said about Dini derivatives but we shall limit ourselves to a few points that are relevant for the sequel.
\\

\begin{x}{\small\bf THEOREM} \ 
Let $f: [a,b] \ra \R$ $-$then for any real number $r$, each of the following sets is at most countable:

\[
\begin{matrix*}[l]
&\{x: (D_+ f)(x) \geq r \quad \text{and} \quad (D^- f)(x) < r \},
\\[11pt]
&\{x: (D_- f)(x) \geq r \quad \text{and} \quad (D^+f)(x) < r \},
\\[11pt]
&\{x: (D^+ f)(x) \leq r \quad \text{and} \quad (D_- f)(x) > r \},
\\[11pt]
&\{x: (D^- f)(x) \leq r \quad \text{and} \quad (D_+f)(x) > r \}.
\end{matrix*}
\]
\\[-.75cm]
\end{x}

\begin{x}{\small\bf APPLICATION} \ 
Let $f: [a,b] \ra \R$ $-$then up to an at most countable set, 
\[
\begin{cases}
\  (D^+ f)(x) \geq (D_- f)(x) \\[4pt]
\ (D^- f)(x) \geq  (D_+ f)(x)
\end{cases}
.
\]
\\[-.75cm]
\end{x}

\begin{x}{\small\bf THEOREM} \ 
Let $f: [a,b] \ra \R$ be a Lebesgue measurable  function $-$then its 
Dini derivatives are Lebesgue measurable  functions.
\\[-.25cm]
\end{x}

To fix the ideas, let us consider a special case.  
So suppose that $f: [a,b] \ra \R$ is a Lebesgue measurable  function and $E \subset [a,b[$ is a Lebesgue measurable subset of $[a,b]$.  
Assume: \ $D_r \hsy f$ exists on $E$ $-$then $D_r \hsy f$ is a Lebesgue measurable  function on $E$.
\\[-.25cm]

To establish this, extend the definition of $f$ to $\R$ by setting $f = 0$ in $\R - [a,b]$.
Define a sequence $g_1, g_2, \ldots$ of Lebesgue measurable  functions via the prescription 
\[
g_n(x) 
\ = \ 
n \hsy (f \big(x + \frac{1}{n}\big) - f(x)).
\]
Let $D_e$ be the subset of $\R$ comprised of those $x$ such that $\lim\limits_{n \ra \infty} \ g_n(x)$ exists in $[-\infty, +\infty]$
$-$then $D_e$ is a Lebesgue measurable  set and 
\[
\lim\limits_{n \ra \infty} \ g_n \hsx : \hsx D_e \ra [-\infty, +\infty] 
\]
is a Lebesgue measurable  function.  
Take now an $x \in \ E$ and write
\begin{align*}
(D_r f ) (x) \ 
&=\ 
\lim\limits_{h \downarrow 0} \ \frac{f(x + h) - f(x)}{h}
\\[15pt]
&= \ 
\lim\limits_{n \ra \infty} \ \frac{f \big(x + \frac{1}{n}\big) - f(x)}{\frac{1}{n}}
\\[15pt]
&= \ 
\lim\limits_{n \ra \infty} \ g_n(x).
\end{align*}
Consequently $E \subset D_e$ and
\[
D_r f 
\ = \ 
\lim\limits_{n \ra \infty} \ g_n
\]
in $E$, hence $D_r f$ is a Lebesgue measurable  function on $E$.
\\

\begin{x}{\small\bf \un{N.B.}} \ 
Analogous considerations apply to $D_\ell \hsy f$ and $f^\prime$.
\end{x}

\chapter{
$\boldsymbol{\S}$\textbf{10}.\quad  DIFFERENTIATION}
\setlength\parindent{2em}
\setcounter{theoremn}{0}
\renewcommand{\thepage}{\S10-\arabic{page}}

\qquad
We shall first review some fundamental points.
\\[-.25cm] 

\begin{x}{\small\bf FACT} \ 
Let $f:[a,b] \ra \R$ be an increasing function $-$then $f$ is differentiable in $]a,b[ \hsx - E$, 
where $E$ is a set of Lebesgue measure 0 contained in $]a,b[$.  
\\[-.25cm] 

\un{Note}: \ 
Bear in mind that ``differentiable'' means that at $x \in \ ]a,b[ \hsx - E$, $f^\prime(x)$ exists and is finite.  
Moreover $f^\prime(x) = +\infty$ is possible only on a set of Lebesgue measure 0.]
\\[-.25cm] 
\end{x}

\begin{x}{\small\bf \un{N.B.}} \ 
\[
f^\prime:[a,b] - E \ra \R_{\geq \hsy 0} 
\]
is a Lebesgue measurable function.
\\[-.25cm] 
\end{x}

\begin{x}{\small\bf REMARK} \ 
If $E \subset \hsx ]a,b[$ is a set of Lebesgue measure 0, then it can be shown that there exists a continuous increasing 
function $f$ which is not differentiable at any point of $E$.
\\[-.25cm] 
\end{x}

\begin{x}{\small\bf RAPPEL} \ 
If $\phi$ is a Lebesgue measurable function and if $\psi = \phi$ almost everywhere, then $\psi$ is a Lebesgue measurable function. 
\\[-.25cm] 
\end{x}

\begin{x}{\small\bf FACT} \ 
Let $f:[a,b] \ra \R$ be an increasing function $-$then $f^\prime$ is integrable on $[a,b]$ and 
\[
\int\limits_a^b \ f^\prime 
\ \leq \ 
f(b) - f(a).
\]

[Note: \ This estimate can be sharpened to 
\[
\int\limits_a^b \ f^\prime 
\ \leq \ 
f(b-) - f(a+). ] 
\]
\\[-1cm]
\end{x}

\begin{x}{\small\bf EXAMPLE} \ 
One can construct a function $f:[a,b] \ra \R$ that is continuous and strictly increasing in $[a,b]$ such that 
$f^\prime = 0$ almost everywhere, hence
\[
0
\ = \ 
\int\limits_a^b \ f^\prime 
\ < \ 
f(b) - f(a).
\]
\\[-1cm]

Such functions are referred to as \un{totally singular}. 
Denjoy\footnote[2]{\vspace{.11 cm}Sur une fonction r\'eelle de Minkowski, \textit{J. Math. Pures Appl.} \textbf{17} (1915), 204-209.}
gave the first example. 
\\[-.5cm]

The existence of Lebesgue type singular functions (cf. \S12 \#5) (not stricly monotonic) are much simpler (intuitive) to construct.  
The existence of strictly increasing singular functions lies deeper and such functions are more difficult to construct and the proofs that they are singular 
and strictly increasing more involved.   
\\[-.5cm]

[Note: \ 
The above inequality also holds for the Lebesgue  ingular function.]
\\[-.5cm]

A famous example, the ?$(x)$ function, was given by 
Minkowski\footnote[3]{\vspace{.11 cm} “Zur Geometrie der Zahlen,”\textit{Gesammelte Abhandlungen} \textbf{2} (1911), 50–51.}.  
Minkowski's objective was to establish a one-one correspondence between the
rational numbers of (0, 1) and the quadratic irrationals of (0, 1).   
It was 
Denjoy\footnote[4]{\vspace{.11 cm}''Sur une fonction r´eelle de Minkowski,'' \textit{J. Math. Pures Appl.} \textbf{17} (1938), 105–151.}
who first proved that ?$(x)$ is totally singular.  
A survey article by
Salem\footnote[5]{\vspace{.11 cm}"On some singular monotonic functions which are strictly increasing," 
\textit{Proc. Amer. Math. Soc.} \textbf{53} (1943), 427–439.}
discusses singular monotonic functions more generally.
\\[-.25cm]
\end{x}

\begin{x}{\small\bf FACT} \ 
Given an $f \in \ \Lp^1[a,b]$, put
\[
F(x) 
\ = \ 
\int\limits_a^x \ f \hps (a \leq x \leq b).
\]
Then $F \in \ \AC[a,b]$ and $F^\prime = f$ almost everywhere.
\\[-.25cm] 
\end{x}

\begin{x}{\small\bf FACT} \ 
Suppose that $f:[a,b] \ra \R$ is absolutely continuous $-$then
\[
f(x) 
\ = \ 
f(a) + \int\limits_a^x \ f^\prime 
\hps (a \leq x \leq b).
\]
\\[-1cm]
\end{x}

\begin{x}{\small\bf FUBINI'S LEMMA} \ 
Let $\{f_n\}$ $(n = 1, 2, \ldots)$ be a sequence of increasing functions in $[a,b]$.  
Assume that the series
\[
\sum\limits_{n = 1}^\infty \ f_n(x)
\]
converges pointwise in $[a,b]$ to a function $F$ $-$then $F$ is differentiable almost everywhere in $[a,b]$ and 
\[
F^\prime (x) 
\ = \ 
\sum\limits_{n = 1}^\infty \ f_n^\prime(x)
\]
off of a set of Lebesgue measure 0.
\\[-.25cm]

PROOF \ 
Without loss of generality, take $f_k(a) = 0$ for all $k$ and observing that $F$ is increasing, 
let $E$ be the set of points $x \in \ ]a,b[$ such that the derivatives 
$F^\prime(x)$, $f_1^\prime(x)$, $f_2^\prime(x), \ldots$ all exist and are finite $-$then $[a,b] - E$ has Lebesgue measure 0.  
Let
\[
F_n (x) 
\ = \ 
\sum\limits_{k = 1}^n \ f_k(x).
\]
Suppose that $x \in \ E$ and $h$ is chosen small enough to ensure that $x + h \in \ [a,b]$ $-$then
\[
\frac{F(x + h) - F(x)}{h} 
\ = \ 
\sum\limits_{k = 1}^\infty \ \frac{f_k(x + h) - f_k(x)}{h} 
\]

$\implies$
\[
\frac{F(x + h) - F(x)}{h} 
\ \geq \ 
\sum\limits_{k = 1}^n \ \frac{f_k(x + h) - f_k(x)}{h} 
\]

$\implies$
\allowdisplaybreaks
\begin{align*}
F^\prime (x) \ 
&\geq \ 
\sum\limits_{k = 1}^n \ f_k^\prime(x)
\\[15pt]
&= \ 
F_n^\prime(x).
\end{align*}
The $f_k^\prime$ are nonnegative and the sequence
\[
\{F_n^\prime (x)\} \hps (n = 1, 2, \ldots)
\]
is bounded above by $F^\prime(x)$, hence is convergent.  It remains to establish that
\[
\lim\limits_{n \ra \infty} \ F_n^\prime 
\ = \ 
F^\prime
\]
almost everywhere in $[a,b]$.  Since 
\[
\lim\limits_{n \ra \infty} \ F_n (b)
\ = \ 
F (b), 
\]
there esists a subsequence $\{F_{n_j}(b)\}$ such that 
\allowdisplaybreaks
\begin{align*}
F(a) - F_{n_j} (a) \ 
&=\ 
0
\\[8pt]
&\leq\ 
F(b) - F_{n_j} (b) 
\\[8pt]
&\leq \ 
2^{-j}.
\end{align*}
But $F - F_{n_j}$ is an increasing function, thus
\allowdisplaybreaks
\begin{align*}
0 \ 
&\leq\ 
F(x) - F_{n_j} (x) 
\\[8pt]
&\leq \ 
2^{-j}
\end{align*}
for all $x \in \ [a,b]$ and so the series
\[
\sum\limits_{j = 1}^\infty \ (F^\prime - F_{n_j}^\prime)
\]
is a pointwise convergent series of increasing functions.  
Reasoning as above, we conclude that the series
\[
\sum\limits_{j = 1}^\infty \ (F^\prime - F_{n_j}^\prime)
\]
is convergent almost everywhere in $[a,b]$ and from this it follows that
\[
F^\prime (x) - F_n^\prime (x) \hsx \ra \hsx 0
\]
as $n \ra \infty$ for almost all $x \in \ [a,b]$.
\\[-.25cm]
\end{x}

\begin{x}{\small\bf APPLICATION} \ 
Suppose that $f:[a,b] \ra \R$ is increasing and let $s_f:[a,b] \ra \R$ be the saltus function attached to $f$ 
$-$then $s_f^\prime = 0$ almost everywhere.
\\[-.25cm]

[In general, $s_f$ is not continuous.  
Still, a \un{continuous singular function} is a continuous function whose derivative exists and is zero almost everywhere.  
To illustrate, write 
\begin{align*}
f \ 
&=\ 
(f - s_f) + s_f 
\\
&=\ 
r_f + s_f,
\end{align*}
where by construction $r_f$ is increasing and continuous.  
And almost everywhere
\begin{align*}
f^\prime  \ 
&=\ 
r_f^\prime + s_f^\prime
\\
&=\ 
r_f^\prime.
\end{align*}
Introduce $F$ by the rule 
\[
F(x) 
\ = \ 
\int\limits_a^x \ f^\prime
\]
and set 
\[
f_\text{cs} 
\ = \ 
r_f - F.
\]
Then almost everywhere
\begin{align*}
f_\text{cs}^\prime 
&=\
r_f^\prime - F^\prime 
\\
&=\
f ^\prime- f^\prime 
\\
&=\
0.
\end{align*}
Therefore $f_\text{cs}$ is a continuous singular increasing function and 
\begin{align*}
f \ 
&=\ 
r_f + s_f
\\
&=\ 
F + f_\text{cs} + s_f.]
\end{align*}

The fact that an $f \in \ \BV[a,b]$ can be represented as the difference of two increasing functions implies that $f$ 
is differentiable almost everywhere. 
\\[-.25cm]

[Note: \ 
Therefore a continuous nowhere differentiable function is not of bounded variation.]
\\[-.25cm]
\end{x}

\begin{x}{\small\bf THEOREM} \ 
Suppose that $f \in \ \BV[a,b]$ $-$then for almost all $x \in \ [a,b]$, 
\[ 
\abs{f^\prime (x)}
\ = \ 
\tT_f^\prime [a,x].
\]

PROOF \ 
Given $n \in \ \N$, choose a partition $P_n \in \ \sP[a,b]$ such that 
\[
\sum\limits_k \ 
\abs{f(x_k) - f(x_{k -1})} 
\ > \ 
\tT_f[a,b] - 2^{-n}.
\]
In the segment $x_{k - 1} \leq x \leq x_k$ of $P_n$, let
\[
\begin{cases}
\ f_n(x) = f(x) + c_n^+ \hspace{0.7cm} \text{if $f(x_k) - f(x_{k -1}) \geq 0$}\\
\hps \text{or} \\
\ f_n(x) = -f(x) + c_n^- \hspace{0.7cm} \text{if $f(x_k) - f(x_{k -1}) \leq 0$}
\end{cases}
,
\]
where the constants are chosen so that $f_n(a) = 0$ and the values of $f_n$ at $x_k$ agree $-$then
\[
f_n(x_k) - f_n(x_{k-1})
\ = \ 
\abs{f_n(x_k) - f_n(x_{k-1})},
\]
so 
\allowdisplaybreaks
\begin{align*}
\tT_f[a,b]  - f_n(b) \ 
&= \ 
\tT_f[a,b]  - \sum\limits_k \ (f_n(x_k) - f_n(x_{k-1}))
\\[15pt]
&=\ 
\tT_f[a,b]  - \sum\limits_k \ \abs{f(x_k) - f(x_{k-1})}
\\[15pt]
&\leq \ 
2^{-n}.
\end{align*}
On the other hand, the function 
\[
x \ra \tT_f[a,x] - f_n(x)
\]
is increasing, hence
\allowdisplaybreaks
\begin{align*}
\tT_f[a,x]  - f_n(x) \ 
&\leq \ 
\tT_f[a,b]  - f_n(b)
\\[15pt] 
&\leq \ 
2^{-n}
\end{align*}

$\implies$
\begin{align*}
\sum\limits_{n = 1}^\infty \ (\tT_f[a,x]  - f_n(x)) \ 
&\leq \ 
\sum\limits_{n = 1}^\infty \ 2^{-n}
\\[15pt] 
&< \ 
+\infty.
\end{align*}
The series
\[
\sum\limits_{n = 1}^\infty \ (\tT_f[a,x]  - f_n(x))
\]
is therefore pointwise convergent, thus by Fubini's lemma, the derived series converges almost everywhere, thus
\[
\tT_f^\prime [a,x] - f_n^\prime (x)  \lra 0
\]
almost everywhere.  But
\[
f_n^\prime (x) 
\ = \ 
\pm f^\prime (x).
\]
Since $\tT_f^\prime[a,x] \geq 0$ ($\tT_f[a,x]$ being increasing), the upshot is that
\[
\abs{f^\prime (x) }
\ = \ 
\tT_f^\prime [a,x]
\]
almost everywhere. 
\\[-.25cm]
\end{x}

\begin{x}{\small\bf APPLICATION} \ 
\[
f \in \ \BV[a,b] 
\implies 
f^\prime \in \ \Lp^1[a,b].
\]

[For 
\allowdisplaybreaks
\begin{align*}
\int\limits_a^b \ \abs{f^\prime} \ 
&=\ 
\int\limits_a^b \ \tT_f^\prime [a,-]
\\[15pt]
&\leq \ 
\tT_f[a,b] - \tT_f[a,a]
\\[15pt]
&= \ 
\tT_f[a,b]
\\[15pt]
&< \ 
+\infty.]
\end{align*}
\\[-1cm]
\end{x}

\begin{x}{\small\bf THEOREM} \ 
Given $f \in \ \Lp^1[a,b]$, put
\[
F(x) 
\ = \ 
\int\limits_a^x \ f.
\]
Then 
\[
\tT_F[a,b] 
\ = \ 
\norm{f}_{\Lp^1}.
\]

PROOF \ 
Given a $P \in \ \sP[a,b]$, 

\allowdisplaybreaks
\begin{align*}
\sum\limits_{k = 1}^n \ \abs{F(x_k) - F(x_{k - 1})} \ 
&= \ 
\sum\limits_{k = 1}^n \
\bigg| \hsx \int\limits_{x_{k - 1}}^{x_k} \ f \ \bigg|
\\[15pt]
&\leq \ 
\int\limits_a^b \ \abs{f} 
\\[15pt]
&< \ 
+\infty
\end{align*}

$\implies$
\[
\tT_f[a,b] 
\ \leq \ 
\norm{f}_{\Lp^1}.
\]
To reverse this, recall that $F \in \ \AC[a,b]$, that $F^\prime = f$ almost everywhere, and that 
\[
\abs{F^\prime}
\ = \ 
\tT_F^\prime [a,-]
\]
almost everywhere.  
Therefore
\allowdisplaybreaks
\begin{align*}
\norm{f}_{\Lp^1} \ 
&=\ 
\int\limits_a^b \ \abs{F^\prime}
\\[15pt]
&=\ 
\int\limits_a^b \ \tT_F^\prime [a,-]
\\[15pt]
&\leq \ 
\tT_F[a,b] - \tT_F[a,a]
\\[15pt]
&=\ 
\tT_F[a,b].
\end{align*}
\\[-1cm]
\end{x}

\begin{x}{\small\bf LEMMA} \ 
Suppose that $f:[a,b] \ra \R$ is increasing $-$then $f  \in \ \AC[a,b]$ iff
\[
\int\limits_a^b \ f^\prime
\ = \ 
fb) - f(a).
\]

PROOF \ 
If $f  \in \ \AC[a,b]$, then
\[
f(x) 
\ = \ 
f(a) + \int\limits_a^x \ f^\prime
\hps (a \leq x \leq b)
\]

$\implies$
\[
f(b) - f(a) 
\ = \ 
\int\limits_a^b \ f^\prime.
\]
Conversely, write
\[
f(x) 
\ = \ 
\int\limits_a^x \ f^\prime + f_\text{cs}(x) + s_f(x).
\]
Then
\[
f(x) 
\ = \ 
f(a) + \int\limits_a^x \ f^\prime + g(x), 
\]
where
\allowdisplaybreaks
\begin{align*}
&\hspace{1.5cm}
 f_\text{cs}(x) + s_f(x) \ = \ f(a) + g(x), 
\\
&\implies
\\
&\hspace{1.5cm}
 f_\text{cs}(a) + s_f(a) \ = \ f(a) + g(a) 
\\
&\implies
\\
&\hspace{1.5cm}
r_f(a) - F(a) + s_f(a)  \ = \ f(a) + g(a)
\\
&\implies
\\
&\hspace{1.5cm}
r_f(a) + s_f(a)  \ = \ f(a) + g(a)
\\
&\implies
\\
&\hspace{1.5cm}
(f - s_f)(a) + s_f(a) \ = \ f(a) + g(a)
\\
&\implies
\\
&\hspace{1.5cm}
f(a) = f(a) + g(a)
\\
&\implies
\\
&\hspace{1.5cm}
g(a) = 0.
\end{align*}
In addition, the assumption that

\[
\int\limits_a^b \ f^\prime
\ = \ 
f(b) - f(a)
\]
implies that
\allowdisplaybreaks
\begin{align*}
g(b) \ 
&=\ 
f(b) - f(a) - \int\limits_a^b \ f^\prime
\\[15pt]
&= \ 
0.
\end{align*}
Since $g$ is increasing, it follows that $g(x) = 0$ for  all $x \in \ [a,b]$, hence
\[
f(x) 
\ = \ 
f(a) + \int\limits_a^x \ f^\prime.
\]
\\[-1cm]
\end{x}

\begin{x}{\small\bf THEOREM} \ 
Suppose that $f \in \ \BV[a,b]$ $-$then $f \in \ \AC[a,b]$ iff 
\[
\tT_f[a, b] 
\ = \ 
\int\limits_a^b \ \abs{f^\prime}.
\]

PROOF \ 
On the one hand, 
\begin{align*}
f \in \ \AC[a,b] 
&\implies
f^\prime \in \ \Lp^1[a,b]
\\[15pt]
&\implies
\tT_f[a, b] = \int\limits_a^b \ \abs{f^\prime}.
\end{align*}
On the other hand, assume the stated relation.  
Since for almost all $x$ in $[a,b]$, 
\[
\abs{f^\prime(x)}
\ = \ 
\tT_f^\prime[a, x],
\]
we have
\[
\tT_f[a, b]
\ = \ 
\int\limits_a^b \ \tT_f^\prime[a, -]
\]
or still, 
\[
\tT_f[a, b] -  \tT_f[a, a] 
\ = \ 
\int\limits_a^b \ \tT_f^\prime[a, -].
\]
But $\tT_f[a, -]$ is increasing, thus in view of the lemma, $\tT_f[a, -]$ is absolutely continuous, 
which in turn implies that $f$ is absolutely continuous.
\end{x}

\chapter{
$\boldsymbol{\S}$\textbf{11}.\quad  ESTIMATE OF THE IMAGE}
\setlength\parindent{2em}
\setcounter{theoremn}{0}
\renewcommand{\thepage}{\S11-\arabic{page}}

\begin{x}{\small\bf RAPPEL} \ 
\[
\begin{cases}
\ \lambda \ = \ \text{Lebesgue measure}\\
\ \lambda^* \ = \ \text{outer Lebesgue measure}
\end{cases}
.
\]
\\[-1cm]
\end{x}

\begin{x}{\small\bf LEMMA} \ 
Let $f:[a,b] \ra \R$.  
Suppose that $E \subset [a,b]$ is a subset in which $f^\prime$ exists, subject to $\abs{f^\prime} \leq K$ $-$then
\[
\lambda^*(f(E))
\ \leq \ 
K \hsx \lambda^*(E).
\]

The proof will be carried out in seven steps.
\\[-.25cm]

\un{Step 1:} \ 
Given $x \in E$, $f^\prime(x)$ exists and
\[
\abs{f^\prime(x)} 
\ = \ 
\bigg|
\lim\limits_{y \ra x} \ \frac{f(y) - f(x)}{y - x} 
\bigg|
\ \leq \ 
K.
\]
So, $\forall \ x \in E$, $\exists \ \delta > 0$: 
\[
\abs{f(y) - f(x)} 
\ \leq \ 
K \hsx \abs{y - x} 
\hps (y \in \hsx ]x - \delta, x + \delta[ \ \cap  \ [a,b]).
\]
If now for $n = 1, 2, \ldots, $
\[
E_n 
\ = \ 
\bigg\{
x \in E \hsx : \hsx \abs{f(y) - f(x)}  \leq K \hsx \abs{y - x}  
\hps (y \in \hsx ]x - \frac{1}{n}, x + \frac{1}{n}[ \hsx)
\bigg\},
\]
then each $x \in E$ belongs to $E_n$ $(n \gg 0)$, hence 
\[
E 
\ \subset \ 
\bigcup\limits_{n = 1}^\infty \ E_n.
\]

\noindent
On the other hand, $\forall \ n$, $E_n \subset E$ and $\{E_n\}$ is increasing.  
Therefore

\[
E 
\ = \ 
\bigcup\limits_{n = 1}^\infty \ E_n
\ = \ 
\lim\limits_{n \ra \infty} \ E_n.
\]

\un{Step 2:} \ Consequently 
\[
\lim\limits_{n \ra \infty} \ \lambda^* (E_n) 
\ = \ 
\lambda^* (E).
\]
But
\begin{align*}
f(E) \ 
&= \ 
f \big(  \bigcup\limits_{n = 1}^\infty \ E_n \big)
\\[15pt]
&= \ 
\bigcup\limits_{n = 1}^\infty \ f(E_n)
\\[15pt]
&= \ 
\lim\limits_{n \ra \infty} \ f(E_n)
\end{align*}

$\implies$
\[
\lim\limits_{n \ra \infty} \ \lambda^*(f(E_n))
\ = \ 
\lambda^*(f(E)).
\]
\\[-1cm]

\un{Step 3:} \ 
Let $\varepsilon > 0$ be given and let $I_{n, k}$ $(k = 1, 2, \ldots)$ be a sequence of open intervals such that 

\[
\begin{cases}
\ \ds\lambda(I_{n, k}) \ < \ \frac{1}{n} \\[11pt]
\ E_n \subset \bigcup\limits_{k = 1}^\infty \ I_{n, k}
\end{cases}
,
\]
and
\[
\sum\limits_{k = 1}^\infty \ \lambda(I_{n, k}) 
\ \leq \ 
\lambda^*(E_n) + \varepsilon.
\]
\\[-1cm]

\un{Step 4:} \ 
\[
E_n 
\ = \ 
\bigcup\limits_{k = 1}^\infty \ (E_n \cap I_{n, k})
\]
and
\[
f(E_n) 
\ = \ 
\bigcup\limits_{k = 1}^\infty \ f(E_n \cap I_{n, k}).
\]
\\[-1cm]

\un{Step 5:} \ 
If $x_1$, $x_n \in E_n \cap I_{n, k}$, then 
\[
\abs{f(x_1) - f(x_2)} 
\ \leq \ 
K \hsx \abs{x_1 - x_2} 
\ \leq \ 
K \hsx \lambda (I_{n,k})
\]

$\implies$
\[
\lambda^*(f(E_n \cap I_{n, k}))
\ \leq \ 
K \hsx \lambda (I_{n,k}).
\]
\\[-1cm]

\un{Step 6:} \ 
\begin{align*}
\lambda^*(f(E_n)) \ 
&=\ 
\lambda^* \bigg( \bigcup\limits_{k = 1}^\infty \ f(E_n \hsx \cap \hsx I_{n, k})\bigg) 
\\[15pt]
&\leq \ 
\sum\limits_{k = 1}^\infty  \ \lambda^* \big( f(E_n \hsx \cap \hsx I_{n, k}) \big)
\\[15pt]
&\leq \ 
\sum\limits_{k = 1}^\infty  \ K \hsx \lambda(I_{n, k})
\\[15pt]
&\leq \
K ( \lambda^* (E_n) + \varepsilon).
\end{align*}
\\[-1cm]

\un{Step 7:} \ 
\allowdisplaybreaks
\begin{align*}
\lambda^*(f(E)) \ 
&=\ 
\lim\limits_{n \ra \infty} \ \lambda^*(f(E_n))
\\[15pt]
&\leq \
K (\lim\limits_{n \ra \infty} \  \lambda^* (E_n) + \varepsilon)
\\[15pt]
&= \
K ( \lambda^* (E) + \varepsilon)
\end{align*}

$\implies$
\[
\lambda^*(f(E))
\ \leq \ 
K \hsx \lambda^* (E) 
\hps (\varepsilon \downarrow 0), 
\]
the assertion of the lemma.
\end{x}

\begin{x}{\small\bf THEOREM} \ 
Let $f: [a,b] \ra \R$ be Lebesgue measurable.  
Suppose that $E \subset [a,b]$ is a Lebesgue measurable subset in which $f$ is differentiable $-$then
\[
\lambda^*(f(E))
\ \leq \ 
\int\limits_E \ \abs{f^\prime (x)}.
\]

PROOF \ 
Note that $f^\prime : E \ra \R$ is a Lebesgue measurable function.  
This said, to begin with, assume that in $E$, $\abs{f^\prime} < M$ (a positive integer).  
Let
\[
E_k^n 
\ = \ 
\bigg\{x \in  E \hsx : \hsx \frac{k-1}{2^n} \hsx \leq \hsx \abs{f^\prime (x)} \hsx < \hsx \frac{k}{2^n}  \bigg\},
\]
where
\[
k = 1, 2, \ldots, M \hsy 2^n, 
\quad 
n = 1, 2, \ldots \hsx .
\]
Then for each $n$, 
\allowdisplaybreaks
\begin{align*}
\lambda^* (f(E)) \ 
&=\ 
\lambda^* 
\big(
f
\big( \ 
\bigcup\limits_k \ 
E_k^n
\big)
\big)
\\[15pt]
&=\ 
\lambda^* 
\big(
\bigcup\limits_k \ 
f(E_k^n)
\big)
\\[15pt]
&\leq\ 
\sum\limits_k \ 
\lambda^* (f(E_k^n)) 
\\[15pt]
&\leq\ 
\sum\limits_k \ \frac{k}{2^n}  \ 
\lambda(E_k^n)
\\[15pt]
&=\ 
\sum\limits_k \ \frac{k - 1}{2^n} \ \lambda(E_k^n)
\hsx + \hsx 
\frac{1}{2^n} \ \sum\limits_k \ \lambda(E_k^n).
\end{align*}
Therefore
\allowdisplaybreaks
\begin{align*}
\lambda^*(f(E)) \ 
&\leq \ 
\lim\limits_{n \ra \infty} \ 
\big(
\sum\limits_k \ \frac{k - 1}{2^n} \ \lambda(E_k^n)
\hsx + \hsx 
\frac{1}{2^n} \ \sum\limits_k \ \lambda(E_k^n)
\big)
\\[15pt]
&=\ 
\int\limits_E \ \abs{f^\prime}.
\end{align*}
To treat the case of unbounded $f^\prime$, let  
\[
A_k 
\ = \ 
\{x \in  E \hsx : \hsx  k - 1 \hsx \leq \hsx \abs{f^\prime (x)} \hsx < \hsx k\} 
\hps (k = 1, 2, \ldots).
\]
Then
\allowdisplaybreaks
\begin{align*}
\lambda^*(f(E)) \ 
&=\ 
\lambda^* 
\big( 
f
\big( 
\bigcup\limits_k \ A_k
\big)
\big)
\\[15pt]
&\leq \
\lambda^*
\big( 
\bigcup\limits_k \ f(A_k)
\big)
\\[15pt]
&\leq \
\sum\limits_k \ \lambda^*( f(A_k))
\\[15pt]
&\leq \
\sum\limits_k \int\limits_{A_k} \ \abs{f^\prime}
\\[15pt]
&= \
\int\limits_E \ \abs{f^\prime}.
\end{align*}

[Note: \ In point of fact, $f(E)$ is Lebesgue measurable, so
\[
\lambda^*(f(E)) 
\ = \  
\lambda(f(E)) .]
\]
\end{x}

\begin{x}{\small\bf \un{N.B.}} \ 
It follows that
\[
\lambda^*(f(E)) 
\ = \  
0
\]
\vspace{-1.25cm}

\noindent
if $f^\prime = 0$.
\\[-.25cm]

[It can be shown conversely that
\[
\lambda^*(f(E)) 
\ = \  
0
\]
implies that $f^\prime = 0$ almost everywhere in $E$.]
\\[-.5cm]
\end{x}

\begin{x}{\small\bf SCHOLIUM} \ 
Suppose that $f$ has a finite derivative on a set $E$ $-$then 
$\lambda^*(f(E)) = 0$ iff $f^\prime = 0$ almost everywhere on $E$.
\end{x}

\chapter{
$\boldsymbol{\S}$\textbf{12}.\quad  ABSOLUTE CONTINUITY II}
\setlength\parindent{2em}
\setcounter{theoremn}{0}
\renewcommand{\thepage}{\S12-\arabic{page}}

\begin{x}{\small\bf THEOREM} \ 
If $f:[a,b] \ra \R$ is absolutely continuous and if $f^\prime (x) = 0$ almost everywhere, then $f$ is a constant function.
\\[-.5cm]

[Let
\[
E 
\ = \ 
\{x \in [a,b] \ : \ f^\prime(x) = 0\}
\]
and let
\[
E^\prime 
\ = \ 
[a,b] - E.
\]
The assumption that $f \in \AC [a,b]$ implies that $f$ has property (N) which in turn implies that $f$ sends 
Lebesgue measurable sets to Lebesgue measurable sets.  
In particular: $f(E)$, $f(E^\prime)$ are Lebesgue measurable and 
\[
\lambda (f ([a,b]) 
\ \leq \ 
\lambda (f(E)) + \lambda (f(E^\prime)).
\]
So first
\begin{align*}
\lambda (f(E) \ 
&\leq \ 
0  \hsx \lambda (E) 
\qquad \text{(\S11 \#2  (``K = 0''),  $E$ measurable)}
\\
&=\ 
0.
\end{align*}
And second, $E^\prime$ is a set of Lebesgue measure 0, hence the same is true of $f(E^\prime)$.
\\[-.5cm]

\noindent All told then
\[
\lambda (f ([a,b]) 
\ = \ 
0 .
\]
Owing now to the continuity of $f$, the image $f([a,b])$ is a point or a closed interval.  
But the latter is a non-sequitur, thus $f([a,b])$ is a singleton.]
\\[-.25cm]
\end{x}

\begin{x}{\small\bf MAIN THEOREM} \ 
Let $f:[a,b] \ra \R$. $-$then $f$ is absolutely continuous iff 
the following four conditions are satisfied: 
\\[-.5cm]

\qquad (1) \quad $f$ is continuous.
\\[-.5cm]

\qquad (2) \quad $f^\prime$ exists almost everywhere.
\\[-.5cm]

\qquad (3) \quad $f^\prime \in \Lp^1 [a,b].$
\\[-.5cm]

\qquad (4) \quad $f$ has property (N).
\\[-.25cm]

PROOF  \ 
An absolutely continuous function has these properties.  
Conversely, assume that $f$ satisfies the stated conditions.   
Owing to (3), given $\varepsilon > 0$, there exists $\delta > 0$ such that 
\[
E \subset [a,b] \quad \& \quad \lambda(E) < \delta 
\implies 
\int\limits_E \ \abs{f^\prime} \ < \ \varepsilon.
\]
Fix
\[
a 
\ \leq \ 
a_1
\ < \ 
b_1
\ \leq \ 
a_2
\ < \ 
b_2
\ \leq \ 
\cdots
\ \leq \ 
a_n 
\ < \ 
b_n
\ \leq \ 
b
\]
with
\[
\sum\limits_{k = 1}^ n \  (b_k - a_k) 
 \ < \ \delta.
\]
Then
\[
\sum\limits_{k = 1}^ n \quad \int\limits_{[a_k, b_k]} \ \abs{f^\prime}
\ < \ 
\varepsilon.
\]
Let
\[
A_k 
\ = \ 
\{x \in [a_k, b_k] \ : \ f^\prime (x) \ \text{exists}\}.
\]
Thanks to (2), $[a_k, b_k] - A_k$ is a set of Lebesgue measure 0, hence thanks to (4), 
$f([a_k, b_k] - A_k)$ is a set of Lebesgue measure 0.  
Therefore

\vspace{-.65cm}
\allowdisplaybreaks
\begin{align*}
\sum\limits_{k = 1}^n \  \abs{f(b_k) - f(a_k)} \ 
&\leq \ 
\sum\limits_{k = 1}^n \  \lambda(f([a_k, b_k]))
\qquad \text{(by (1))}
\\[8pt]
&=\ 
\sum\limits_{k = 1}^n \ \lambda (f(A_k))
\\[8pt]
&\leq \  
\sum\limits_{k = 1}^n \ \int\limits_{A_k} \ \abs{f^\prime}
\\[8pt]
&= \  
\sum\limits_{k = 1}^n \ \int\limits_{[a_k, b_k]} \ \abs{f^\prime}
\\[8pt]
 &< \ 
\varepsilon.
\end{align*}
\\[-1.5cm]
\end{x}


\begin{x}{\small\bf SCHOLIUM} \ 
If $f \in \BV [a,b]$ is continuous and possesses property (N), then $f \in \AC [a,b]$.
\\[-.5cm]

[One has only to note that if $f$ is of bounded variation, then $f^\prime$ exists almost everywhere and 
$f^\prime \in \Lp^1[a,b]$.]
\\[-.25cm]
\end{x}

\begin{x}{\small\bf LEMMA} \ 
If $f:[a,b] \ra \R$ has a finite derivative at every point $x \in [a,b]$, then $f$ has property (N).
\\[-.5cm]

PROOF \ 
Suppose that $\lambda(E) = 0$ $(E \subset [a,b])$.  
For each positive integer $n$, let
\[
E_n 
\ = \ 
\{x \in E : \abs{f^\prime (x)} \leq n\}.
\]
Then $\lambda(E_n) = 0$ and 
\begin{align*}
\lambda^* (f(E_n)) \ 
&\leq \ 
n \hsy \lambda^* (f(E_n))
\\[11pt]
&=\ 
n \hsy \lambda(E_n) 
\\[11pt]
&=\ 
0
\end{align*}

$\implies$ 
\[
\lambda (f(E_n )) 
\ = \ 
0.
\]
Since
\[
E 
\ = \ 
\bigcup\limits_{n = 1}^\infty \ E_n
\]
and
\[
f(E) 
\ = \ 
f\left( \bigcup\limits_{n = 1}^\infty \ E_n \right) 
\ = \ 
\bigcup\limits_{n = 1}^\infty \ f(E_n),
\]
the conclusion is that

\allowdisplaybreaks
\begin{align*}
\lambda^* (f(E)) \ 
&\leq \ 
\sum\limits_{n = 1}^ \infty \ \lambda^* (f(E_n)) 
\\[15pt]
&=\ 
\sum\limits_{n = 1}^ \infty \ \lambda (f(E_n)) 
\\[15pt]
&=\ 
0.
\end{align*}
I.e.: \  $\lambda (f(E)) = 0$.
\\[-.25cm]
\end{x}

\begin{x}{\small\bf EXAMPLE} \ 
One can construct a continuous function $f:[a,b] \ra \R$ with a finite derivative almost everywhere which fails 
to have property (N).
\\[-.25cm]

In particular, one such example is the \un{Lebesgue (singular) function}.

\noindent
Let \CantorSet be the usual Cantor set: 
\\[-.5cm]

\qquad $P_0 = [0,1]$.  
\\[-.5cm]

\qquad
$P_1 \ = \ \bigg[0, \frac{1}{3}\big] \cup \bigg[\frac{2}{3}, 1\big]
\qquad \text{(deleting $\big]\frac{1}{3}, \frac{2}{3}\big[$ from $P_0$)}
$
\\[-.5cm]

\qquad
Etc.
\\[-.25cm]

\noindent
Get a sequence of compact sets
\[
P_1 \supset P_2 \supset \ldots 
\]
where
\[
P_n 
\ = \ 
J_{n, 1} \cup \cdots \cup J_{n, 2^n}
\quad
\text{and}  
\quad
\text{length} (J_{n,k}) 
\ = \ 
\frac{1}{3^n}.
\]
Then
\[
P 
\ \equiv \  
\bigcap \ P_n.]
\]
\\[-.5cm]
Let $\{I_{n, k}\}$ be the deleted open intervals associated with the construction of $\CantorSet.$  
\\[-.25cm]

\noindent
So
\[
[0,1] \ - \ \CantorSet 
\ = \ 
\bigcup\limits_{n = 1}^\infty \ 
\bigcup\limits_{k = 1}^{2 n - 1} \ 
I_{n, k}.
\]
\\[-.5cm]

\noindent
For each $n$ define
\\[-.25cm]

\[
\LebesgueFunction_n : [0, 1] \lra [0,1]
\qquad 
\begin{cases}
\ \text{linearly on each $J_{n, k}$ increasing by $\ds\frac{1}{2^n}$}
\\[8pt]
\text{constant on each $I_{n, k}$}
\end{cases}
.
\]

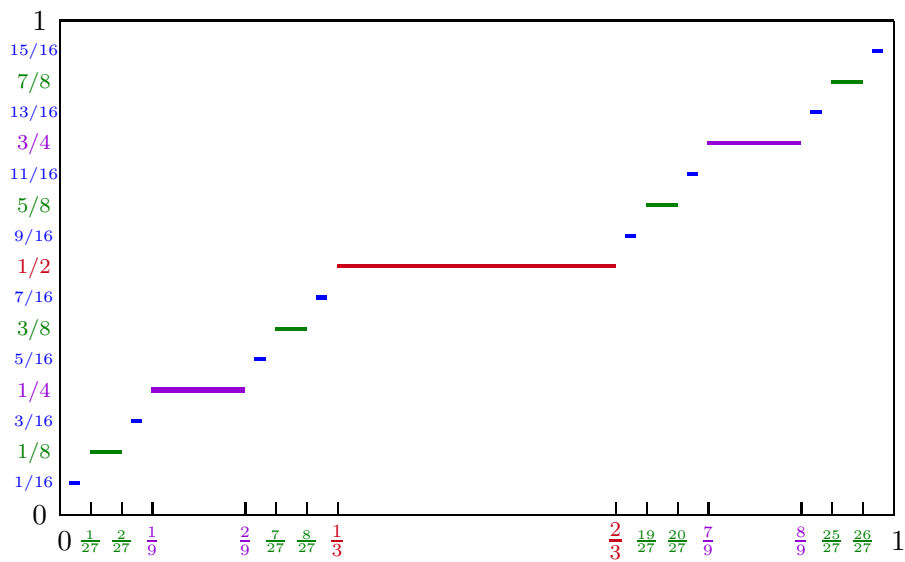
\setlength{\unitlength}{0.240900pt}
\ifx\plotpoint\undefined\newsavebox{\plotpoint}\fi
\sbox{\plotpoint}{\rule[-0.200pt]{0.400pt}{0.400pt}}%
\begin{picture}(1500,900)(0,0)
\sbox{\plotpoint}{\rule[-0.200pt]{0.400pt}{0.400pt}}%
\put(130.0,82.0){\rule[-0.200pt]{4.818pt}{0.400pt}}
\put(110,82){\makebox(0,0)[r]{ \text{\small$0$}}}

\put(1419.0,82.0){\rule[-0.200pt]{4.818pt}{0.400pt}}

\put(90,131){\color{blue}{\makebox(0,0){\text{\tiny$1/16$}}}}

\put(90,179){\color{officegreen}{\makebox(0,0){\text{\scriptsize$1/8$}}}}

\put(90,228){\color{blue}{\makebox(0,0){\text{\tiny$3/16$}}}}

\put(90,276){\color{darkviolet}{\makebox(0,0){\text{\scriptsize$1/4$}}}}

\put(90,325){\color{blue}{\makebox(0,0){\text{\tiny$5/16$}}}}

\put(90,373){\color{officegreen}{\makebox(0,0){\text{\scriptsize$3/8$}}}}

\put(90,422){\color{blue}{\makebox(0,0){\text{\tiny$7/16$}}}}

\put(90,471){\color{harvardcrimson}{\makebox(0,0){\text{\scriptsize$1/2$}}}}

\put(90,519){\color{blue}{\makebox(0,0){\text{\tiny$9/16$}}}}

\put(90,567){\color{officegreen}{\makebox(0,0){\text{\scriptsize$5/8$}}}}

\put(90,616){\color{blue}{\makebox(0,0){\text{\tiny$11/16$}}}}

\put(90,665){\color{darkviolet}{\makebox(0,0){\text{\scriptsize$3/4$}}}}

\put(90,713){\color{blue}{\makebox(0,0){\text{\tiny$13/16$}}}}

\put(90,761){\color{officegreen}{\makebox(0,0){\text{\scriptsize$7/8$}}}}

\put(90,810){\color{blue}{\makebox(0,0){\text{\tiny$15/16$}}}}

\put(130.0,859.0){\rule[-0.200pt]{4.818pt}{0.400pt}}
\put(110,859){\makebox(0,0)[r]{ \text{\small$1$}}}

\put(178.0,82.0){\rule[-0.200pt]{0.400pt}{4.818pt}}
\put(178,41){\color{officegreen}{\makebox(0,0){\text{\tiny$\frac{1}{27}$}}}}

\put(227.0,82.0){\rule[-0.200pt]{0.400pt}{4.818pt}}
\put(227,41){\color{officegreen}{\makebox(0,0){\text{\tiny$\frac{2}{27}$}}}}

\put(275.0,82.0){\rule[-0.200pt]{0.400pt}{4.818pt}}
\put(275,41){\color{darkviolet}{\makebox(0,0){\text{\scriptsize$\frac{1}{9}$}}}}

\put(421.0,82.0){\rule[-0.200pt]{0.400pt}{4.818pt}}
\put(421,41){\color{darkviolet}{\makebox(0,0){\text{\scriptsize$\frac{2}{9}$}}}}

\put(469.0,82.0){\rule[-0.200pt]{0.400pt}{4.818pt}}
\put(469,41){\color{officegreen}{\makebox(0,0){\text{\tiny$\frac{7}{27}$}}}}

\put(518.0,82.0){\rule[-0.200pt]{0.400pt}{4.818pt}}
\put(518,41){\color{officegreen}{\makebox(0,0){\text{\tiny$\frac{8}{27}$}}}}

\put(566.0,82.0){\rule[-0.200pt]{0.400pt}{4.818pt}}
\put(566,41){\color{harvardcrimson}{\makebox(0,0){\text{\footnotesize$\frac{1}{3}$}}}}

\put(1003.0,82.0){\rule[-0.200pt]{0.400pt}{4.818pt}}
\put(1003,41){\color{harvardcrimson}{\makebox(0,0){\text{$\frac{2}{3}$}}}}

\put(1051.0,82.0){\rule[-0.200pt]{0.400pt}{4.818pt}}
\put(1051,41){\color{officegreen}{\makebox(0,0){\text{\tiny$\frac{19}{27}$}}}}

\put(1100.0,82.0){\rule[-0.200pt]{0.400pt}{4.818pt}}
\put(1100,41){\color{officegreen}{\makebox(0,0){\text{\tiny$\frac{20}{27}$}}}}

\put(1148.0,82.0){\rule[-0.200pt]{0.400pt}{4.818pt}}
\put(1148,41){\color{darkviolet}{\makebox(0,0){\text{\scriptsize$\frac{7}{9}$}}}}

\put(1294.0,82.0){\rule[-0.200pt]{0.400pt}{4.818pt}}
\put(1294,41){\color{darkviolet}{\makebox(0,0){\text{\scriptsize$\frac{8}{9}$}}}}

\put(1342.0,82.0){\rule[-0.200pt]{0.400pt}{4.818pt}}
\put(1342,41){\color{officegreen}{\makebox(0,0){\text{\tiny$\frac{25}{27}$}}}}

\put(1391.0,82.0){\rule[-0.200pt]{0.400pt}{4.818pt}}
\put(1391,41){\color{officegreen}{\makebox(0,0){\text{\tiny$\frac{26}{27}$}}}}

\put(1439.0,82.0){\rule[-0.200pt]{0.400pt}{4.818pt}}
\put(1439,41){\makebox(0,0){ \text{\small$1$}}}

\put(130,41){\makebox(0,0){ \text{\small$0$}}}

\put(130.0,82.0){\rule[-0.200pt]{0.400pt}{187.179pt}}
\put(130.0,82.0){\rule[-0.200pt]{315.338pt}{0.400pt}}
\put(1439.0,82.0){\rule[-0.200pt]{0.400pt}{187.179pt}}
\put(130.0,859.0){\rule[-0.200pt]{315.338pt}{0.400pt}}
\put(130.0,859.0){\color{red}{\rule[-0.200pt]{315.338pt}{0.400pt}}}

\put(146.0,131.0){\color{blue}{\rule[-0.200pt]{3.854pt}{1.2000pt}}} 

\put(178.0,179.0){\color{officegreen}{\rule[-0.200pt]{11.804pt}{1.2000pt}}} 

\put(243.0,228.0){\color{blue}{\rule[-0.200pt]{3.854pt}{1.2000pt}}} 

\put(275.0,276.0){\color{darkviolet}{\rule[-0.200pt]{35.171pt}{1.800pt}}} 

\put(437.0,325.0){\color{blue}{\rule[-0.200pt]{3.854pt}{1.2000pt}}} 

\put(469.0,373.0){\color{officegreen}{\rule[-0.200pt]{11.804pt}{1.200pt}}} 

\put(534.0,422.0){\color{blue}{\rule[-0.200pt]{3.854pt}{1.2000pt}}} 

\put(566.0,471.0){\color{harvardcrimson}{\rule[-0.200pt]{105.273pt}{1.200pt}}} 

\put(1019.0,519.0){\color{blue}{\rule[-0.200pt]{3.854pt}{1.200pt}}} 

\put(1051.0,567.0){\color{officegreen}{\rule[-0.200pt]{11.804pt}{1.200pt}}} 

\put(1116.0,616.0){\color{blue}{\rule[-0.200pt]{3.854pt}{1.2000pt}}} 

\put(1148.0,665.0){\color{darkviolet}{\rule[-0.200pt]{35.171pt}{1.200pt}}} 

\put(1310.0,713.0){\color{blue}{\rule[-0.200pt]{3.854pt}{1.2000pt}}} 

\put(1342.0,761.0){\color{officegreen}{\rule[-0.200pt]{11.804pt}{1.2000pt}}} 

\put(1407.0,810.0){\color{blue}{\rule[-0.200pt]{3.854pt}{1.2000pt}}} 

\put(130.0,82.0){\rule[-0.200pt]{0.400pt}{187.179pt}}
\put(130.0,82.0){\rule[-0.200pt]{315.338pt}{0.400pt}}
\put(130.0,859.0){\rule[-0.200pt]{315.338pt}{0.400pt}}
\end{picture}
\\ 

Fix $x \in [0,1]$.  
Then $n > m$ 
\\[-.25cm]

\qquad $\implies$ 
\[
\abs{\LebesgueFunction_n(x) - \LebesgueFunction_m(x)} 
\ < \ 
\frac{1}{2^m}
\]

\qquad $\implies$ 
\[
\LebesgueFunction (x) \ \equiv \ \lim\limits_n \ \LebesgueFunction_n 
\hspace{1cm} \text{exists $\forall \ x \in [0,1]$}.
\]

$\LebesgueFunction$ has the following properties: 
\\[-.25cm]

\qquad 1) \quad 
$\LebesgueFunction$ is continuous ($\LebesgueFunction_n \lra \LebesgueFunction$ uniformly).
\\[-.25cm]

\qquad 2) \quad 
$\LebesgueFunction$ is (not strictly) increasing.  
\\[-.25cm]

\qquad 3) \quad 
$
\begin{cases}
\ \LebesgueFunction(0) = 0
\\
\ \LebesgueFunction(1) = 1
\end{cases}
.
$
\\

\qquad 4) \quad 
$\LebesgueFunction^\prime = 0$ almost everywhere 
(being locally constant on  $[0,1] - \CantorSet$).
\\[-.5cm]

\qquad 5) \quad 
$\LebesgueFunction$ does not have property (N) ($\lambda(\LebesgueFunction (\CantorSet)) = 1$).
\\[-.25cm]
\end{x}

\begin{x}{\small\bf THEOREM} \ 
Let $f:[a,b] \ra \R$.  
Assume: $f^\prime(x)$ exists and is finite for all $x \in [a,b]$ and that $f^\prime$ is integrable there $-$then $f$ is 
absolutely continuous.
\\[-.5cm]

PROOF \ 
Condition (1) of the Main Theorem is satisfied 
(``differentiability'' $\implies$ ``continuity''), 
conditions (2) and (3) are given, and (4) is satisfied in view of the previous lemma.
\end{x}

The composition of two absolutely continuous functions need not be absolutely continuous.  
However: 

\begin{x}{\small\bf FACT} \ 
Suppose that $f:[a,b] \ra [c,d]$ and $g:[c,d] \ra \R$ are absolutely continuous $-$then 
$g \circ f \in \AC [a,b]$ iff $(g^\prime \circ f) \hsy f^\prime$ is integrable.
\\[-.25cm]

[Note: \ 
Interpret $g^\prime (f(x)) f^\prime(x)$ to be zero whenever $f^\prime (x) = 0$.] 
\end{x}

\chapter{
$\boldsymbol{\S}$\textbf{13}.\quad  MULTIPLICITIES}
\setlength\parindent{2em}
\setcounter{theoremn}{0}
\renewcommand{\thepage}{\S13-\arabic{page}}

\qquad 
Let $f:[a,b] \ra \R$ be a continuous function.  
Put 
\[
m \ = \ \min\limits_{[a,b]} \ f, 
\quad
M \ = \ \max\limits_{[a,b]} \ f.
\]
\\[-1.25cm]

\begin{x}{\small\bf NOTATION} \ 
Define a function $N(f; -) : \ ]-\infty, +\infty[ \hsy \ra \R$ by stipulating that $N(f; y)$ is the number of times that $f$
assumes the value $y$ in $[a,b]$, i.e., 
the number of solutions of the equation
\[
f(x) 
\ = \ 
y
\qquad 
(a \leq x \leq b).
\]

[Note: \ 
$N(f; y)$ is either 0, or a positive integer, or $+\infty$.]
\\[-.25cm]
\end{x}

\begin{x}{\small\bf DEFINITION} \ 
$N(f; -) $ is the \un{multiplicity function} attached to $f$.
\\[-.25cm]
\end{x}

\begin{x}{\small\bf THEOREM} \ 
$N(f; -)$ is a Borel measurable function and
\[
\int\limits_{-\infty}^{+\infty} \ N(f; -) 
\ = \ 
\tT_f [a,b].
\]
\\[-1cm]

PROOF \ 
Subdivide $[a,b]$ into $2^n$ equal parts, let
\[
I_{n \hsy i} 
\ = \ 
[a, a + (b - a) / 2^n], \quad i = 1, 
\]
and let
\[
I_{n \hsy i} 
\ = \ 
]a + (i - 1) (b - a) / 2^n , a + i (b - a) / 2^n],
i = 2, 3, \ldots, 2^n.
\]
Then $f$ maps each $I_{n \hsy i} $ to a segment (closed or not), viz. the segment from $m_i$ to $M_i$, where
\[
m_i \ = \ \inf\limits_{I_{n \hsy i}} \ f, 
\quad 
M_i \ = \ \sup\limits_{I_{n \hsy i}} \ f.
\]
The characteristic function $\chisubni$ of the set $f(I_{n \hsy i})$ is zero for $y > M_i$ \& $y < m_i$, 
one for $m_i < y < M_i$, while it may be zero or one at the two endpoints.  
Therefore
$\chisubni$ is Borel measurable, thus so is the function
\[
\chisubn (y)
\ = \ 
\sum\limits_{i = 1}^{2^n} \ \chisubni (y) 
\qquad 
(-\infty < y < +\infty).
\]
And
\begin{align*}
\int\limits_{-\infty}^{+\infty} \ \chisubn \ 
&=\ 
\sum\limits_{i = 1}^{2^n} \ \int\limits_{-\infty}^{+\infty} \ \chisubni
\\[15pt]
&=\ 
\sum\limits_{i = 1}^{2^n} \ (M_i - m_i)
\\[15pt]
&=\ 
\sum\limits_{i = 1}^{2^n} \ \osc (f; I_{n \hsy i}).
\end{align*}
Moreover
\[
\chisubn
\ \geq \ 
0, 
\quad 
\chisubn
\ \leq \ 
\chisubnplusone,
\]
which implies that
\[
\chi 
\ \equiv \ 
\lim\limits_{n \ra \infty} \ \chisubn
\]
is Borel measurable.  
Pass then to the limit:
\[
\int\limits_{-\infty}^{+\infty} \ \chi
\ = \ 
\lim\limits_{n \ra \infty} \ 
\int\limits_{-\infty}^{+\infty} \
\chisubn
\ = \ 
\tT_f [a,b],
\]
$f$ being continuous.  
Matters thereby reduce to establishing that
\[
\chi 
\ = \ 
N(f; -).
\]
First
\[
\forall \ n, 
\chisubn
\ \leq \ 
N(f; -) 
\implies 
\chi 
\ \leq \ 
N(f; -).
\]
Let now $q$ be a natural number not greater than $N(f; y)$, giving rise to $q$ distinct
roots
\[
x_1 
\ < \ 
x_2 
\ < \
\cdots 
\ < \
x_q
\]
of the equation
\[
f(x) 
\ = \ 
y 
\qquad 
(a \leq x \leq b).
\]
Upon choosing $n \gg 0$: 
\[
\frac{b - a}{2^n} 
\ < \ 
\min (x_{i + 1} - x_i), 
\]
it follows that all $q$ roots will fall into distinct intervals $I_{n \hsy i}$, hence
\[
\chisubn 
\ \geq \ 
q 
\implies 
\chi 
\ \geq \ 
q .
\]
If $N(f; y) = +\infty$, $q$ can be chosen arbitrarily large, thus $\chi (y) = +\infty$.  
On the other hand, if $N(f; y)$ is finite, take $q = N(f; y)$ to get 
\[
\chi (y) 
\ \geq \ 
N(f; y)
\implies 
\chi 
\ \geq \ 
N(f; -).
\]
\\[-1cm]
\end{x}

\begin{x}{\small\bf SCHOLIUM} \ 
A continuous function $f[a,b] \ra \R$ is of bounded variation iff its multiplicity function $N(f; -)$ is integrable.
\\[-.25cm]
\end{x}

\begin{x}{\small\bf \un{N.B.}}  \ 
If $f \in \BV [a,b] \cap C [a,b]$, then
\[
\{y : N(f; y) = +\infty\}
\]
is a set of Lebesgue measure 0.
\\[-.5cm]

[In fact, $N(f; -)$ is integrable, thus is finite almost everywhere.]
\\[-.25cm]
\end{x}

Maintain the assumption that $f[a,b] \ra \R$ is continuous.
\\

\begin{x}{\small\bf NOTATION} \ 
Given $J = [c,d] \subset [a,b]$, write
\[
\Phi(f; J, y) \ = \ 
\begin{cases}
\ +1 \ \ \text{if } \ f(c) < y < f(d)\\
\ -1 \ \ \text{if } \ f(c) > y > f(d)\\
\ \ 0 \quad \  \text{otherwise}
\end{cases}
, \quad \text{where $-\infty < y < +\infty$.}
\]
\\[-1cm]
\end{x}


\begin{x}{\small\bf LEMMA} \ 
If 
\[c 
\ = \ 
y_0 
\ < \ 
y_1 
\ < \ 
\cdots 
\ <  \ 
y_m 
\ = \ 
d
\]
is a partition of $J = [c,d]$ into the $m$ intervals $J_j = [y_{j - 1}, y_j]$ and 
$f(y_j) \neq y$ for $j = 0, 1, \ldots, m$, then
\[
\Phi(f; J, y) 
\ = \ 
\sum\limits_{j = 1}^m \ \Phi(f; J_j, y).
\]
\\[-1cm]
\end{x}

\begin{x}{\small\bf NOTATION} \ 
Given a finite system $S$ of nonoverlapping intervals $J = [c,d]$ in $[a,b]$, put 
\[
c \hsy N(f; y) 
\ = \ 
\sup\limits_S \ \sum\limits_{J \in S} \ \abs{\Phi(f; J, y) }.
\]
\\[-1cm]
\end{x}

\begin{x}{\small\bf DEFINITION} \ 
$c \hsy N(f; y)$ is the \un{corrected multiplicity function} attached to $f$.
\\[-.25cm]
\end{x}

Obviously
\[
0 
\ \leq \ 
c \hsy N(f; -) 
\ \leq \ 
+\infty.
\]
\\[-1cm]

\begin{x}{\small\bf THEOREM} \ 
$\forall \ y$, $-\infty, < y < +\infty$, 
\[
0
\ \leq \
c \hsy N(f; y) 
\ \leq \  
N(f; y) 
\]
and
\[
c \hsy N(f; y) 
\ = \  
N(f; y) 
\]
for all but countably many $y$.
\\[-.25cm]

Therefore
\[
\tT_f [a,b] 
\ = \ 
\int\limits_{-\infty}^{+\infty} \ N(f; -) 
\ = \ 
\int\limits_{-\infty}^{+\infty} \ c \hsy N(f; -).
\]
\end{x}

\chapter{
$\boldsymbol{\S}$\textbf{14}.\quad  LOWER SEMICONTINUITY}
\setlength\parindent{2em}
\setcounter{theoremn}{0}
\renewcommand{\thepage}{\S14-\arabic{page}}

\begin{x}{\small\bf EXAMPLE} \ 
(Fatou's Lemma) \ 
Suppose given a measure space $(\XX, \mu)$ and a sequence $\{f_n\}$ of nonnegative integrable functions such that 
$f_n \ra f$ almost everywhere $-$then

\[
\int\limits_\XX \ f \ \td\mu 
\ \leq \ 
\liminf \limits_{n \ra \infty} \ \int\limits_\XX \ f_n  \td\mu.
\]
\\[-.75cm]
\end{x}

\begin{x}{\small\bf THEOREM} \ 
Suppose that $f_n : [a,b] \ra \R$ $(n = 1, 2, \ldots)$  is a sequence of continuous functions that converges pointwise to 
$f : [a,b] \ra \R$ $-$then
\[
\tT_f [a,b] 
\ \leq \ 
\liminf \limits_{n \ra \infty} \ \tT_{f_n} [a,b].
\]
\\[-1cm]

PROOF \ 
Given $\varepsilon > 0$, there exists a partition $P = \{x_0, \ldots, x_m\}$ of $[a,b]$ such that 
\begin{align*}
\bigvee\limits_a^b \ (f;P) \ 
&= \ 
\sum\limits_{j = 1}^m \ \abs{f(x_j) - f(x_{j - 1})}
\\[15pt]
&> \ 
\tT_f [a,b]  - 2^{-1} \hsy \varepsilon
\end{align*}

\noindent
if $\tT_f [a,b]  < +\infty$ or $> \varepsilon^{-1}$ if $\tT_f [a,b]  = +\infty$.  
Since $f_n(x_j) \ra f(x_j)$ at each of the $m + 1$ points $x_0, \ldots, x_m$, there is an $n_\varepsilon$ such that 
\[
\abs{f(x_j) - f_n(x_j)}
\ < \ 
4^{-1} \hsy m^{-1} \hsy \varepsilon
\]

\noindent
for all $n \geq n_\varepsilon$ and $j = 0, \ldots, m$, hence if $n \geq n_\varepsilon$, 
\begin{align*}
\abs{f(x_j) - f(x_{j - 1})} \ 
&=\ 
\abs{f(x_j) - f_n(x_j) + f_n(x_j) - f_n(x_{j - 1}) - f(x_{j - 1}) + f_n(x_{j - 1})}
\\[15pt]
&\leq \ 
\abs{f(x_j) - f_n(x_j) } + \abs{ f(x_{j - 1}) + f_n(x_{j - 1})} + \abs{f_n(x_j) - f_n(x_{j - 1})}
\end{align*}

$\implies$

\[
\sum\limits_{j = 1}^m \ \abs{f(x_j) - f(x_{j - 1})}
\ \leq \ 
4^{-1} \hsy \varepsilon  \hsx + \hsx 4^{-1} \hsy \varepsilon  
\hsx + \hsx 
\sum\limits_{j = 1}^m \ \abs{f_n(x_j) - f_n(x_{j - 1})}
\]
or still, 
\begin{align*}
\sum\limits_{j = 1}^m \ \abs{f(x_j) - f(x_{j - 1})} - 2^{-1} \hsy \varepsilon \ 
&\leq \ 
\sum\limits_{j = 1}^m \ \abs{f_n(x_j) - f_n(x_{j - 1})}
\\[15pt]
&\leq \ 
\tT_{f_n} [a,b].
\end{align*}

\un{Case 1:} \ 
$\tT_f [a,b] \ < \  +\infty$ $-$then
\begin{align*}
\sum\limits_{j = 1}^m \ \abs{f(x_j) - f(x_{j - 1})} \ 
&> \ 
\tT_f [a,b] - 2^{-1} \hsy \varepsilon - 2^{-1} \hsy \varepsilon
\\[15pt]
&= \ 
\tT_f [a,b] - \varepsilon
\end{align*}

$\implies$

\[
\tT_f [a,b] - \varepsilon
\ < \ 
\tT_{f_n} [a,b] 
\qquad (n \geq n_\varepsilon)
\]

$\implies$

\[
\tT_f [a,b] - \varepsilon
\ \leq \ 
\liminf\limits_{n \ra \infty} \ \tT_{f_n} [a,b]
\]

$\implies (\varepsilon \downarrow 0)$

\[
\tT_f [a,b] 
\ \leq \ 
\liminf \limits_{n \ra \infty} \ \tT_{f_n} [a,b].
\]
\\[-.75cm]

\un{Case 2:} \ 
$\tT_f [a,b] = +\infty$ $-$then
\[
\sum\limits_{j = 1}^m \ \abs{f(x_j) - f(x_{j - 1})} \hsx - \hsx 2^{-1} \hsy \varepsilon
\ > \ 
\varepsilon^{-1} - 2^{-1} \varepsilon
\]

$\implies$

\[
\varepsilon^{-1} - 2^{-1} \varepsilon
\ < \ 
\tT_{f_n} [a,b] 
\qquad (n \geq n_\varepsilon)
\]

$\implies$

\[
+\infty
\ = \ 
\tT_f [a,b]
\ = \ 
\liminf \limits_{n \ra \infty} \ \tT_{f_n} [a,b].
\]
\\[-1cm]
\end{x}

\begin{x}{\small\bf REMARK} \ 
One cannot in general replace pointwise convergence by convergence almost everywhere, i.e., 
it can happen that under such circumstances
\[
\liminf \limits_{n \ra \infty} \ \tT_{f_n} [a,b]
\ < \ 
\tT_f [a,b].
\]
\\[-1cm]
\end{x}

\begin{x}{\small\bf EXAMPLE} \ 
Work on $[0, 2 \hsy \pi]$ and take
\[
f_n(x) 
\ = \ 
\frac{1}{n} \hsx \sin (n x), 
\]
so $f(x) = 0$ $-$then $f_n \ra f$ uniformly, 
\[
\tT_f [0, 2 \pi]  
\ = \ 
0, 
\quad 
\tT_{f_n} [0, 2 \pi]  
\ = \ 
4.
\]
\\[-1cm]
\end{x}

\begin{x}{\small\bf EXAMPLE} \ 
Work on $[0, 2 \hsy \pi]$ and take
\[
f_n(x) 
\ = \ 
\frac{1}{n} \hsx \sin (n^2 x), 
\]
so $f(x) = 0$ $-$then $f_n \ra f$ uniformly, 
\[
\tT_f [0, 2 \pi] 
\ = \ 
0, 
\quad 
\tT_{f_n} [0, 2 \pi] 
\ = \ 
+\infty.
\]
\\[-1cm]
\end{x}

\begin{x}{\small\bf THEOREM} \ 
Let $f : [a,b] \ra \R$  be a continuous function $-$then $c N (f;-)$ is lower semicontinuous in $]-\infty, +\infty[$, 
i.e., $\forall \ y_0$, 
\[
c N(f; y_0) 
\ \leq \ 
\liminf \limits_{y \ra y_0} \ c N (f; y).
\]
\\[-1cm]
\end{x}

\begin{x}{\small\bf THEOREM} \ 
Suppose that $f_n : [a,b] \ra \R$ is a sequence of continuous functions that converges pointwise to 
$f : [a,b] \ra \R$ $-$then $\forall \ y$, 
\[
c N(f; y) 
\ \leq \ 
\liminf \limits_{n \ra \infty} \ c N (f_n; y).
\]
\\[-1cm]
\end{x}

\begin{x}{\small\bf REMARK} \ 
These statements ensure that $c N$ is lower semicontinuous w.r.t $f$ and w.r.t. $y$ separately.  
More is true:  
$c N$ is lower semicontinuous w.r.t the pair $(f, y)$, i.e., if $f_n \ra f$, $y \ra y_0$, then
\[
c N(f; y_0) 
\ \leq \ 
\liminf \limits_{n \ra \infty} \ c N (f_n; y)
\]
as $f_n \ra f$, $y \ra y_0$.
\\[-.25cm]
\end{x}

\begin{x}{\small\bf \un{N.B.}} \ 
In the foregoing, one cannot in general replace $c N$ by $N$.
\end{x}

\chapter{
$\boldsymbol{\S}$\textbf{15}.\quad  FUNCTIONAL ANALYSIS}
\setlength\parindent{2em}
\setcounter{theoremn}{0}
\renewcommand{\thepage}{\S15-\arabic{page}}

\begin{x}{\small\bf THEOREM} \ 
$\BV[a,b]$ is a Banach space under the norm
\[
\norm{f}_\BV 
\ = \ 
\abs{f(a)} + \tT_f [a,b].
\]

[Note: \ 
$\tT_f[a,b]$ is not a norm since a constant function $f$ has zero total variation, 
hence the introduction of $\abs{f(a)}$.  
Recall, however, that
\[
\tT_{f + g} [a,b]
\ \leq \ 
\tT_f [a,b] + \tT_g [a,b]
\]
and
\[
\tT_{c \hsy f} [a,b]
\ = \ 
\abs{c} \hsx \tT_f [a,b] \hsx .]
\]
\\[-.75cm]

As a preliminary to the proof, consider a Cauchy sequence $\{f_k\}$ in $\BV[a,b]$.  
Given $\varepsilon > 0$, there exists $C_\varepsilon \in \ \N$ such that 
\[
\norm{f_k - f_\ell}_\BV 
\ = \ 
\abs{f_k (a) - f_\ell (a)} + \tT_{f_k - f_\ell} [a,b]
\ \leq \ 
\epsilon
\]
for all $k$, $\ell \geq C_\varepsilon$.  
Therefore
\[
\norm{f_k - f_\ell}_\infty 
\ \leq \ 
\varepsilon,
\]
thus the sequence $\{f_k\}$ converges uniformly to a bounded function $f:[a,b] \ra \R$, the claim being that $f \in \ \BV[a,b]$.  

This said, take a partition $P \in \ \sP[a,b]$ and note that
\begin{align*}
\sum\limits_{i = 1}^n \ \abs{(f_k - f_\ell) (x_i) - (f_k - f_\ell)(x_{i - 1})} \ 
&\leq \ 
\tT_{f_k - f_\ell} [a,b] 
\\[15pt]
&\leq \ 
\varepsilon
\end{align*}

\noindent
for all $k$, $\ell \geq C_\varepsilon$.  
From here, send $\ell$ to $+\infty$ to get
\\[-1cm]

\[
\sum\limits_{i = 1}^n \ \abs{(f_k - f) (x_i) - (f_k - f)(x_{i - 1})} 
\ \leq \ 
\varepsilon
\]
\\[-1.25cm]

\noindent 
for all $k \geq C_\varepsilon$, hence

\[
\tT_{f_k - f} [a,b] 
\ \leq \ 
\varepsilon
\]
for all $k \geq C_\varepsilon$.   
And
\[
\abs{f_k (a) - f_\ell (a)} 
\ \ra \ 
\abs{f_k (a) - f (a)}
\ \leq \ 
\varepsilon
\hps (\ell \ra +\infty).
\]
Therefore
\[
\norm{f_k - f}_\BV 
\ \leq \ 
2 \hsy \varepsilon
\]
for all $k \geq C_\varepsilon$.  
Moreover
\begin{align*}
\tT_f[a,b] \ 
&\leq 
\tT_{f - f_k}[a,b] + \tT_{f_k}[a,b]
\\[11pt]
&< \ 
+\infty.
\end{align*}
So $f \in \ \BV[a,b]$ and $f_k \ra f$ in $\BV[a,b]$.
\\[-.25cm]
\end{x}

\begin{x}{\small\bf REMARK} \ 
$\BV[a,b]$, equipped with the norm $\norm{\hsx \cdot \hsx}_\BV$, is not separable.  
\\[-.25cm]

[Take $[a,b] = [0,1]$ and for $f \in \ \BV[0,1]$, $r > 0$, let
\[
\tS(f,r) 
\ = \ 
\{g \in \ \BV[0,1] \hsx : \hsx \norm{g - f}_\BV < r\}.
\]
Call $\chisubt$ $(0 < t < 1)$ the characteristic function of $\{t\}$ $-$then for $t_1 \neq t_2$, 
\begin{align*}
\normx{\chisubtsubone - \chisubtsubtwo}_\BV
&=\ 
(\chisubtsubone - \chisubtsubtwo) (a) + \tT_{\chisubtsubone - \chisubtsubtwo} \hsx [0,1]
\\[11pt]
&=\ 
0 + \tT_{\chisubtsubone - \chisubtsubtwo} \hsx [0,1]
\\[11pt]
&=\ 
4.
\end{align*}
But this implies that 
\[
\tS(\chisubtsubone, 1) \hsx \cap \hsx \tS(\chisubtsubtwo, 1)
\ \neq \ 
\emptyset.
\]
In fact
\[
\begin{cases}
\ \normx{h - \chisubtsubone}_\BV \ < 1
\\[8pt]
\ \normx{h - \chisubtsubtwo}_\BV \ < 1
\end{cases}
\]

$\implies$
\begin{align*}
\normx{\chisubtsubone - \chisubtsubtwo}_\BV
&=\ 
\normx{\chisubtsubone - h + h - \chisubtsubtwo}_\BV
\\[8pt]
&\leq \ 
\normx{\chisubtsubone - h}_\BV + \normx{\chisubtsubtwo - h}_\BV
\\[8pt]
&< \ 
1 + 1
\\[8pt]
&= \ 
2.
\end{align*}
Accordingly there exists a continuum of disjoint spheres $\tS(\chisubt, 1) \subset \tS(0,3)$, 
hence an arbitrary sphere $\tS(f,r)$ contains a continuum of disjoint spheres 
$\tS(r \hsy \chisubt/3 + f, r/3)$.]
\\[-.25cm]
\end{x}

\begin{x}{\small\bf THEOREM} \ 
$\BV[a,b]$ is a complete metric space under the distance function
\[
\td_\BV (f,g) 
\ = \ 
\int\limits_a^b \ \abs{f - g} + \abs{\tT_f[a,b] - \tT_g[a,b]}.
\]
\\[-1cm]

The issue is completeness and for this, 
it suffices to establish that the balls $B_M$ of radius $M$ centered at 0 are compact, 
the claim being that every sequence $\{f_n\} \subset B_M$ has a subsequence converging to a limit in $B_M$.
\end{x}

\begin{x}{\small\bf \un{N.B.}} \ 
Spelled out, $B_M$ is the set of functions $f \in \ \BV[a,b]$ satisfying the condition
\[
\td_\BV (f,0) 
\ = \ 
\int\limits_a^b \ \abs{f} + \tT_f [a,b]
\ \leq \ 
M.
\]
\\[-1cm]
\end{x}

\begin{x}{\small\bf HELLY'S SELECTION THEOREM} \ 
Let $\sF$ be an infinite family of functions in $\BV[a,b]$.  
Assume that there exists a point $x_0 \in \ [a,b]$ and a constant $K > 0$
such that $\forall \ f \in \ \sF$,
\[
\abs{f(x_0)} + \tT_f [a,b] 
\ \leq \ 
K.
\]
Then there exists a sequence $\{f_n\} \subset \sF$ and a function $g \in \ \BV[a,b]$ such that
\[
f_n \ra g \hps (n \ra \infty)
\]
pointwise in $[a,b]$.
\\[-.25cm]
\end{x}

\begin{x}{\small\bf LEMMA} \ 
$\forall \ f \in \ B_M$, 
\[
\abs{f(a)}
\ \leq \ 
M \hsx \left(1 + \frac{1}{b - a}\right).
\]

PROOF \ 
Write
\[
f(a) 
\ = \ 
f(a) - f(x) + f(x)
\]

$\implies$
\\[-2cm]

\allowdisplaybreaks
\begin{align*}
\abs{f(a)} \ 
&\leq \ 
\abs{f(a) - f(x) } + \abs{f(x)}
\\[11pt]
&\leq \ 
\tT_f [a,b] + \abs{f(x)}
\end{align*}

$\implies$
\\[-1.75cm]

\allowdisplaybreaks
\begin{align*}
\abs{f(a)} \ \int\limits_a^b \ 1 \ 
&\leq \ 
 \int\limits_a^b \  \tT_f [a,b] \ + \   \int\limits_a^b \ \abs{f}
 \\[11pt]
&\leq \ 
M \hsx (b - a) + M
\end{align*}

$\implies$
\\[-1.55cm]

\[
\abs{f(a)}
\ \leq \ 
M \hsx \left(1 + \frac{1}{b - a}\right).
\]
\\[-.75cm]

In the HST, take $\sF = \{f_n\}$, $x_0 = a$, and
\[
K
\ = \ 
M \hsx \left(1 + \frac{1}{b - a}\right)+ M.
\]
Then there exists a subsequence $\{f_{n_k}\}$ and a function $g \in \ \BV[a,b]$ such that
\[
f_{n_k} \ra g \hps (k \ra \infty)
\]
pointwise in $[a,b]$.
\end{x}


\begin{x}{\small\bf LEMMA} \ 
$\forall$ $n_k$, $\forall$ $x \in \ [a,b]$, 
\[
\abs{f_{n_k} (x)} 
\ \leq \ 
\abs{f_{n_k} (a)} + T_{f_{n_k}} [a,b] 
\ < \ 
+\infty.
\]

The $f_{n_k}$ are therefore bounded, hence by dominated convergence, 
\[
f_{n_k} \ra g \hps (k \ra \infty)
\]
in $\Lp^1[a,b]$.
\\[-.25cm]

Consider now the numbers
\[
T_{f_{n_k}} [a,b] \hps (k = 1, 2, \ldots).
\]
They constitute a bounded set, 
hence there exists a subsequence $\{T_{f_{n_k}} [a,b]\}$ (not relabeled) which converges to a limit  $\tau$.  
Since $f_{n_k}$ tends to $g$ pointwise, on the basis of lower semicontinuity, it follows that
\[
\tT_g [a,b] 
\ \leq \
\lim\limits_{k \ra \infty} \ T_{f_{n_k}} [a,b],
\] 
which implies that
\[
\tT_g [a,b] 
\ \leq \
\tau.
\]
Adjusting $g$ at $a$ if necessary, matters can be arranged so as to ensure that $\tT_g[a,b] = \tau$.
\\[-.25cm]

Consequently
\begin{align*}
\td_\BV(f_{n_k}, g) \ 
&=\ 
\int\limits_0^1 \ 
\abs{f_{n_k} - g} \ + \ \abs{T_{f_{n_k}} [a,b] - \tT_g [a,b]}.
\\
&
\hspace{1.5cm}
\downarrow \ (k \ra \infty) 
\hspace{1cm}
\downarrow \ (k \ra \infty) 
\\
&
\hspace{1.6cm}
0 
\hspace{2.50cm}
\abs{\tau -\tau}
\end{align*}
I.e.:
\[
\lim\limits_{k \ra \infty} \ \td_\BV(f_{n_k}, g) 
\ = \ 
0.
\]
\\[-1cm]

The final detail is the verification that $g \in B_M$.  
To this end, fix $\varepsilon > 0$ $-$then for $k \gg 0$, 
\begin{align*}
\td_\BV (g, 0) \ 
&\leq \ 
\td_\BV(g, f_{n_k}) + \td_\BV(f_{n_k}, 0) 
\\[11pt]
&\leq \ 
\varepsilon + M.
\end{align*}
\\[-1cm]
\end{x}

\begin{x}{\small\bf LEMMA} \ 
In the $\td_\BV$ metric, $\BV[a,b]$ is separable.
\\[-.25cm]
\end{x}

\begin{x}{\small\bf LEMMA} \ 
$\forall \ a \in \ \R$, $\forall \ f, \hsx g \in \ \BV[a,b]$, 
\[
\td_\BV (a \hsy f, a \hsy g) 
\ = \ 
\abs{a} \hsy \td_\BV (f, g).
\]
\\[-1cm]
\end{x}

\begin{x}{\small\bf THEOREM} \ 
Let $\alpha \in \ \Lp^1[a,b]$ $-$then the assignment
\[
f \ \ra \ \int\limits_a^b \ f \hsx \alpha \ \equiv \ \Lambda_\alpha (f)
\]
is a continuous linear functional on $\BV[a,b]$ when equipped with the $\td_\BV$ metric.
\\[-.25cm]

PROOF \ 
To establish the continuity, take an $f \in \ \BV[a,b]$ and suppose that $\{f_n\}$ is a sequence in $\BV[a,b]$ such that
\[
\td_\BV (f_n, f) \ \ra \ 0 \hps (n \ra \infty),
\]
the objective being to show that if $\varepsilon > 0$ be given, then
\[
\abs{\Lambda_\alpha (f_n) - \Lambda_\alpha (f)} \ < \ \varepsilon
\]
provided $n \gg 0$.

So fix a constant $C > 0$: $\forall$ $n$, 
\[
\int\limits_0^1 \ \abs{f_n - f} + \abs{\tT_{f_n} [a,b] - \tT_f [a,b]} 
\ \leq \ 
C.
\]
For each $n$ choose a point $\bar{x}_n$ such that 
\[
\abs{f_n (\bar{x}_n) - f(\bar{x}_n)} \ \leq \ C
\]
and note that for all $x \in \ [a,b]$, 
\[
\begin{cases}
\ \abs{f_n (x) - f_n(\bar{x}_n)} \ \leq \ \tT_{f_n} [a,b]
\\[8pt]
\ \abs{f_n (x) - f(\bar{x}_n)}  \ \leq \ \tT_{f} [a,b]
\end{cases}
\]
and 
\[
\tT_{f_n} [a,b]
\ \leq \ 
\tT_{f} [a,b] + C
\]

$\implies$
\allowdisplaybreaks
\begin{align*}
\abs{f_n (x)  - f (x)} \ 
&\leq \ 
\abs{f_n (x) - f_n(\bar{x}_n) +  f_n(\bar{x}_n) -  f(\bar{x}_n) +  f(\bar{x}_n) - f(x)}
\\[11pt]
&\leq \ 
\abs{f_n (x) - f_n(\bar{x}_n} + \abs{f(x) - f(\bar{x}_n)} + \abs{f_n(\bar{x}_n) - f(\bar{x}_n)}
\\[11pt]
&\leq \ 
\tT_{f_n} [a,b] + \tT_{f} [a,b] + \abs{f_n(\bar{x}_n) - f(\bar{x}_n)}
\\[11pt]
&\leq \ 
\tT_f[a,b] + C + \tT_f[a,b] + C
\\[11pt]
&\leq \ 
2 \hsx \tT_f[a,b] + 2 \hsx C
\\[11pt]
&\equiv \ 
K.
\end{align*}
On general grounds (absolute continuity of the integral), given $\varepsilon > 0$ there exists $\delta > 0$ such that 
\[
\int\limits_E \ K \hsx \abs{\alpha} 
\ < \ 
\frac{\varepsilon}{2}
\]
if $\lambda(E) < \delta$.  
Take now $N \gg 0$: 
\[
\lambda(E_N) < \delta \hsx 
\hspace{.7cm} (E_N = \{x: \abs{\alpha(x)} > N\}).
\]
Then 
\begin{align*}
\abs{\Lambda_\alpha (f_n) - \Lambda_\alpha (f)} \ 
&=\ 
\bigg| \hsx \int\limits_a^b \ f_n \hsx \alpha \ - \ \int\limits_a^b \ f \hsx \alpha \hsx \bigg| 
\\[15pt]
&\leq \ 
\int\limits_a^b \ \abs{f_n - f} \hsx \abs{\alpha}
\\[15pt]
&=\ 
\int\limits_{E_N} \ \abs{f_n - f} \hsx \abs{\alpha} 
\hsx + \hsx 
\int\limits_{E_N^c} \ \abs{f_n - f} \hsx \abs{\alpha}
\\[15pt]
&\leq \
\int\limits_{E_N} \ K \hsx \abs{\alpha} 
\hsx + \hsx 
\int\limits_{E_N^c} \ \abs{f_n - f} \hsx \abs{\alpha}
\\[15pt]
&< \
\frac{\varepsilon}{2} \hsx + \hsx 
\int\limits_{E_N^c} \ \abs{f_n - f} \hsx \abs{\alpha}.
\end{align*}
And
\[
x \in \ E_N^c \implies \abs{\alpha(x)} \leq N
\]

$\implies$
\allowdisplaybreaks
\begin{align*}
\int\limits_{E_N^c} \ \abs{f_n - f} \hsx \abs{\alpha} \ 
&\leq \ 
N \ \int\limits_{E_N^c} \ \abs{f_n - f} 
\\[15pt]
&\leq \ 
N \ \int\limits_a^b \ \abs{f_n - f}
\\[15pt]
&< \ 
\frac{\varepsilon}{2}  \hps (n \gg 0).
\end{align*}
Therefore in the end

\[
\abs{\Lambda_\alpha (f_n) - \Lambda_\alpha (f)} 
\ < \ 
\frac{\varepsilon}{2} + \frac{\varepsilon}{2} 
\ = \ 
\varepsilon
\]
for all $n$ sufficiently large.
\\[-.25cm]
\end{x}

\begin{x}{\small\bf \un{N.B.}} \ 
\[
\Lambda_{\alpha_1} 
\ = \ 
\Lambda_{\alpha_2}
\]
iff $\alpha_1 = \alpha_2$ almost everywhere.
\\[-.25cm]

[Suppose that $\Lambda_{\alpha_1} = \Lambda_{\alpha_2}$.
Define $f_t \in \ \BV[a,b]$ by the prescription
\[
f_t(x) \quad
\begin{cases}
\ 1 \hps (0 \leq x \leq t)\\
\ 0 \hps (t < x \leq 1)
\end{cases}
.
\]
Then
\[
\int\limits_a^b \ f_t \hsx \alpha_1 
\ = \ 
\int\limits_a^b \ f_t \hsx \alpha_2
\]

$\implies$
\[
\int\limits_0^t \ \alpha_1 
\ = \ 
\int\limits_0^t \ \alpha_2
\]

$\implies$
\[
\alpha_1 
\ = \ 
\alpha_2
\]
almost everywhere.]
\end{x}

\chapter{
$\boldsymbol{\S}$\textbf{16}.\quad  DUALITY}
\setlength\parindent{2em}
\setcounter{theoremn}{0}
\renewcommand{\thepage}{\S16-\arabic{page}}

\qquad
In the abstract theory, take $\XX = [a,b]$ $-$then there is an isometric isomorphism
\[
\Lambda \hsy : \hsy \sM([a,b]) \ra C [a,b]^*,
\]
viz. the rule that sends a finite signed measure $\mu$ to the bounded linear functional
\[
f \ra \int\limits_{[a, b]} \ f \ \td \mu.
\]
On the other hand, it is a point of some importance that there is another description of $C[a,b]^*$ 
which does not involve any measure theory at all.  
\\[-.25cm]

\begin{x}{\small\bf RAPPEL} \ 
If $f$ is continuous on $[a,b]$ and if $g \in \BV[a,b]$, then the Stieltjes integral
\[
\int\limits_a^b \ f(x) \ \td g(x)
\]
exists.
\\[-.25cm]
\end{x}

\begin{x}{\small\bf NOTATION} \ 
$C[a,b]$ is the set of continuous functions on $[a,b]$ equipped with the supremum norm
\[
\norm{f}_\infty
\ = \ 
\sup\limits_{[a,b]} \ \abs{f},
\]
and $C[a,b]^*$  is its dual.
\\[-.25cm]
\end{x}

\begin{x}{\small\bf LEMMA} \ 
Let $g \in \BV[a,b]$ $-$then the assignment
\[
f  \ra 
\int\limits_a^b \ f(x) \ \td g(x)
\]
defines a bounded linear functional $\Lambda_g \in C[a,b]^*$.
\\[-.5cm]

[Note: \ 
\[
\forall \ f, \quad
\abs{\Lambda_g (f)} 
\ \leq \ 
\tT_g[a,b] \hsy \norm{f}_\infty,
\]
hence
\[
\norm{\Lambda_g}
\ \leq \ 
\tT_g [a,b].]
\]
\\[-1cm]
\end{x}

\begin{x}{\small\bf RIESZ REPRESENTATION THEOREM} \ 
If $\Lambda$ is a bounded linear functional on $C [a,b]$, then there exists a $g \in \BV [a,b]$ such that
\[
\Lambda (f) 
\ = \ 
\int\limits_a^b \ f(x) \ \td g(x) 
\quad (= \Lambda_g (f))
\]
for all $f \in C [a,b]$.  
And: 
\[
\norm{\Lambda}
\ = \ 
\tT_g [a,b].
\]
\\[-1cm]

PROOF \ 
Extend $\Lambda$ to $\Lp^\infty [a,b] \supset C [a,b]$ without increasing its norm (Hahn-Banach).  
\\[-.5cm]

\noindent
Given $x \in [a,b]$, let
\[
u_x(t) \ = \quad
\begin{cases}
\ 1 \hspace{0.7cm} (a \leq t \leq x)\\
\ 0 \hspace{0.7cm} (x < t \leq b)
\end{cases}
\]
and put
\[
g(x) 
\ = \ 
\Lambda (u_x).
\]

Claim: \ $g \in \BV[a,b]$ and in fact
\[
\tT_g [a,b]
\ \leq \ 
\norm{\Lambda}.
\]
Thus take a partition $P \in \sP[a,b]$ and let
\[
\varepsilon_i 
\ = \ 
\sgn (g(x_i) - g(x_{i - 1})) 
\quad (i = 1, \ldots, n).
\]

Then
\begin{align*}
\sum\limits_{i = 1}^n \ \abs{g(x_i) - g(x_{i - 1})} \ 
&= \ 
\sum\limits_{i = 1}^n \ \varepsilon_i  \ (g(x_i) - g(x_{i - 1}))
\\[15pt]
&= \ 
\sum\limits_{i = 1}^n \ \varepsilon_i  \ (\Lambda(u_{x_i}) - \Lambda(u_{x_{i - 1}}))
\\[15pt]
&= \ 
\Lambda \hsx \
\big(
\sum\limits_{i = 1}^n \ \varepsilon_i  \ (u_{x_i} - u_{x_{i - 1}})
\big)
\\[15pt]
&\leq 
\norm{\Lambda} \ 
\normx
{
\sum\limits_{i = 1}^n \ \varepsilon_i  \ (u_{x_i} - u_{x_{i - 1}})
}
\\[15pt]
&\leq 
\norm{\Lambda}.
\end{align*}
Therefore
\[
\tT_g [a,b] 
\ \leq \ 
\norm{\Lambda} 
\ < \ 
+\infty 
\implies 
g \in \BV[a,b].
\]
Suppose next that $f \in C[a,b]$ and let
\[
x_i 
\ = \ 
a + \frac{i \hsy (b - a)}{n} \hspace{0.7cm} (i = 0, \ldots, n).
\]
Define
\[
f_n(x) 
\ = \ 
\sum\limits_{i = 1}^n \ 
f(x_i) \hsx (u_{x_i} (x) - u_{x_{i - 1}} (x)).
\]
Then
\begin{align*}
\norm{f - f_n}_\infty \ 
&=\ 
\sup\limits_{[a,b]} \ \abs{f - f_n}
\\[11pt]
&\leq \ 
\max\limits_{1 \leq i \leq n} \ 
\sup \ \{\abs{f(x) - f(x_i)} \hsx : \hsx x_{i -1} \leq x \leq x_i\}.
\end{align*}
Invoking uniform continuity, it follows that
\[
\norm{f - f_n}_\infty \ra 0 \quad (n \ra +\infty),
\]
i.e., 
\begin{align*}
f_n \ra f \implies \Lambda (f) \ 
&= \ 
\lim\limits_{n \ra \infty} \ \Lambda (f_n)
\\[15pt]
&= \ 
\lim\limits_{n \ra \infty} \ \sum\limits_{i = 1}^n \ f(x_i) \hsx \big(\Lambda (u_{x_i}) - \Lambda (u_{x_{i - 1}})\big)
\\[15pt]
&= \ 
\lim\limits_{n \ra \infty} \ \sum\limits_{i = 1}^n \ f(x_i) \hsx (g(x_i) - g(x_{i - 1}))
\\[15pt]
&= \ 
\int\limits_a^b \ f(x) \ \td g(x) 
\\[15pt]
&= \ 
\Lambda_g(f).
\end{align*}
From the above, 
\[
\tT_g [a,b] 
\ \leq \ 
\norm{\Lambda}
\]
and
\[
\norm{\Lambda}
\ \leq \ 
\tT_g [a,b].
\]
So
\[
\norm{\Lambda}
\ = \ 
\tT_g [a,b], 
\]
as contended.
\\[-.25cm]
\end{x}

The ``$g$'' that figures in this theorem is definitely not unique.  
To remedy this, proceed as follows.
\\[-.25cm]

\begin{x}{\small\bf DEFINITION} \ 
$g \in \BV[a,b]$ is \un{normalized} if $g(a) = 0$ and $g(x+) = g(x)$ when $a < x < b$.
\\[-.5cm]

[Note: \ 
Since $g(a) = 0$, 
\[
\norm{g}_\BV 
\ = \ 
\tT_g [a,b].
\]
Observe too that by definition, the right continuous modification $g_r$ of $g$ in $]a,b[$ is given by the formula
\[
g_r(x) 
\ = \ 
g(x+), 
\]
so the assumption is that $g_r = g$, i.e., in $]a,b[$, $g$ is right continuous.]
\\[-.25cm]
\end{x}

\begin{x}{\small\bf NOTATION} \ 
Write $\NBV [a,b]$ for the linear subspace of $\BV [a,b]$ whose elements are normalized.
\\[-.25cm]
\end{x}

\begin{x}{\small\bf THEOREM} \ 
The arrow
\[
\NBV [a,b] \ \ra \ C [a,b]^*
\]
that sends $g$ to $\Lambda_g$ is an isometric isomorphism:
\[
\norm{g}_\BV 
\ = \ 
\tT_g [a,b] 
\ = \ 
\norm{\Lambda_g}.
\]

Here is a sketch of the proof.  
\\[-.25cm]

\qquad \un{Step 1:} \ 
Define an equivalence relation in $\BV[a,b]$ by writing $g_1 \sim g_2$ iff $\Lambda_{g_1} = \Lambda_{g_2}$.
\\

\qquad \un{Step 2:} \ 
Note that 
\begin{align*}
g \sim 0 \
&\implies
0 \ = \ \int\limits_a^b \ \td g(x) \ = \ g(b) - g(a)
\\
&\implies
g(a) \ = \  g(b).
\end{align*}

\qquad \un{Step 3:} \ 
Establish that 
\[
g \sim 0 
\]

\qquad $\implies$
\[
g(a) \ = \ g(c+) \ = \ g(c-) \ = \ g(b)
\]
if $a < c < b$.
\\[-.5cm]

[Suppose that 
\[
a \leq c < b, 
\quad 
0 < h < b - c
\]
and define 
\[
f(x) \ = \quad 
\begin{cases}
\ 1 \hspace{2.55cm} (a \leq x \leq c)
\\[8pt]
\ 1 - \ds\frac{x - c}{h} \hspace{1cm} (c \leq x \leq c + h)
\\[8pt]
\ 0 \hspace{2.55cm} (c + h \leq x \leq b)
\end{cases}
.
\]
Then
\[
g \sim 0 
\]

\qquad $\implies$
\[
0 
\ = \ 
\int\limits_a^b \ f(x) \ \td g(x) 
\ = \ 
g(c) - g(a) + \int\limits_c^{c + h} \ f(x) \ \td g(x).
\]
Integrate
\[
 \int\limits_c^{c + h} \ f(x) \ \td  g(x)
\]
by parts to get
\[
- g(c) \hsx + \hsx \frac{1}{h} \ \int\limits_c^{c + h} \ g(x) \ \td x
\]

\qquad $\implies$ $(h \ra 0)$
\[
0 
\hsx = \hsx
g(c) - g(a) - g(c) - g(c+)
\]

\qquad $\implies$
\[
g(a) 
\ = \
g(c+).
\]
\\[-1cm]

Analogously
\[
a \hsx < \hsx c \hsx \leq \hsx b 
\implies
g(b) \hsx = \hsx g(c-).]
\]
\\[-1cm]

\qquad \un{Step 4:} \ 
Establish that if $g \in \BV[a,b]$ and if 
\[
g(a) 
\ = \ 
g(c+)
\ = \ 
g(c-) 
\ = \ 
g(b)
\]
when $a < c < b$, then $g \sim 0$.
\\[-.5cm]

[In fact, $g(x) = g(a)$ at $x = a$, $x = b$, and at all interior points of $[a,b]$ at which $g$ is continuous, 
thus $\forall$ $f \in C [a,b]$, 
\[
\int\limits_a^b \ f(x) \ \td  g(x)
\ = \ 
\int\limits_a^b \ f(x) \ \td  h(x)
\ = \ 
0, 
\]
where $h(x) \equiv g(a)$.]
\\[-.25cm]

\qquad \un{Step 5:} \ 
Every equivalence class contains at most one normalized function.
\\[-.25cm]

[If $g_1$, $g_2 \in \NBV[a,b]$ and if $g_1 \sim g_2$, then 
$g \equiv g_1 - g_2 \sim 0.$  
By hypothesis, $g_1(a) = 0$, $g_2(a) = 0$, so 
\\[-1.25cm]

\begin{align*}
(g_1 - g_2) (a) = 0 
&\implies 
(g_1 - g_2) (b) = 0 
\\
&\implies 
g_1(b) - g_2(b) = 0
\\
&\implies
g_1(b) = g_2(b).
\end{align*}

Moreover
\begin{align*}
g(c+) = g(a) = 0
&\implies 
g_1(c+) - g_2(c+) = 0
\\
&\implies 
g_1(c+) = g_2(c+).
\end{align*}
On the other hand, 
\[
\begin{cases}
\ g_1 \in \NBV[a,b] \implies g_1(c+) = g_1(c) 
\\[11pt]
\ g_2 \in \NBV[a,b] \implies g_2(c+) = g_2(c) 
\end{cases}
\implies g_1(c) = g_2(c).
\]
I.e.: \ $g_1 = g_2$.]
\\[-.25cm]

\qquad \un{Step 6:} \ 
Every equivalence class contains at least one normalized function.  
\\[-.5cm]

[Given $g \in \NBV[a,b]$, define $g^* \in \NBV[a,b]$ as follows: 
\\[-.75cm]

\[
\begin{cases}
\ g^*(a) = 0
\\[8pt]
\ g^*(b) = g(b) - g(a) 
\end{cases}
\]
\[
g^*(x) = g(x+) - g(a) \quad (a < x < b).
\]
Then $g^* \in \NBV[a,b]$ and $g^* \sim g$.  
The verification that $g^* \in \NBV[a,b]$ is immediate.  
\\[-.5cm]

\noindent
There remains the claim that $g^* - g \sim 0.$
\\[-.25cm]

\qquad \textbullet \quad $(g^* - g) (a) = g^*(a) - g(a) = -g(a)$.
\\[-.5cm]

\qquad \textbullet \quad $(g^* - g) (b) = g^*(b) - g(b) = g(b) - g(a) - g(b) = -g(a)$.
\\[-.25cm]

\noindent When $a < x < b$, 
\[
g^*(x) 
\ = \ 
g_r(x) - g(a).
\]
And for $c \in \ ]a,b[$, 
\[
\begin{cases}
\ \lim\limits_{x \downarrow c} \ g_r(x) \  = \ \lim\limits_{x \downarrow c} \ g(x)\\[11pt]
\ \lim\limits_{x \uparrow c} \ g_r(x) \  = \ \lim\limits_{x \uparrow c} \ g(x)
\end{cases}
.
\]
\allowdisplaybreaks
\begin{align*} 
\text{\textbullet \quad  $(g^* - g) (c+)$} \ 
&=\ 
g^*(c+) - g(c+) \hspace{6cm}
\\[8pt]
&= \ 
g_r(c+) - g(a) - g(c+)
\\[8pt]
&= \ 
\lim\limits_{x \downarrow c} \ g_r(x) - g(a) - g(c+)
\\[8pt]
&= \ 
\lim\limits_{x \downarrow c} \ g(x) - g(a) - g(c+)
\\[8pt]
&= \
g(c+) - g(a) - g(c+)
\\[8pt]
&= \
-g(a).
\end{align*}
\\[-2.25cm]

\allowdisplaybreaks
\begin{align*}
\text{\textbullet \quad  $(g^* - g) (c-)$} \ 
&=\ 
g^*(c-) - g(c-) \hspace{6cm}
\\[8pt]
&= \ 
g_r(c-) - g(a) - g(c-)
\\[8pt]
&= \ 
\lim\limits_{x \uparrow c} \ g_r(x) - g(a) - g(c-)
\\[8pt]
&= \ 
\lim\limits_{x \uparrow c} \ g(x) - g(a) - g(c-)
\\[8pt]
&= \
g(c-) - g(a) - g(c-)
\\[8pt]
&= \
-g(a).
\end{align*}
Therefore
\[
g^* - g  \sim 0 
\implies 
g^* \sim g.]
\]
\\[-1cm]

\qquad \un{Step 7:} \ 
\[
\tT_{g^*} [a,b]
\ \leq \ 
\tT_g [a,b].
\]
\\[-1cm]

[Let $P \in \sP [a,b]$: 
\[
a \ = \ x_0 \ < \ x_1  \ < \ \cdots  \ < \ x_n \ = \ b.
\]
Given $\varepsilon > 0$, choose points $y_1, \ldots, y_{n-1}$ at which $g$ is continuous with $y_i$ so close to 
$x_i$ (on the right) that
\[
\abs{g(x_i+) - g(y_i)} 
\ < \ 
\frac{\varepsilon}{2 \hsy n}.
\]
Taking $y_0 = a$, $y_n = b$, there follows
\begin{align*}
\sum\limits_{i = 1}^n \ \abs{g^*(x_i) - g^*(x_{i - 1})} \ 
&= \ 
\sum\limits_{i = 1}^n \ \abs{g(x_i+) - g(a) - g(x_{i - 1}+) + g(a) } 
\\[15pt]
&\leq \ 
\sum\limits_{i = 1}^n \ \abs{g(x_i+) - g(y_i)} 
+ 
\sum\limits_{i = 1}^n \ \abs{g(x_{i - 1}+) - g(y_{i - 1})} 
\\[15pt]
&\hspace{2.6cm}
+ 
\sum\limits_{i = 1}^n \ \abs{g(y_i) - g(y_{i - 1})} 
\\[15pt]
&\leq \ 
\sum\limits_{i = 1}^n \ \abs{g(y_i) - g(y_{i - 1})} + \varepsilon
\end{align*}

$\implies$
\[
\tT_{g^*} [a,b]
\ \leq \ 
\tT_g [a,b] + \varepsilon
\]

$\implies$ $(\varepsilon \ra 0)$
\[
\tT_{g^*} [a,b]
\ \leq \ 
\tT_g [a,b].
\]
\\[-1cm]

Consider now the arrow
\[
\NBV [a,b] \ra C [a,b]^*
\]
that sends $g$ to $\Lambda_g$.  
To see that it is surjective, let $\Lambda \in C [a,b]^*$ and choose a $g \in \BV [a,b]$ such that
\[
\Lambda_g 
\ = \ 
\Lambda.
\]
The equivalence class to which $g$ belongs contains a unique normalized element $g^*$, 
so $g^* \sim g$
\\[-.5cm] 

$\implies$
\[
\Lambda_{g^*}
\ = \ 
\Lambda_g 
\ = \ 
\Lambda.
\]
Finally, as regards the norms,
\begin{align*}
\norm{\Lambda} \ 
&= \ 
\norm{\Lambda_g}
\\[8pt]
&= \ 
\norm{\Lambda_{g^*}}
\\[8pt]
&\leq \ 
\tT_{g^*} [a,b]
\\[8pt]
&\leq \ 
\tT_g [a,b]
\\[8pt]
&= \ 
\norm{\Lambda}.
\end{align*}
Meanwhile

\[
\tT_{g^*} [a,b]
\ = \ 
\norm{g^*}_\BV 
\implies
\norm{\Lambda}
\ = \ 
\norm{g^*}_\BV .
\]
\end{x}

\chapter{
$\boldsymbol{\S}$\textbf{17}.\quad  INTEGRAL MEANS}
\setlength\parindent{2em}
\setcounter{theoremn}{0}
\renewcommand{\thepage}{\S17-\arabic{page}}

\qquad 
To simplify the notation, work in $[0,1]$ (the generalization to $[a,b]$ being straightforward).
\\

\begin{x}{\small\bf NOTATION} \ 
$I = [0,1]$, $0 < \delta < 1$, $I_\delta = [0, 1 - \delta]$ 
($\implies  1 - \delta > 0$), 
$0 < h < \delta$ $(\implies 1 - h > 1 - \delta)$. 
\\[-.25cm]
\end{x}

\begin{x}{\small\bf DEFINITION} \ 
Let $f \in \BV [0,1]$ and suppose that $f$ is continuous $-$then its 
\un{integral mean} is the function $f^h$ on $[0, 1 - \delta]$ defined by the prescription
\[
f^h (x)
\ = \ 
\frac{1}{h} \ \int\limits_0^h \ f(x + t) \ \td t 
\qquad (0 \leq x \leq 1 - \delta).
\]
\\[-1cm]
\end{x}

\begin{x}{\small\bf LEMMA} \ 
$f^h \in C [I_\delta]$ and 
\[
f^h \ra f \quad (h \ra 0)
\]
uniformly in $I_\delta$.
\\[-.25cm]
\end{x}

\begin{x}{\small\bf LEMMA} \ 
The derivative of $f^h$ exists on $]0, 1 - \delta[$ and is given there by the formula

\[
(f^h)^\prime (x) 
\ = \ 
\frac{f(x + h) - f(x)}{h}.
\]

[Note: \ 
Therefore $f^h$ has a continuous first derivative in the interior of $I_\delta$.]
\\[-.25cm]
\end{x}

\begin{x}{\small\bf LEMMA} \ 
\[
f^h \in \AC [0, 1 - \delta].
\]

PROOF \ 
Let
\[
M 
\ = \ 
\sup\limits_{[0,1]} \ \abs{f}.
\]
Then for fixed $h$, 
\begin{align*}
\abs{(f^h)^\prime (x)} \ 
&=\ 
\abs{\frac{f(x + h) - f(x)}{h}} 
\hspace{1cm} (0 < x < 1 - \delta)
\\[11pt]
&\leq \ 
\frac{2 \hsy M}{h}.
\end{align*}
Choose $a < b$ such that 
\[
0 
\ < \ 
a 
\ < \ 
b 
\ < \ 
1 - \delta.
\]
Then
\[
f^h (b) - f^h(a) 
\ = \ 
\int\limits_a^b \ (f^h)^\prime (x) \td x
\]

$\implies$
\[
\abs{f^h (b) - f^h(a) } 
\ \leq \ 
\frac{2 \hsy M \hsy (b - a)}{h}
\qquad (0 < a < b - 1 - \delta)
\]
or still, by continuity, 
\[
\abs{f^h (b) - f^h(a)} 
\ \leq \ 
\frac{2 \hsy M \hsy (b - a)}{h}
\qquad (0 \leq a < b \leq 1 - \delta).
\]
And this implies that $f^h$ is absolutely continuous.
\\[-.25cm]

[In the usual notation, 
\[
\sum\limits_{k = 1}^n \ \abs{f^h (b_k) - f^h(a_k)} 
\ \leq \ 
\frac{2 \hsy M}{h} \ 
\sum\limits_{k = 1}^n \ (b_k - a_k).]
\]
\\[-.75cm]
\end{x}

\begin{x}{\small\bf LEMMA} \ 
Let
\[
[a,b] \subset I_\delta.
\]
Then
\[
\tT_{f^h} [a,b] 
\ \leq \  
\tT_f [a,b + \delta] 
\quad (0 < h < \delta).
\]
\\[-.75cm]

PROOF \ 
Take a finite system of intervals $[a_i, b_i]$ $(1 \leq i \leq n)$ without common interior points in $[a,b]$ $-$then
\[
\sum\limits_{i = 1}^n \ \abs{f(b_i + t) - f(a_i + t)} 
\ \leq \ 
\tT_f [a,b + \delta] 
\]

$\implies$
\begin{align*}
\sum\limits_{i = 1}^n \ \abs{f^h (b_i) - f^h(a_i)} \ 
&\leq \ 
\frac{1}{h} \ \int\limits_0^h \ \tT_f [a,b + \delta] \  \td t
\\[15pt]
&=\ 
\tT_f [a,b + \delta] 
\\[-1cm]
\end{align*}

$\implies$
\[
\tT_{f^h} [a,b]
\leq \ 
\tT_f [a,b + \delta] 
\quad (0 < h < \delta).
\]
\\[-1cm]
\end{x}

\begin{x}{\small\bf THEOREM} \ 
Let
\[
[a,b] \subset I_\delta.
\]
Then
\[
\tT_{f^h} [a,b] \ra \tT_f [a,b] 
\quad (0 < h \ra 0).
\]

PROOF \ 
\[
\tT_{f^h} [a,b] 
\ \leq \ 
\tT_f [a,b + \delta] 
\quad (0 < h < \delta)
\]

$\implies$
\[
\limsup\limits_{h \ra 0} \ \tT_{f^h} [a,b] 
\ \leq \ 
\tT_f [a,b  + \delta].
\]
Since
\[
\tT_f [a,b + \delta] 
\ra 
\tT_f [a,b] 
\quad (\delta \ra 0), 
\]
it follows that
\[
\limsup\limits_{h \ra 0} \ \tT_{f^h} [a,b] 
\ \leq \ 
\tT_f [a,b] .
\]
\\[-.75cm]

\noindent
By hypothesis, $[a,b] \subset I_\delta$ and in $I_\delta$, 
\[
f^h \ra f \quad (h \ra 0)
\]
uniformly, hence pointwise.  
Therefore
\[
\liminf\limits_{h \ra 0} \ \tT_{f^h} [a,b]
\ \geq \ 
\tT_f [a,b].
\]
\\[-1cm]
\end{x}

\begin{x}{\small\bf SCHOLIUM} \ 
Owing to the absolute continuity of $f^h$ in $I_\delta$, for any $[a,b] \subset I_\delta$, we have
\begin{align*}
\tT_{f^h} [a,b] \ 
&=\ 
\int\limits_a^b \ \abs{(f^h)^\prime (x)} \ \td x
\\[15pt]
&=\ 
\int\limits_a^b \ \abs{\frac{f(x + h) - f(x)}{h}} \ \td x 
\end{align*}
and 
\[
\int\limits_a^b \ \abs{\frac{f(x + h) - f(x)}{h}} \ \td x 
\ra 
\tT_f [a,b] 
\quad (0 < h \ra 0).
\]
\end{x}

\chapter{
$\boldsymbol{\S}$\textbf{18}.\quad  ESSENTIAL VARIATION}
\setlength\parindent{2em}
\setcounter{theoremn}{0}
\renewcommand{\thepage}{\S18-\arabic{page}}

\begin{x}{\small\bf DEFINITION} \ 
$\BV\Lp^1 ]a,[b$ is the subset of $\Lp^1 ]a,b[$ consisting of those $f$ whose distributional derivative 
$\tD f$ is represented by a finite signed Radon measure in $]a,b[$ of finite total variation, i.e., if
\[
\int\limits_{]a,b[} \ f \hsy \phi^\prime 
\ = \ 
-
\int\limits_{]a,b[} \ \phi   \ \td \tD f 
\qquad (\forall \ \phi \in C_c^\infty ]a,b[)
\]
for some finite signed Radon measure $\tD f$ with 
\[
\abs{\tD f} \hsx ]a,b[ 
\ < \ 
+\infty.
\]
\\[-1.25cm]

[Note: \ 
Two $\Lp^1$ -functions which are equal almost everywhere define the same distribution 
(and so have the same distributional derivative).]
\\[-.25cm]
\end{x}

\begin{x}{\small\bf \un{N.B.}} \ 
A smoothing argument shows that the integration by parts formula is still true for all $\phi \in  C_c^1 \hsx ]a,b[$.
\\[-.25cm]
\end{x}

Of course it may happen that $\tD f$ is a function, say $\tD f = g \hsy \td \hsy x$, hence $\forall \ \phi \in C_c^1 \hsx ]a,b[$,
\[
\int\limits_{]a,b[} \ f \hsy \phi^\prime 
\ = \ 
-
\int\limits_{]a,b[} \ \phi \hsx g \hsy \td x. 
\]
\\[-.75cm]

\begin{x}{\small\bf EXAMPLE} \ 
Work in $]0,2[$ and let
\[
f(x) \ = \ 
\begin{cases}
\ x \hspace{0.5cm} (0 < x \leq 1)
\\[4pt]
\ 1 \hspace{0.5cm} (1 < x < 2)
\end{cases}
.
\]
Put
\[
g(x) \ = \ 
\begin{cases}
\ 1 \hspace{0.5cm} (0 < x \leq 1)
\\[4pt]
\ 0 \hspace{0.5cm} (1 < x < 2)
\end{cases}
.
\]
\\[-1cm]

Then $\tD f = g \td x$.  
In fact, $\forall \ \phi \in C_c^1 ]0,2[$, 
\begin{align*}
\int\limits_0^2 \ 
f \hsy \phi^\prime 
\ \td x 
&=\ 
\int\limits_0^1 \ 
x \phi^\prime 
\ \td x 
\hsx + \hsx 
\int\limits_1^2 \ 
\phi^\prime 
\ \td x 
\\[15pt]
&=\ 
-
\int\limits_0^1 \ 
\phi
\ \td x 
\hsx + \hsx 
\phi (1) - \phi (1)
\\[15pt]
&=\ 
-
\int\limits_0^1 \
\phi
\ \td x 
\\[15pt]
&=\ 
-
\int\limits_0^2 \ 
\phi  \hsy g
\ \td x.
\end{align*}
\\[-.75cm]
\end{x}

\begin{x}{\small\bf EXAMPLE} \ 
Let $\mu$ be a finite signed Radon measure in $]a,b[$.  
Put $f(x) = \mu (]a,x[)$ $-$then the distributional derivative of $f$ is $\mu$.  
\\[-.25cm]

[$\forall \ \phi \in  C_c^1 \hsx ]a,b[$,
\begin{align*}
\int\limits_{]a,b[} \ 
f(x) \hsy \phi^\prime (x) 
\ \td x \ 
&=\ 
\int\limits_{]a,b[} \ 
\int\limits_{]a,x[} \ 
\phi^\prime (x) 
\ \td \mu (y) \hsx \td x
\\[15pt]
&=\ 
\int\limits_{]a,b[} \ 
\int\limits_{]y,b[} \ 
\phi^\prime (x) 
\ \td x  \hsx \td \mu (y)
\\[15pt]
&=\ 
-
\int\limits_{]a,b[} \
\phi (y) 
\ \td \mu (y).]
\end{align*}
\\[-.75cm]
\end{x}

\begin{x}{\small\bf NOTATION} \ 
Let $f : ]a,b[ \ra \R$ $-$then the \un{total variation} $\tT_f \hsx ]a,b[$ of $f$ in $]a,b[$ is the supremum of the total variations of $f$ 
in the closed subintervals of $]a,b[$.
\\[-.25cm]
\end{x}

\begin{x}{\small\bf FACT} \ 
If $f : [a,b] \ra \R$ $-$then
\[
\tT_f [a,b] 
\ = \ 
\tT_f \hsx ]a,b[ \ + \hsx \abs{f(a+) - f(a)} 
\hsx + \hsx 
\abs{f(b-) - f(b)}.
\]
\\[-1cm]
\end{x}

\begin{x}{\small\bf \un{N.B.}} \ 
Therefore
\[
\tT_f [a,b]
\ = \ 
\tT_f \hsx ]a,b[
\]
whenever $f$ is continuous.  
\\[-.25cm]
\end{x}


\begin{x}{\small\bf DEFINITION} \ 
A function $f : ]a, b[ \hsx \ra \R$ is of \un{bounded variation} in $]a, b[$ provided
\[
\tT_f \hsx ]a,b[ \ < \ +\infty.
\]
\\[-1cm]
\end{x}

\begin{x}{\small\bf NOTATION} \ 
$\BV ]a,b[$ is the set of functions of bounded variation in $]a,b[$.
\\[-.25cm]
\end{x}

\begin{x}{\small\bf \un{N.B.}} \ 
Elements of $\BV ]a,b[$ are bounded, hence are integrable:
\[
\BV ]a,b[ \subsetx \Lp^1 ]a,b[ \hsx .
\]
Moreover, $\forall \ f \in \BV ]a,b[$, 
\[
\begin{cases}
& f(a+) \\
& f(b-)
\end{cases}
\quad \text{exist}.
\]
\\[-.75cm]
\end{x}

\begin{x}{\small\bf EXAMPLE} \ 
Take $]a,b[ \ = \  ]0, 1[$ $-$then
\[
f(x) 
\ = \ 
\frac{1}{1 - x}
\]
is increasing and of bounded variation in every closed subinterval of $]0,1[$, yet $f \notin \BV ]0,1[\hsy.$
\\[-.25cm]
\end{x}

The initial step in the theoretical development is to characterize the elements of $\BVL^1 ]a,b[\hsy.$
\\

\begin{x}{\small\bf FACT} \ 
Let $\mu$ be a finite signed Radon measure in $]a,b[$ $-$then for any open set $S \subset ]a,b[$, 
\[
\abs{\mu} (S)
\ = \ 
\sup\limits_{\substack{\phi \hsy \in \hsy C_c(S), \\[0.15cm]\norm{\phi}_\infty \hsy\leq \hsy 1}} \ 
\bigg\{
\ 
\int\limits_{]a,b[} \ \phi \ \td \mu 
\bigg\}.
\]
\\[-.75cm]
\end{x}

\begin{x}{\small\bf DEFINITION} \ 
Given $f \in \Lp^1 ]a,b[$, let
\[
\sV (f; ]a,b[) 
\ = \ 
\sup\limits_{\substack{\phi \hsy \in \hsy C_c^1]a,b[, \\[0.15cm] \norm{\phi}_\infty \hsy\leq \hsy 1}} \ 
\bigg\{
\ 
\int\limits_{]a,b[} \ f \hsy \phi^\prime
\bigg\}.
\]
\\[-.75cm]
\end{x}

\begin{x}{\small\bf THEOREM} \ 
Let $f \in \Lp^1 ]a,b[$ $-$then $f \in \BVL^1 ]a,b[$ iff 
\[
\sV (f; ]a,b[) 
\ < \ 
+\infty.
\]
And when this is so, 
\[
\sV (f; ]a,b[) 
\ = \ 
\abs{\tD f} \hsx ]a,b[.
\]
\\[-1cm]

PROOF \ 
Suppose first that $f \in \BVL^1 ]a,b[$ $-$then
\begin{align*}
\sV (f; ]a,b[)  \ 
&=\ 
\sup\limits_{\substack{\phi \hsy \in \hsy C_c^1]a,b[, \\[0.15cm] \norm{\phi}_\infty \hsy\leq \hsy 1}} \ 
\bigg\{
-
\int\limits_{]a,b[} \ 
\phi 
\ \td \tD f 
\bigg\}
\\[15pt]
&=\ 
\sup\limits_{\substack{\phi \hsy \in \hsy C_c]a,b[, \\[0.15cm] \norm{\phi}_\infty \hsy\leq \hsy 1}} \ 
\bigg\{
-
\int\limits_{]a,b[} \ 
\phi 
\ \td \tD f 
\bigg\}
\\[15pt]
&=\ 
\abs{-\tD f} \hsx ]a,b[
\\[15pt]
&=\ 
\abs{\tD f} \hsx ]a,b[
\\[15pt]
&< \ 
+\infty.
\end{align*}
Conversely assume that
\[
\sV (f; ]a,b[) 
\ < \ 
+\infty.
\]

\noindent
Then
\[
\bigg| \hsx \int\limits_{]a,b[} \ f \hsx \phi^\prime \hsx \bigg|
\ \leq \ 
\sV (f; ]a,b[) \hsx \norm{\phi}_\infty.
\]
Since $C_c^1 ]a,b[$ is dense in $C_0 ]a,b[$, the linear functional
\[
\Lambda : C_c^1 ]a,b[ \ \ra \R
\]
defined by the rule
\[
\phi \ra \int\limits_{]a,b[} \ f \hsx \phi^\prime
\]
can be extended uniquely to a continuous linear functional
\[
\Lambda : C_0 ]a,b[ \ \ra \R,
\]

\noindent
where
\[
\norm{\Lambda}^* 
\ \leq \ 
\sV (f; ]a,b[).
\]
Thanks to the ``$C_0$'' version of the RRT, there exists a finite signed Radon measure $\mu$ in $]a,b[$ such that
\[
\norm{\Lambda}^* 
\ = \ 
\abs{\mu} \hsy (]a,b[)
\]
and
\[
\Lambda (\phi) 
\ = \ 
\int\limits_{]a,b[} \ \phi \ \td \mu 
\qquad (\forall \ \phi \in C_0 ]a,b[).
\]
Definition: 
\[
\tD f 
\ = \ 
\mu
\]

$\implies$
\\[-1.7cm]

\allowdisplaybreaks
\begin{align*}
\abs{\tD f} \hsx ]a,b[ \ 
&=\ 
\abs{\mu} \hsy (]a,b[)
\\[11pt]
&=\ 
\norm{\Lambda}^* 
\\[11pt]
&\leq \ 
\sV (f; ]a,b[) 
\\[11pt]
&< \ 
+\infty.
\end{align*}
\\[-1cm]
\end{x}

\begin{x}{\small\bf LEMMA} \ 
The map
\[
f \ra \sV (f; ]a,b[) 
\]
is lower semicontinuous in the $\Lp_\loc^1 ]a,b[$ topology.
\\[-.25cm]
\end{x}

\begin{x}{\small\bf APPLICATION} \ 
The map
\[
f \ra \abs{\tD f} \hsx ]a,b[
\]
is lower semicontinuous in the $\Lp_\loc^1 ]a,b[$ topology.
\\[-.25cm]
\end{x}

\begin{x}{\small\bf SUBLEMMA} \ 
Any element of $\BV \hsx ]a,b[$ can be represented as the difference of two bounded increasing functions.
\\[-.25cm]
\end{x}

\begin{x}{\small\bf LEMMA} \ 
$\forall \ f \in \BV \hsx ]a,b[$, 
\[
\sV (f; ]a,b[)  
\ \leq \ 
\tT_f \hsx ]a,b[ 
\qquad ( < +\infty).
\]
\\[-1cm]

PROOF \ 
Construct a sequence $\chisubn$ of step functions such that 
\[
\chisubn \ra f 
\quad (n \ra \infty)
\]
in $\Lp_\loc^1 ]a,b[$ and $\forall \ n$, 
\[
\sV (\chisubn; ]a,b[)
\ \leq \ 
\tT_f \hsx ]a,b[ \hsx .
\]
Thanks now to lower semicontinuity, 
\begin{align*}
\sV (f; ]a,b[)  \ 
&\leq \ 
\liminf\limits_{n \ra \infty} \ \sV (\chisubn; ]a,b[)  \
\\[11pt]
&\leq \ 
\tT_f \hsx ]a,b[  \hsx .
\end{align*}

\end{x}

\begin{x}{\small\bf SCHOLIUM} \ 
\[
\BV ]a,b[ 
\subsetx 
\BVL^1 ]a,b[.
\]
\\[-1cm]

[Note :  \ 
If $f : [a,b] \ra \R$ is in $\BV [a,b]$, then its restriction to $]a,b[$ is in $\BV ]a,b[$, hence is in $\BVL^1 ]a,b[$.]
\\[-.25cm]
\end{x}

\begin{x}{\small\bf DEFINITION} \ 
Let $f \in \Lp^1 \hsx ]a,b[$ $-$then the \un{essential variation} of $f$, denoted $e - \tT_f \hsx ]a,b[$, is the set
\[
\inf \ \{T_g \hsx ]a,b[ \ : \hsx g = f \ \text{almost everywhere}\}.
\]

[Note: \ 
If $f$,  $f_2 \in \Lp^1 ]a,b[$ and if $f_1 = f_2$ almost everywhere, then 
\[
e - T_{f_1} \hsx ]a,b[ 
\hsx \ = \ 
e - T_{f_2} \hsx ]a,b[ \hsx .]
\]
\\[-1cm]
\end{x}

\begin{x}{\small\bf LEMMA} \ 
Let $f \in \Lp^1 \hsx ]a,b[$ $-$then
\[
e - \tT_f \hsx ]a,b[
\hsx \ = \ 
\sV (f; ]a,b[).
\]
\\[-1.5cm]
\end{x}

Consequently
\\[-.25cm]

\begin{x}{\small\bf THEOREM} \ 
Let $f \in \Lp^1 ]a,b[$  $-$then
\[
e - \tT_f \hsx ]a,b[
\ < \ 
+\infty 
\iff 
f \in \BVL^1 ]a,b[.
\]
And then
\[
\abs{\tD f} \hsy ]a,b[ 
\hsx \ = \ 
e - \tT_f \hsx ]a,b[.
\]
\\[-1cm]
\end{x}

\begin{x}{\small\bf LEMMA} \ 
Let $f \in \BVL^1 ]a,b[$.  
Assume: \ 
$\tD f = 0$ $-$then $f$ is (equivalent to) a unique constant.  
\\[-.25cm]
\end{x}

Assuming still that $f \in \BVL^1 ]a,b[$, let $\mu = \tD f$ and put $w (x) = \mu (]a,x[)$ $-$then $\tD w = \mu$, 
thus $\tD (f - w) = 0$, so there exists a unique constant $C$ such that 
\[
f 
\ = \ C + w
\]
almost everywhere.
\\[-.25cm]

\begin{x}{\small\bf LEMMA} \ 
\[
\tT_{C + w} ]a,b[ 
\ = \ 
e - \tT_f \hsx ]a,b[.
\]
\\[-1cm]

PROOF \ 
Take points
\[
x_0 < x_1 < \cdots < x_n
\]
in $]a,b[$ $-$then

\[
\sum\limits_{i = 1}^n \ 
\abs{(C + w) (x_i) - (C + w) (x_{i - 1})}
\ \leq \ 
\abs{\mu} (]a,b[)
\]

$\implies$ 
\begin{align*}
\tT_{C + w} ]a,b[  
&\leq \ 
\sV (f; ]a,b[)
\\[11pt]
&= \ 
e - \tT_f \hsx ]a,b[.
\end{align*}
\\[-1cm]
\end{x}

\begin{x}{\small\bf DEFINITION} \ 
Given $f \in \BVL^1 ]a,b[$, a function $g \in \Lp^1 ]a,b[$ such that $g = f$ almost everywhere is \un{admissible} if
\[
\tT_g ]a,b[ 
\ = \ 
e - \tT_f \hsx ]a,b[.
\]
\\[-1cm]

[Note: \ Since
\[
e - \tT_f \hsx ]a,b[ 
\ < \ 
+\infty 
\implies 
\tT_g ]a,b[
\ < \ 
+\infty,
\]
this says that $f$ is equivalent to $g$, where $g \in \BV \hsx ]a,b[$.]
\\[-.25cm]
\end{x}

So, in this terminology, $C + w$ is admissible, i.e., 
\[
f^\ell (x) 
\ \equiv \ 
C + \tD f ]a,x[
\]
is admissible, the same being the case of 
\[
f^r (x) 
\ \equiv \ 
C + \tD f ]a,x[.
\]
\\[-.75cm]

\begin{x}{\small\bf LEMMA} \ 
\[
\begin{cases}
\ f^\ell \ \text{is left continuous}
\\[4pt]
\ f^r \ \text{is right continuous}
\end{cases}
.
\]
\\[-1cm]
\end{x}

\begin{x}{\small\bf REMARK} \ 
\[
\begin{cases}
\ f^\ell (x) - f^\ell (y) \ = \ \tD f[y,x[
\\[4pt]
\ f^r (x) - f^r (y) \ = \ \tD f]y,x]
\end{cases}
\quad (a < y < x < b).
\]
\\[-1cm]
\end{x}

\begin{x}{\small\bf THEOREM} \ 
A function $g \in \Lp^1 \hsx ]a, b[$ is admissible iff
\[
g \in \{ \theta \hsy f^\ell + (1 - \theta) f^r \hsy : \hsy 0 \leq \theta \leq 1\}.
\]
\\[-1cm]
\end{x}

\begin{x}{\small\bf \un{N.B.}} \ 
Denote by $\AT_f$ the atoms of the theory, i.e., the $x \in \hsx ]a,b[$ such that $\tD f(\{x\}) \neq 0$ $-$then 
$f^\ell = f^r$ in $]a,b[ \ - \ \AT_f$ and every admissible $g$ is continuous in $]a,b[ \ - \ \AT_f$.
\\[-.25cm]
\end{x}

\begin{x}{\small\bf LEMMA} \ 
Suppose that $g \in \Lp^1 \hsx ]a, b[$ is admissible $-$then $g$ is differentiable almost everywhere and its derivative 
$g^\prime$ is the density of $\tD f$ w.r.t. Lebesgue measure.
\\[-.25cm]


There is a characterization of the essential variation which is purely internal.
\\[-.25cm]
\end{x}

\begin{x}{\small\bf NOTATION} \ 
Given an $f \in \Lp^1 \hsx ]a, b[$, let $C_\ap (f)$ stand for its set of points of approximate continuity.
\\[-.25cm]

[Recall that $C_\ap (f)$ is a subset of $]a,b[$ of full measure.]
\\[-.25cm]
\end{x}

\begin{x}{\small\bf LEMMA} \ 
\[
e - \tT_f \hsx ]a,b[ 
\ \  = \ 
\sup\ \sum\limits_{i = 1}^n \ \abs{f(x_i) - f(x_{i - 1})},
\]
where the supremum is taken over all finite collections of points $x_i \in C_\ap (f)$ subject to
\[
a 
\hsx < \hsx 
x_0
\hsx < \hsx 
x_1
\hsx < \hsx 
\cdots
\hsx < \hsx 
x_n
\hsx < \hsx 
b.
\]
\end{x}

\chapter{
$\boldsymbol{\S}$\textbf{19}.\quad  BVC}
\setlength\parindent{2em}
\setcounter{theoremn}{0}
\renewcommand{\thepage}{\S19-\arabic{page}}

\begin{x}{\small\bf NOTATION} \ 
Given a subset $M \subset \ ]a,b[$ of Lebesgue measure 0, denote by 
$\sP_M ]a,b[$ the set of all sequences
\[
P \hsx : \hsx 
x_0 \hsx < \hsx x_1 \hsx < \hsx  \cdots \hsx < \hsx  x_n,
\]
where
\[
\begin{cases}
\ a < x_0 \\
\ x_n < b
\end{cases}
\]
and
\[
x_i \in \ ]a,b[ \hsx - \hsx M 
\quad (i = 0, 1, \ldots, n).
\]

[Note: \ 
The possibility that $M = \emptyset$ is not excluded.]
\\[-.25cm]
\end{x}

\begin{x}{\small\bf NOTATION} \ 
Given a function $f : \hsx ]a,b[ \ \ra \ \R$, let $f_M$ be the restriction of $f$ to $]a,b[ \ - \ M$.
\\[-.25cm]
\end{x}

\begin{x}{\small\bf NOTATION} \ 
Given an element $P \in \sP_M \hsx ]a,b[$, put
\[
\bigvee\limits_a^b \ (f_M; P)
\ = \ 
\sum\limits_{i = 1}^n \ \abs{f_M(x_i) - f_M(x_{i - 1})}.
\]
\\[-1cm]
\end{x}

\begin{x}{\small\bf NOTATION} \ 
Given a function $f : \hsx ]a,b[ \ \ra \ \R$, put
\[
\tT_{f_M} \hsx ]a,b[
\ \ = \ 
\sup\limits_{P \in \sP_M \hsx ]a,b[} \ \bigvee\limits_a^b \ (f_M; P).
\]
\\[-1cm]
\end{x}

\begin{x}{\small\bf DEFINITION} \ 
$\tT_{f_M} \hsx ]a,b[$ is the \un{total variation} of $f_M$ in $]a,b[ \hsx - \hsx M$.
\\[-.25cm]
\end{x}

\begin{x}{\small\bf DEFINITION} \ 
A function $f \in \Lp^1 ]a,b[$\  is said to be of 
\un{bounded variation} \un{in the sense of Cesari} 
if there exists a subset $M \subset \ ]a.b[$ of Lebesgue measure 0 such that 
\[
\tT_{f_M} \hsx ]a,b[
\ < \ 
+\infty.
\]
\\[-1cm]
\end{x}

\begin{x}{\small\bf NOTATION} \ 
$\BVC \hsx ]a,b[$ is the set of functions which are of bounded variation in the sense of Cesari.
\\[-.25cm]
\end{x}

\begin{x}{\small\bf EXAMPLE} \ 
\[
\BV \hsx ]a,b[  
\subsetx
\BVC \hsx ]a,b[ 
\qquad (M = \emptyset).
\]
\\[-1cm]
\end{x}

\begin{x}{\small\bf THEOREM} \ 
\[
\BVC \hsx ]a,b[ 
\ = \ 
\BVL^1 \hsx ]a,b[ \hsx .
\]
\\[-1.5cm]
\end{x}

Proceed via a couple of lemmas. 
\\[-.25cm]

\begin{x}{\small\bf LEMMA} \ 
Suppose that $f \in \BVL^1 ]a,b[$ \ $-$then $f \in \BVC ]a,b[\hsy.$
\\[-.25cm]

PROOF \ 
The assumption that
\[
f \in \BVL^1 ]a,b[ 
\implies 
e - \tT_f ]a,b[ \ < +\infty.
\]
So there exists a $g$: $g = f$ almost everywhere and
\[
\tT_g ]a,b[
\ < \ 
+\infty.
\]
Take now for $M$ the set of $x$ such that $g(x) \neq f(x)$, the complement $]a,b[ \ - \hsx M$ being 
the set of $x$ where $g(x) = f(x)$.  
Consider a typical sum
\[
\sum\limits_{i = 1}^n \ \abs{f_M(x_i) - f_M(x_{i - 1})}
\]
which is equal to 
\[
\sum\limits_{i = 1}^n \ \abs{g(x_i) - g(x_{i - 1})}
\]
which is less than or equal to 
\[
\tT_g ]a,b[
\ < \ 
+\infty.
\]
Therefore $f \in \BVC]a,b[$.
\\[-.25cm]
\end{x}

\begin{x}{\small\bf SUBLEMMA} \ 
If $\tT_{f_M} ]a,b[  \ < \ +\infty$, then there exists a $g : ]a,b[ \ \ra \ \R$ such that $g_M = f_M$ and 
\[
\tT_g ]a,b[
\ = \ 
\tT_{f_M} ]a,b[ \hsx .
\]
\\[-1cm]
\end{x}

\begin{x}{\small\bf LEMMA} \ 
Suppose that $f \in \BVC ]a,b[$ $-$then $f \in \BVL^1 ]a,b[$\hsx .
\\[-.25cm]

PROOF \ 
The assumption that $f \in \BVC ]a,b[$ produces an ''$M$'' and from the preceding consideration, 
\[
g_M = f_M
\ \implies \ 
\restr{g}{]a,b[} - M 
\ = \ 
\restr{f}{]a,b[} - M,
\]
hence $g = f$ almost everywhere.  
But
\begin{align*}
\tT_{f_M} ]a,b[ \ < \ +\infty\ 
&\implies \ 
\tT_g ]a,b[ \ < \ +\infty\
\\[15pt]
&\implies \ 
g \in \BV ]a,b[
\\[15pt]
&\implies \ 
g \in \BVL^1 ]a,b[ \hsy .
\end{align*}
Since $g = f$ almost everywhere, they have the same distributional derivative, thus $f \in \BVL^1 ]a,b[$.
\\[-.25cm]
\end{x}

Let $\sM$ be the set of all subsets of $]a,b[$ of Lebesgue measure 0.
\\[-.25cm]

\begin{x}{\small\bf NOTATION} \ 
Given an $f \in \BVL^1 ]a,b[$, put
\[
\phi(f) 
\ = \ 
\inf\limits_{M \in \sM} \ \tT_{f_M} \hsx ]a,b[ \hsx .
\]
\\[-1cm]
\end{x}

\begin{x}{\small\bf THEOREM} \ 
\[
e - \tT_f ]a,b[ 
\ = \ 
\phi(f).
\]

PROOF \ 
To begin with, 
\[
f \in \BVL^1 ]a,b[ 
\implies 
e - \tT_f ]a,b[ 
\ < \ 
+\infty.
\]
On the other hand, $f \in \BVC ]a,b[$, so there exists $M \in \sM$: 
\[
\tT_{f_M} ]a,b[  
\ < \ 
+\infty
\ \implies \ 
\phi(f) 
\ < \ 
+\infty.
\]

\qquad \textbullet \quad 
$e - \tT_f ]a,b[ \ \ \leq \ \phi(f)$.
\\[-.25cm]

[Denote by $\sM_f$ the subset of $\sM$ consisting of those $M$ such that 
$\tT_{f_M} ]a,b[ \ < \ +\infty$.  
Assign to each $M \in \sM_f$ a function 
$g : \ ]a,b[ \ \ra \R$ such that $g_M = f_M$ and 
\[
\tT_g ]a,b[
\ \ = \ 
\tT_{f_M} ]a,b[.
\]
Therefore
\[
\{ \tT_{f_M} ]a,b[ \ : \hsx M \in \sM_f\} 
\ \subset \ 
\{\tT_g ]a,b[ \ : g = f \ \text{almost eveywhere}\}
\]

$\implies$
\begin{align*}
\phi(f) \ 
&=\ 
\inf\limits_{M \in \sM_f} \  \tT_{f_M} \ ]a,b[
\\[15pt]
&\geq \ 
e - \tT_f ]a,b[ \hsx.]
\end{align*}

\qquad \textbullet \quad
$\phi(f) \ \leq \  e - \tT_f ]a,b[$.
\\[-.25cm]

[Denote by $\sM_E$ the subset of $\sM$ consisting of those $M$ that arise from the elements 
$\tT_g ]a,b[$ in the set defining $e - \tT_f ]a,b[$ (i.e., per the requirement that $g = f$ almost everywhere) $-$then
\[
\tT_{f_M} \ ]a,b[
\ \  \leq \ 
\tT_g ]a,b[ 
\qquad (M \in \sM_E),
\]
hence
\begin{align*}
\phi(f) \ &=\ 
\inf\limits_{M \in \sM} \  \tT_{f_M} \ ]a,b[
\\[15pt]
&\leq \ 
\inf\limits_{M \in \sM_E} \  
\tT_{f_M} \ ]a,b[
\\[15pt]
&\leq \ 
\inf \ \{T_g ]a,b[ \ : \hsx g = f \ \text{almost everywhere}\}
\\[15pt]
&= \ 
e - \tT_f ]a,b[ \hsx.]
\end{align*}
\\[-1cm]
\end{x}

\begin{x}{\small\bf THEOREM} \ 
Let $f \in \BVL^1 ]a,b[$ $-$then there exists a $g \in \BV ]a,b[$ which is equal to $f$ almost everywhere and has the property that 
\[
\phi(f) 
\ = \ 
\tT_g ]a,b[  \hsx .
\]

PROOF \ Take $g$ admissible:
\[
\tT_g ]a,b[ 
\ \ = \ 
e - \tT_f ]a,b[
\ = \ 
\phi(f).
\]
\end{x}

\chapter{
$\boldsymbol{\S}$\textbf{20}.\quad  ABSOLUTE CONTINUITY III}
\setlength\parindent{2em}
\setcounter{theoremn}{0}
\renewcommand{\thepage}{\S20-\arabic{page}}

\begin{x}{\small\bf DEFINITION} \ 
A function $f : \ ]a,b[ \ \ra \R$ is said to be 
\un{absolutely} \un{continuous} in $]a,b[$ if for every $\varepsilon > 0$ there exists a $\delta > 0$ 
such that for any collection of non overlapping closed intervals
\[
[a_1, b_1] 
\subsetx 
]a,b[ , \ldots, [a_n, b_n] 
\subsetx 
]a,b[,
\]
then
\[
\sum\limits_{k = 1}^n \ 
(b_k - a_k) 
\ < \ 
\delta
\implies 
\sum\limits_{k = 1}^n \ 
\abs{f(b_k) -f(a_k)}
\ < \ 
\varepsilon.
\]
\\[-1cm]
\end{x}

\begin{x}{\small\bf NOTATION} \ 
$\AC \hsx ]a,b[$ is the set of absolutely continuous functions in $]a,b[$.
\\[-.5cm]
\end{x}

\begin{x}{\small\bf \un{N.B.}} \ 
An absolutely continuous function $f : \  ]a,b[ \ \ra \R$ is uniformly continuous.
\\[-.5cm]
\end{x}

\begin{x}{\small\bf RAPPEL} \ 
A uniformly continuous function $f : \  ]a,b[ \ \ra \R$ can be extended uniquely to $[a,b]$ in such a way that the 
extended function remains uniformly continuous. 
\\[-.5cm]
\end{x}

\begin{x}{\small\bf LEMMA} \ 
If $f \in \AC \hsx ]a,b[$, then its extension to $[a,b]$ belongs to $\AC [a,b]$.
\\[-.25cm]
\end{x}

\begin{x}{\small\bf THEOREM} \ 
Let 
$f : \ ]a,b[ \ \ra \R$  
$-$then $f$ is absolutely continuous iff the following four conditions are satisfied.
\\[-.5cm]

\qquad (1) \quad 
$f$ is continuous.
\\[-.5cm]

\qquad (2) \quad 
$f^\prime$ exists almost everywhere.
\\[-.5cm]

\qquad (3) \quad 
$f^\prime\in \Lp^p \hsx ]a,b[$ 
for some 
$1 \leq p < +\infty$.
\\[-.5cm]

\qquad (4) \quad 
$\forall \ x, \hsy x_0 \in \ ]a,b[$, 
\[
f(x) 
\ = \ 
f(x_0) + 
\int\limits_{x_0}^x \ 
f^\prime
\ \td \Lm^1.
\]
\\[-.5cm]
Here (and infra), $\Lp^1$ is Lebesgue measure on $]a,b[$.
\end{x}

\begin{x}{\small\bf \un{N.B.}} \ 
For the record, 
\[
\Lp^p \hsx ]a,b[ 
\subsetx 
\Lp^1 \hsx ]a,b[ 
\qquad 
(1 \leq p < +\infty).
\]
\\[-1.5cm]
\end{x}

\begin{x}{\small\bf DEFINITION} \ 
Let $1 \leq p < +\infty$ $-$then a function 
$f \in \Lp_\locx^1 \hsx ]a,b[$ admits a 
\un{weak derivative} in $\Lp^p \hsx ]a,b[$ if there exists a function 
$\ds\frac{\td f}{\td x} \in \Lp^p \hsx ]a,b[$ such that 
$\forall \ \phi \in C_c^\infty ]a,b[$, 
\[
\int\limits_{]a,b[} \ 
\phi \frac{\td f}{\td x}
\ \td \Lm^1 
\ = \ 
-
\int\limits_{]a,b[} \ 
\phi^\prime f
\ \td \Lm^1.
\]
\\[-1cm]
\end{x}

\begin{x}{\small\bf CRITERION} \ 
If $f \in \Lp_\locx^1 \hsx ]a,b[$ and if $\forall \ \phi \in C_c^\infty ]a,b[$,
\[
\int\limits_{]a,b[} \ 
\phi \hsy  f
\ \td \Lm^1 
\ = \ 
0,
\]
then $f = 0$ almost everywhere.
\\[-.25cm]
\end{x}

\begin{x}{\small\bf SCHOLIUM} \ 
A weak derivative of $f$ in $\Lp^p \hsy ]a,b[$, if it exists at all, is unique up to a set of Lebesgue measure 0.  
For suppose that you 
have two weak derivatives $u, \hsx v \in \Lp^p \hsy ]a,b[$, thus $\forall \ \phi \in C_c^\infty ]a,b[$,
\[
\begin{cases}
\ds
\int\limits_{]a,b[} \ 
\phi u
\ \td \Lm^1
\ = \ 
-
\int\limits_{]a,b[} \ 
\phi^\prime f
\ \td \Lm^1
\\[26pt]
\ds
\int\limits_{]a,b[} \ 
\phi v
\ \td \Lm^1
\ = \ 
-
\int\limits_{]a,b[} \ 
\phi^\prime f
\ \td \Lm^1
\end{cases}
\]
\qquad $\implies$
\[
\int\limits_{]a,b[} \ 
\phi (u - v)
\ \td \Lm^1
\ = \ 
0
\]
and so $u = v$ almost everywhere, $\phi \in C_c^\infty ]a,b[$ being arbitrary.
\\[-.25cm]
\end{x}

\begin{x}{\small\bf \un{N.B.}} \ 
If $f, \hsx g \in \Lp_\locx^1 \hsx ]a,b[$ are equal almost everywhere, then they have 
the ``same'' weak derivative. 
\\[-.5cm]

[$\forall \ \phi \in C_c^\infty ]a,b[$, 
\begin{align*}
\int\limits_{]a,b[} \ 
\phi \frac{\td f}{\td x}
\ \td \Lm^1 \ 
&=\ 
-
\int\limits_{]a,b[} \ 
\phi^\prime \hsy f
\ \td \Lm^1 \ 
\\[15pt]
&=\ 
-
\int\limits_{]a,b[} \ 
\phi^\prime \hsy g
\ \td \Lm^1 \ 
\\[15pt]
&=\ 
\int\limits_{]a,b[} \ 
\phi \frac{\td g}{\td x}
\ \td \Lm^1,
\end{align*}
so
\[
\frac{\td f}{\td x} \ = \  \frac{\td g}{\td x}
\]
almost everywhere.]
\\[-.25cm]
\end{x}

\begin{x}{\small\bf LEMMA} \ 
Let $f, \hsy g \in \Lp_\locx^1 \hsx ]a,b[$ and suppose that each of them admits a weak derivative $-$then 
$f + g$ admits a weak derivative and
\[
\frac{\td}{\td x} (f + g) 
\ = \ 
\frac{\td f}{\td x} + \frac{\td g}{\td x} \hsy.
\]

PROOF 
$\forall \ \phi \in C_c^\infty ]a,b[$,

\begin{align*}
\int\limits_{]a,b[} \ 
\phi 
\left(
\frac{\td f}{\td x} + \frac{\td g}{\td x}
\right)
\ \td \Lm^1 \ 
&=\ 
\int\limits_{]a,b[} \ 
\phi 
\frac{\td f}{\td x}
\ \td \Lm^1 \ 
\ + \ 
\int\limits_{]a,b[} \ 
\phi 
\frac{\td g}{\td x}
\ \td \Lm^1 \ 
\\[15pt]
&=\ 
-
\int\limits_{]a,b[} \
\phi^\prime \hsy f
\ \td \Lm^1 \ 
\ - \ 
\int\limits_{]a,b[} \
\phi^\prime \hsy g
\ \td \Lm^1 \ 
\\[15pt]
&=\ 
-
\int\limits_{]a,b[} \
\phi^\prime \hsy (f + g)
\ \td \Lm^1.
\end{align*}
\end{x}

\begin{x}{\small\bf LEMMA} \ 
If $\psi \in C_c^\infty ]a,b[$ and if $f$ admits a weak derivative $\ds\frac{\td f}{\td x} $, then
$\psi f$ admits a weak derivative and 
\[
\frac{\td}{\td x} (\psi f) 
\ = \ 
\psi^\prime f + \psi \frac{\td f}{\td x}.
\]

PROOF 
$\forall \ \phi \in C_c^\infty ]a,b[$,

\begin{align*}
\int\limits_{]a,b[} \ 
\phi^\prime (\psi f) 
\ \td \Lm^1 \ 
&=\ 
\int\limits_{]a,b[} \ 
\left(
f (\psi \phi)^\prime - f (\psi^\prime \phi)
\right)
\ \td \Lm^1
\\[15pt]
&=\ 
-
\int\limits_{]a,b[} \ 
\phi \hsy
\left(
\psi
\frac{\td f}{\td x} + f \psi^\prime
\right)
\ \td \Lm^1.
\end{align*}
\\[-1cm]
\end{x}

\begin{x}{\small\bf SUBLEMMA} \ 
Given $\phi \in C_c^\infty ]a,b[$, let
\[
\Phi(x) 
\ = \ 
\int\limits_{]a,x[} \ 
\phi
\ \td \Lm^1
\]
and suppose that 
\[
\int\limits_{]a,b[} \ 
\phi
\ \td \Lm^1
\ = \ 
0.
\]
Then $\Phi \in C_c^\infty ]a,b[$.
\\[-.25cm]
\end{x}

\begin{x}{\small\bf LEMMA} \ 
Let $f \in \Lp_\locx^1 \hsy ]a,b[$ and assume that $f$ has weak derivative 0 $-$then $f$ coincides 
almost everywhere in $]a,b[$ with a constant function.  
\\[-.25cm]

PROOF \ 
Fix 
$\psi_0 \in C_c^\infty ]a,b[$ : 
$\int\limits_{]a,b[} \ 
\psi_0
\ \td \Lm^1
\ = \ 1$, and given any 
$\phi \in C_c^\infty ]a,b[$, 
put 
$I(\phi) = 
\int\limits_{]a,b[} \ 
\phi
\ \td \Lm^1$
$-$then 
\[
I(\phi - I(\phi) \psi_0) 
\ = \ 
I(\phi) - I(\phi)I(\psi_0) 
\ = \ 
0,
\]
hence
\[
\Psi (x) 
\ = \ 
\int\limits_{]a,x[} \ 
(\phi - I(\phi) \psi_0) 
\ \td \Lm^1
\in C_c^\infty ]a,b[\hsy.
\]
Since $f$ has weak derivative 0, 
\[
\int\limits_{]a,b[} \ 
\Psi \frac{\td f}{\td x} 
\ \td \Lm^1
\ = \ 
0,
\]
\qquad $\implies$
\begin{align*}
0\ 
&=\ 
\int\limits_{]a,b[} \ 
\Psi^\prime f
\ \td \Lm^1
\\[15pt]
&=\ 
\int\limits_{]a,b[} \ 
\left(
\phi - I(\phi) \psi_0
\right) 
f
\ \td \Lm^1
\\[15pt]
&=\ 
\int\limits_{]a,b[} \ 
\phi 
\hsy f
\ \td \Lm^1
\ - \ 
\bigg(\ 
\int\limits_{]a,b[} \ 
\phi
\ \td \Lm^1
\bigg)
\bigg(\ 
\int\limits_{]a,b[} \ 
f \psi_0
\ \td \Lm^1
\bigg)
\\[15pt]
&=\ 
\int\limits_{]a,b[} \ 
\phi (f - C_0)
\ \td \Lm^1,
\end{align*}
where 
\[
C_0 
\ = \ 
\int\limits_{]a,b[} \ 
f \psi_0
\ \td \Lm^1.
\]
Therefore $f - C_0 = 0$ almost everywhere or still, $f = C_0$ almost everywhere.
\\[-.25cm]
\end{x}

\begin{x}{\small\bf NOTATION} \ 
Let $1 \leq p < +\infty$ $-$then 
$W^{1, p} ]a,b[$ is the set of all functions $f \in \Lp^p ]a,b[$ which possess a weak derivative 
$\ds\frac{\td f}{\td x}$ in $\Lp^p ]a,b[$.
\\[-.25cm]
\end{x}

\begin{x}{\small\bf \un{N.B.}} \ 
$W^{1, 1} ]a,b[$ is contained in $\BVL^1 ]a,b[$.  
\\[-.5cm]

[Take an $f \in W^{1, 1} ]a,b[$ and consider
\[
\tD f (E) 
\ = \ 
\int\limits_E \ 
\frac{\td f}{\td x} 
\ \td \Lm^1 
\qquad 
(E \in \text{BO\hsy $]a,b[$}\quad \text{(Borel sets in $]a,b[\hsy$)}),
\]
i.e., 
\[
\td \tD f 
\ = \ 
\frac{\td f}{\td x} \ \td \Lm^1.
\]
Then $\forall \ \phi \in C_c^\infty ]a,b[$, 
\begin{align*}
\int\limits_{]a,b[} \ 
\phi 
\ \td \tD f \ 
&=\ 
\int\limits_{]a,b[} \ 
\phi \frac{\td f}{\td x}
\ \td \Lm^1
\\[15pt]
&=\ 
- 
\int\limits_{]a,b[} \ 
\phi^\prime \hsy f
\ \td \Lm^1,
\end{align*}
so by definition, $f \in \BVL^1 \hsx ]a,b[\hsy.$]
\\[-.5cm]

[Note: \ 
The containment is strict.]
\\[-.25cm]
\end{x}

\begin{x}{\small\bf THEOREM} \ 
Let $1 \leq p < +\infty$ $-$then a function 
$f : \ ]a,b[ \ \ra \R$ 
belongs to $W^{1, p} ]a,b[$ iff it admits an absolutely continuous representative 
$\bar{f} : \ ]a,b[ \ \ra \R$ 
such that $\bar{f}$ and its ordinary derivative
$\bar{f}^\prime$ 
belongs to 
$\Lp^p \hsx ]a,b[$.
\\[-.25cm]
\end{x}

\begin{x}{\small\bf LEMMA} \ 
If $f \in \AC \hsx ]a,b[$, then $\forall \ \phi \in C_c^\infty ]a,b[$,
\[
\int\limits_{]a,b[} \ 
\phi \hsy f^\prime
\ \td \Lm^1 
\ = \ 
-
\int\limits_{]a,b[} \ 
\phi^\prime \hsy f
\ \td \Lm^1,
\]
there being no boundary term in the (implicit) integration by parts since $\phi$ has compact support in $]a,b[$.
\\[-.25cm]
\end{x}

\begin{x}{\small\bf SCHOLIUM} \ 
If $f$ is absolutely continuous, then its ordinary derivative $f^\prime$ is a weak derivative. 
\\[-.5cm]

One direction of the theorem is immediate.  
For suppose that 
$f : \ ]a,b[ \ \ra \R$
admits an absolutely continuous representative
$\bar{f} : \ ]a,b[ \ \ra \R$ 
such that 
$\bar{f}$ and $\bar{f}^\prime$
are in $\Lp^p \hsx ]a,b[$ $-$then the claim is that $f \in W^{1,p} ]a,b[$.  
Of course, $f \in \Lp^p \hsx ]a,b[$.  
As for the 
\\[-.5cm]

\noindent
existence of the weak derivative $\ds\frac{\td f}{\td x}$, note that $\forall \ \phi \in C_c^\infty ]a,b[$,
\[
\int\limits_{]a,b[} \ 
\phi \hsy \bar{f}^\prime
\ \td \Lm^1 
\ = \ 
-
\int\limits_{]a,b[} \ 
\phi^\prime \hsy \bar{f}
\ \td \Lm^1 
\]
or still, 
\[
\int\limits_{]a,b[} \ 
\phi \hsy \bar{f}^\prime
\ \td \Lm^1 
\ = \ 
-
\int\limits_{]a,b[} \ 
\phi^\prime \hsy f
\ \td \Lm^1,
\]
since $\bar{f} = f$ almost everywhere.  
Therefore $\bar{f}^\prime$ is a weak derivative of $f$ in $\Lp^p \hsy ]a,b[$.  
\\[-.25cm]

Turning to the converse, let $f \in W^{1,p} ]a,b[$, fix a point $x_0 \in \ ]a,b[$, and put
\[
\bar{f} (x) 
\ = \ 
f(x_0) + 
\int\limits_{x_0}^x \ 
\frac{\td f}{\td x}
\ \td \Lm^1 
\qquad (x \in ]a,b[).
\]
Then $\bar{f} \in \AC ]a,b[$ and almost everywhere, 
\[
\bar{f}^\prime
\ = \ 
\frac{\td f}{\td x}
\quad 
(\in \Lp^p \hsy ]a,b[)
\]
i.e., almost everywhere, 
\[
\bar{f}^\prime - \frac{\td f}{\td x}
\ = \ 
0,
\]
or still, almost everywhere, 
\[
\frac{\td }{\td x} (\bar{f} - f) 
\ = \ 
0,
\]
which implies that there exists a constant $C$ such that $\bar{f} - f = C$ almost everywhere, 
thus $f$ has an absolutely continuous representative $\bar{f}$ such that it and its ordinary derivative belong to 
$\Lp^p \hsy ]a,b[$.
\\[-.25cm]
\end{x}

\begin{x}{\small\bf REMARK} \ 
Matters simplify slightly when $p = 1$ : $f \in W^{1, 1} ]a,b[$ iff $f$ admists an absolutely continuous representative $\bar{f}$.
\end{x}

\newpage

\centerline{\textbf{\large REFERENCES}}
\setcounter{page}{1}
\setcounter{theoremn}{0}
\renewcommand{\thepage}{References-\arabic{page}}
\vspace{0.75cm}

\[
\text{BOOKS}
\]

\begin{rf}
Aumann, Georg, 
\textit{Reelle Funktionen}, Springer-Verlag Berlin, 1969. 
\end{rf}

\begin{rf}
Bogachev, V. I.,
\textit{Measure Theory, Volume I}, 
Springer-Verlag, Berlin Heidelberg, 2007.
\end{rf}

\begin{rf}
Bogachev, V. I.,
\textit{Measure Theory, Volume II}, 
Springer-Verlag, Berlin Heidelberg, 2007.
\end{rf}

\begin{rf}
Bruckner, Andrew M., Bruckner, Judith B., and Thomson, Brian S., 
\textit{Real Analysis}, 
Second Edition (2008). 
\end{rf}

\begin{rf}
Chandrasekharan, K., 
\textit{Classical Fourier Transforms}, 
Springer-Verlag, Berlin Heidelberg, 1989. 
\end{rf}

\begin{rf}
Courant, R., 
\textit{Differential and Integral Calculus Vol. I.}, 
Interscience Publishers, Inc., New York, 1937. 
\end{rf}

\begin{rf}
Courant, R., 
\textit{Differential and Integral Calculus Vol. II.}, 
Interscience Publishers, Inc., New York, 1939. 
\end{rf}

\begin{rf}
Dieudonn\'e, J., 
\textit{Foundations of Modern Analysis}, 
Academic Press, 
New York, 1960. 
\end{rf}

\begin{rf}
Flett, T. M., 
\textit{Mathematical Analysis}, 
McGraw-Hill House, England, 1966. 
\end{rf}

\begin{rf}
Hawkins, Thomas, 
\textit{Lebesgue's Theory of Integration}, 
Chelsea Publishing Company, New York, 1975. 
\end{rf}

\begin{rf}
Hewitt, Edwin and Stromberg, Karl, 
\textit{Real and Abstract Analysis}, 
Springer-Verlag, New York, 1969. 
\end{rf}

\begin{rf}
Kannan, R. and Krueger, Carole King, 
\textit{Advanced Analysis on the Real Line}, 
Springer-Verlag, New York, 1996. 
\end{rf}

\begin{rf}
Körner, T. W., 
\textit{Fourier Analysis}, 
Cambridge University Press, 1988. 
\end{rf}

\begin{rf}
Loomis, Lynn H., and Sternberg, Shlomo, 
\textit{Advanced Calculus}, 
Revised Edition, Jones and Bartlett Publishers, Boston, 1990.
\end{rf}

\begin{rf}
McShane, E. J., 
\textit{Unified Integration}, 
Academic Press, 1983. 
\end{rf}

\begin{rf}
Natanson, I. P.,
\textit{Theory of Functions of a Real Variable, Volume I.}, 
Ungar Publishing Co., New York, 1961. 
\end{rf}

\begin{rf}
Natanson, I. P., 
\textit{Theory of Functions of a Real Variable, Volume II}., 
Ungar Publishing Co., New York, 1960. 
\end{rf}

\begin{rf}
Rudin, Walter,
\textit{Principles of Mathematical Analysis}, 
Third Edition, 
McGraw-Hill Book Company, 1976. 
\end{rf}

\begin{rf}
Thomson, Brian S., 
\textit{Theory of the Integral}, 
Classical Real Analysis.com, 
2013. 
\end{rf}

\begin{rf}
Torchinsky, Alberto, 
\textit{Real-Variable Methods in Harmonic Analysis}, 
Academic Press, 1986. 
\end{rf}

\begin{rf}
Yeh, J., Real Analysis, 
\textit{Theory of Measure and Integration}, 
World Scientific Publishing Co. Pte. Ltd., Singapore, 2014. 
\end{rf}


\setcounter{theoremn}{0}

\[
\text{ARTICLES}
\]

\begin{rf}
Avdispahic, M., Concepts of generalized bounded variation and the theory of Fourier series, 
\textit{Internat. J. Math. Sci.} \textbf{9} (1986), 223-244.
\end{rf}

\begin{rf}
Bruckner, A. M., Density-preserving homeomorphisms and the theorem of Maximoff, 
\textit{Quart. J. Math. Oxford} Ser. (2) \textbf{21} (1970), 337-347.
\end{rf}

\begin{rf}
Buczolich, Z., Density points and bi-Lipschitz functions in $\R^m$, 
\textit{Proc. Amer. Math. Soc.} \textbf{116} (1992), 53-59.
\end{rf}

\begin{rf}
Fleming, W. H., Functions whose partial derivatives are measures, 
\textit{Illinois J. Math.} \textbf{4} (1960), 452-478.
\end{rf}

\begin{rf}
Goffman, C., A characterization of linearly continuous functions whose partial derivatives are measures, 
\textit{Acta Math.} \textbf{117} (1967), 165-190.
\end{rf}

\begin{rf}
Goffman, C. and Liu, F-C, Derivative measures, 
\textit{Proc. Amer. Math. Soc.} \textbf{78} (1980), 218-220.
\end{rf}

\begin{rf}
Goffman, C. and Neugebauer, C.J., On approximate derivatives, 
\textit{Proc. Amer. Math. Soc.} \textbf{11} (1960), 962-966.
\end{rf}

\begin{rf}
Hughs, R. E., Functions of BVC type, 
\textit{Proc. Amer. Math. Soc.} \textbf{12} (1961), 698-701.
\end{rf}

\begin{rf}
Serrin, J., On the differentiability of functions of several variables, 
\textit{Arch. Rational Mech. Anal.} \textbf{7} (1961), 359-372.
\end{rf}
 
\printindex
\end{document}